\documentclass[moor]{informs1}              % for a regular run
\usepackage{natbib}
 \NatBibNumeric
 \def\bibfont{\small}%

 \bibpunct[, ]{[}{]}{,}{n}{}{,}%

\usepackage[colorlinks=true,breaklinks=true,bookmarks=true,urlcolor=blue,
     citecolor=blue,linkcolor=blue,bookmarksopen=false,draft=false]{hyperref}

\def\EMAIL#1{\href{mailto:#1}{#1}}% When hyperref is used, otherwise outcomment
\TheoremsNumberedBySection  % (Theorem 1.1, Lema 1.1, Theorem 1.2)

\EquationsNumberedBySection % (1.1), (1.2), ...

\usepackage{bm}
\usepackage{enumitem}
\usepackage{mathtools}
\usepackage{algorithm}
\usepackage{algorithmic}
\usepackage{multirow}
\usepackage{float}
\usepackage{subfigure}
\def \beginproof{\par\noindent {\bf Proof.}\ }
\def \endproof{\hskip .5cm $\Box$ \vskip .5cm}
\definecolor{blue}{rgb}{0,0,0.9}
\definecolor{red}{rgb}{0.9,0,0}
\definecolor{green}{rgb}{0,0.9,0}

\graphicspath{{./images/}}
\begin{document}
%%%%%%%%%%%%%%%%

% Outcomment only when entries are known. Otherwise leave as is and
%   default values will be used.
%\setcounter{page}{1}
%\VOLUME{00}%
%\NO{0}%
%\MONTH{Xxxxx}% (month or a similar seasonal id)
%\YEAR{0000}% e.g., 2005
%\FIRSTPAGE{000}%
%\LASTPAGE{000}%
%\SHORTYEAR{00}% shortened year (two-digit)
%\ISSUE{0000} %
%\LONGFIRSTPAGE{0001} %
%\DOI{10.1287/xxxx.0000.0000}%

% Author's names for the running heads
% Sample depending on the number of authors;
% \RUNAUTHOR{Jones}
% \RUNAUTHOR{Jones and Wilson}
% \RUNAUTHOR{Jones, Miller, and Wilson}
% \RUNAUTHOR{Jones et al.} % for four or more authors
% Enter authors following the given pattern:
%\RUNAUTHOR{}

% Title or shortened title suitable for running heads. Sample:
% \RUNTITLE{Bundling Information Goods of Decreasing Value}
% Enter the (shortened) title:
%\RUNTITLE{}

% Full title. Sample:
\TITLE{Sparse Solutions of a Class of Constrained Optimization Problems}
% Enter the full title:
%\TITLE{}

% Block of authors and their affiliations starts here:
% NOTE: Authors with same affiliation, if the order of authors allows,
%   should be entered in ONE field, separated by a comma.
%   \EMAIL field can be repeated if more than one author
\ARTICLEAUTHORS{%
\AUTHOR{Lei Yang}
\AFF{ Department of Mathematics, National University of Singapore, 10 Lower Kent Ridge Road, Singapore 119076. (\EMAIL{yanglei.math@gmail.com})}
\AUTHOR{Xiaojun Chen}
\AFF{Department of Applied Mathematics, The Hong Kong Polytechnic University, Hung Hom, Kowloon, Hong Kong, China. (\EMAIL{xiaojun.chen@polyu.edu.hk})}
\AUTHOR{Shuhuang Xiang}
\AFF{School of Mathematics and Statistics, INP-LAMA, Central South University, Changsha, Hunan 410083, China. (\EMAIL{xiangsh@mail.csu.edu.cn})}
} % end of the block

\ABSTRACT{In this paper, we consider a well-known sparse optimization problem that aims to find a sparse solution of a possibly noisy underdetermined system of linear equations. Mathematically, it can be modeled in a unified manner by minimizing $\|\bm{x}\|_p^p$ subject to $\|A\bm{x}-\bm{b}\|_q\leq\sigma$ for given $A \in \mathbb{R}^{m \times n}$, $\bm{b}\in\mathbb{R}^m$, $\sigma \geq0$, $0\leq p\leq 1$ and $q \geq 1$. We then study various properties of the optimal solutions of this problem. Specifically, without any condition on the matrix $A$, we provide upper bounds in cardinality and infinity norm for the optimal solutions, and show that all optimal solutions must be on the boundary of the feasible set when $0<p\leq1$. Moreover, for $q \in \{1,\infty\}$, we show that the problem with $0<p<1$ has a finite number of optimal solutions and prove that there exists $0<p^*<1$ such that the solution set of the problem with any $0<p<p^*$ is contained in the solution set of the problem with $p=0$ and there further exists $0<\overline{p}<p^*$ such that the solution set of the problem with any $0<p\leq\overline{p}$ remains unchanged. An estimation of such $p^*$ is also provided. In addition, to solve the constrained nonconvex non-Lipschitz $L_p$-$L_1$ problem ($0<p<1$ and $q=1$), we propose a smoothing penalty method and show that, under some mild conditions, any cluster point of the sequence generated is a stationary point of our problem. Some numerical examples are given to implicitly illustrate the theoretical results and show the efficiency of the proposed algorithm for the constrained $L_p$-$L_1$ problem under different noises.
}%

% Fill in data. If unknown, outcomment the field
\MSCCLASS{Primary: 90C26, 90C30; secondary: 65K05}
\ORMSCLASS{Primary: mathematics, systems solution; secondary: programming, algorithms}
%\HISTORY{}

\maketitle
%%%%%%%%%%%%%%%%%%%%%%%%%%%%%%%%%%%%%%%%%%%%%%%%%%%%%%%%%%%%%%%%%%%%%%
\section{Introduction}

In this paper, we consider a class of sparse optimization problems, which can be modeled in a unified manner as the following constrained $L_p$-$L_q$ problem:
\begin{equation}\label{lplqmodel}
\min\limits_{\bm{x}\in\mathbb{R}^n}~~\|\bm{x}\|_p^p:={\textstyle\sum_{i=1}^{n}} |x_{i}|^p  \quad~~ \mbox{\rm s.t.} \quad~~ \|A\bm{x}-\bm{b}\|_q \leq \sigma,
\end{equation}
where $A \in \mathbb{R}^{m \times n}$, $\bm{b}\in\mathbb{R}^m$, $\sigma \geq0$, $0\le p\le 1$ and $1 \leq q \leq \infty$ are given. We assume that the feasible set of problem \eqref{lplqmodel} is nonempty so that problem \eqref{lplqmodel} is well-defined. With this assumption, one can easily verify that an optimal solution for $p=0$ (namely, a sparsest solution) exists thanks to the discrete and discontinuous nature of $\|\cdot\|_0$ and the closedness of the feasible set. Moreover, for $0<p\leq 1$, since $\|\bm{x}\|_p^p$ is level-bounded, then an optimal solution exists (see \cite[Theorem 1.9]{rw1998variational}). Therefore, the optimal solution set of problem \eqref{lplqmodel}, denoted by ${\rm SOL}(A,\bm{b},\sigma, p,q)$, is nonempty for any $0\le p\le 1$ and $1 \leq q \leq \infty$. We also assume that $\|\bm{b}\|_q > \sigma$ so that $A\neq0$ and $0 \not\in {\rm SOL}(A,\bm{b},\sigma,p,q)$. Obviously, when $p=1$, \eqref{lplqmodel} is a convex optimization problem and when $0<p<1$, \eqref{lplqmodel} yields a nonconvex and non-Lipschitz optimization problem.

Problem \eqref{lplqmodel} aims to find a sparse vector $\bm{x}$ from the corrupted observation $\bm{b}=A\bm{x}+\bm{\xi}$, where $\bm{\xi}$ denotes an unknown noisy vector bounded by $\sigma$ (the noise level) in $L_q$-norm, i.e., $\|\bm{\xi}\|_q\leq\sigma$. This problem arises in many contemporary applications and has been widely studied under different choices of $p$, $q$ and $\sigma$ in the literature; see, for example, \cite{bde2009sparse,cll2011constrained,CRT2006denoise,ct2005decoding,ct2007dantzig,c2007exact,clp2016penalty,cw2018spherical,cdd2009compressed,d2006compressed,det2005stable,dh2001uncertainty,fl2009sparsest,htw2015statistical,pyl2015np,sm2018,vf2008probing,zhang2013theory,zhao2013rsp,ZhaoYB}. Among these studies, the $L_2$-norm is commonly used for measuring the noise and leads to a mathematically tractable problem when the noise exists and comes from a Gaussian distribution \cite{bde2009sparse,CRT2006denoise,clp2016penalty,cw2018spherical,det2005stable,fl2009sparsest,sm2018,vf2008probing}.
In particular, it has been known that a sparse vector can be (approximately) recovered by the solution of the convex optimization problem \eqref{lplqmodel} with $p=1$ and $q=2$ under some well-known recovery conditions such as the restricted isometry property (RIP) \cite{CRT2006denoise}, the mutual coherence condition \cite{bde2009sparse,det2005stable} and the null space property (NSP) \cite{cdd2009compressed,zhang2013theory}. Such convex constrained $L_1$-$L_2$ problem can also be solved efficiently by a spectral projected gradient $L_1$ minimization algorithm (SPGL1) proposed by Van den Berg and Friedlander \cite{vf2008probing}. On the other hand, it is natural to find a sparse vector by solving problem \eqref{lplqmodel} with $0<p<1$ since $\|\bm{x}\|_p^p$ approaches $\|\bm{x}\|_0$ as $p\to0$. Indeed, under certain RIP conditions, Foucart and Lai \cite{fl2009sparsest} showed that a sparse vector can be (approximately) recovered by the solution of the nonconvex non-Lipschitz problem \eqref{lplqmodel} with $0<p<1$ and $q=2$. Chen, Lu and Pong \cite{clp2016penalty} also proposed a penalty method for solving this constrained $L_p$-$L_2$ problem ($0<p<1$) with promising numerical performances. Later, this penalty method and the SPGL1 are further combined to solve \eqref{lplqmodel} with $0<p<1$ and $q=2$ for recovering sparse signals on the sphere in \cite{cw2018spherical}. However, when the noise does not come from the Gaussian distribution but other heavy-tailed distributions (e.g., Student's t-distribution) \textit{or} contains outliers, using $\|A\bm{x}-\bm{b}\|_2$ as the data fitting term is no longer appropriate. In this case, some robust loss functions such as the $L_1$-norm \cite{FanLi2001,wlj2007robust,w2013l1} and the $L_{\infty}$-norm \cite{cll2011constrained,ct2007dantzig} are used to develop robust models. Recently, Zhao, Jiang and Luo \cite{ZhaoYB} also established a fairly comprehensive weak stability theory for problem \eqref{lplqmodel} with $p=1$ and $q\in \{1, 2, \infty\}$ under a so-called weak range space property (RSP) condition. The weak RSP condition can be induced by several existing compressed sensing matrix properties and hence can be the mildest one for the sparse solution recovery. However, it is still not easy to verify this condition in practice.

In this paper, we focus on problem \eqref{lplqmodel} with different choices of $p$ and $q$, and establish the following theoretical results concerning its optimal solutions \textit{without} any condition on the sensing matrix $A$.
\begin{itemize}[leftmargin=0.8cm]
\item[(i)] For any $\bm{x}^*\in{\rm SOL}(A,\bm{b},\sigma,p,q)$ with $0\leq p<1$ and $1 \leq q \leq \infty$, we have $\|\bm{x}^*\|_0 = \mathrm{rank}(A_{\mathcal J})$ and
    \begin{equation*}
    \frac{(\|\bm{b}\|_q-\sigma)\,m^{\min\left\{\frac{1}{2}-\frac{1}{q}, \,0\right\}}}{\sqrt{|\mathcal{J}|\lambda_{\max}(A_{\mathcal{J}}^{\top}A_{\mathcal{J}})}}
    \leq \|\bm{x}^*\|_{\infty}
    \leq \frac{\sigma m^{\max\left\{\frac{1}{2}-\frac{1}{q}, \,0\right\}}+\|\bm{b}\|_2}{\sqrt{\lambda_{\min}(A_{\mathcal{J}}^{\top}A_{\mathcal{J}})}},
    \end{equation*}
    where ${\mathcal J}=\mathrm{supp}(\bm{x}^*)$, $|\mathcal{J}|$ denotes its cardinality, and $\lambda_{\max}(A_{\mathcal{J}}^{\top}A_{\mathcal{J}})$ and $\lambda_{\min}(A_{\mathcal{J}}^{\top}A_{\mathcal{J}})$ are the largest and smallest eigenvalues of $A_{\mathcal{J}}^{\top}A_{\mathcal{J}}$, respectively. Moreover, for any $1 \leq q \leq \infty$, $\|A\bm{x}^*-\bm{b}\|_q=\sigma$ for $0<p\leq1$; and $\|A(\alpha\bm{x}^*)-\bm{b}\|_q=\sigma$ with some $\alpha\in (0,1]$ for $p=0$.

\item[(ii)] For $q\in\{1,\infty\}$, the solution set SOL$(A,\bm{b},\sigma,p,q)$ with $0<p<1$ has a finite number of elements.

\item[(iii)] There exists a $p^*\in (0,\,1]$ such that ${\rm SOL}(A,\bm{b},\sigma,p,1)\subseteq {\rm SOL}(A,\bm{b},\sigma,0,1)$ for any $p\in(0, \,p^*)$. An explicit estimation of such $p^*\in (0,1]$ is also given. Moreover, there exists a $\overline{p}\in (0,\,p^*)$ such that ${\rm SOL}(A,\bm{b},\sigma,p,1)= {\rm SOL}(A,\bm{b},\sigma,\overline{p},1)$ for any $p\in(0,\,\overline{p}]$.

\end{itemize}
Here, we would like to point out that the sparse solution recovery result (iii) is developed without any aforementioned recovery condition on $A$. This not only complements the existing recovery results in the literature, but also shows the potential advantage of using the $L_p$-norm ($0<p<1$) for recovering the sparse solution over the $L_1$-norm ball.

Note that problem \eqref{lplqmodel} is a constrained problem, while, in statistics and computer science, the $L_p$-$L_q$ problem/minimization often refers to the following unconstrained regularized problem \cite{cds2001atomic,cgwy2014complexity,cxy2010lower}:
\begin{equation}\label{BPDN}
\min\limits_{\bm{x}\in\mathbb{R}^n}~\|A\bm{x}-\bm{b}\|_q^q + \lambda \|\bm{x}\|^p_p,
\end{equation}
where $\lambda$ is a positive regularization parameter. Indeed, when $p=1$ and $q=2$, problem \eqref{BPDN} is the well-known $L_1$-regularized least-squares problem (namely, the LASSO problem) and it is known that, in this case, there exists a $\bar{\lambda}>0$ such that, for $\lambda\geq\bar{\lambda}$, the constrained problem \eqref{lplqmodel} is equivalent to the unconstrained problem \eqref{BPDN} regarding solutions; see, for example,
\cite[Section 3.2.3]{bde2009sparse}. However, Example 3.1 in \cite{clp2016penalty} shows that for $0<p<1$ and $q=2$, there does not exist a $\lambda$ so that problems \eqref{lplqmodel} and \eqref{BPDN} have a common global \textit{or} local minimizer. Hence, for $0<p<1$, one cannot expect to solve \eqref{lplqmodel} by solving the regularized problem \eqref{BPDN} with some fixed $\lambda>0$. In view of this, we shall consider a penalty method for solving problem \eqref{lplqmodel} with $0<p<1$, which basically solves the constrained problem \eqref{lplqmodel} by solving a sequence of unconstrained penalty problems. Specifically, we consider the following penalty problem of \eqref{lplqmodel}:
\begin{equation}\label{penaltylplq}
\min\limits_{\bm{x}\in\mathbb{R}^n}~~\|\bm{x}\|_p^p + \lambda \left(\|A\bm{x}-\bm{b}\|^q_q -\sigma^q\right)_+.
\end{equation}
Note that the function $\bm{x}\mapsto\|A\bm{x}-\bm{b}\|^q_q$ is continuously differentiable for $1<q<\infty$. Then, based on problem \eqref{penaltylplq}, one can readily extend the penalty method proposed in \cite{clp2016penalty} for solving problem \eqref{lplqmodel} with $0<p<1$ and $q=2$ to solve problem \eqref{lplqmodel} with $0<p<1$ and $1<q<\infty$. However, for $q\in\{1,\infty\}$, since the function $\bm{x} \mapsto \|A\bm{x}-\bm{b}\|_q$ is nonsmooth, then the approach in \cite{clp2016penalty} cannot be adapted directly. In view of this, in this paper, we propose an alternative smoothing penalty method for solving
\begin{equation}\label{lpl1model}
\min\limits_{\bm{x}\in\mathbb{R}^n}~~\|\bm{x}\|_p^p \quad~~ \mbox{\rm s.t.} \quad~~ \|A\bm{x}-\bm{b}\|_1 \leq \sigma,
\end{equation}
where $0<p<1$. Notice that we omit the case of $q=\infty$ to save space in this paper. Nevertheless, our approach can be extended without much difficulty to solve problem \eqref{lplqmodel} with $0<p<1$ and $q=\infty$, because the $L_{1}$-constrained problem and the $L_{\infty}$-constrained problem have similar properties in the sense that both constraints $\|A\bm{x}-\bm{b}\|_1 \leq \sigma$ and $\|A\bm{x}-\bm{b}\|_{\infty} \leq \sigma$ can be represented as linear constraints, and the functions $\bm{x}\mapsto\|A\bm{x}-\bm{b}\|_1$ and $\bm{x}\mapsto\|A\bm{x}-\bm{b}\|_\infty$ are piecewise linear. We shall show that problem \eqref{penaltylplq} with $q=1$ is the exact penalty problem of problem \eqref{lpl1model} regarding local minimizers and global minimizers. We also prove that any cluster point of a sequence generated by our smoothing penalty method is a stationary point of problem \eqref{lpl1model}. Moreover, some numerical results are reported to show that all computed stationary points have the properties in our theoretical contribution (i) mentioned above. Here, we would like to emphasize that finding a global optimal solution of \eqref{lpl1model} is NP-hard \cite{cgwy2014complexity,gjy2011note}. Thus, it is interesting to see that our smoothing penalty method can efficiently find a `good' stationary point of problem \eqref{lpl1model}, which has important properties of a global optimal solution of problem \eqref{lpl1model}.

The rest of this paper is organized as follows. In Section \ref{secpro}, we rigorously prove properties (i)-(iii) listed above and give a concrete example to verify these properties. In Section \ref{secpen}, we present a smoothing penalty method for solving problem \eqref{lpl1model} and show some convergence results. Some numerical results are presented in Section \ref{secnum}, with some concluding remarks given in Section \ref{secconc}.

\paragraph{\textbf{Notation and Preliminaries}}
In this paper, we use the convention that $\frac{1}{\infty}=0$. For an index set $\mathcal{J}\subseteq\{1,\cdots,n\}$, let $|\mathcal{J}|$ denote its cardinality and $\mathcal{J}^c$ denote its complementarity set. We denote by $\bm{x}_{\mathcal{J}}\in\mathbb{R}^{|\mathcal{J}|}$ the subvector formed from a vector $\bm{x}\in\mathbb{R}^n$ by picking the entries indexed by $\mathcal{J}$ and denote by $A_{\mathcal{J}}\in\mathbb{R}^{m\times |\mathcal{J}|}$ the submatrix formed from a matrix $A\in\mathbb{R}^{m\times n}$ by picking the columns indexed by $\mathcal{J}$. Recall from \cite[Definition 8.3]{rw1998variational} that, for a proper closed function $f$, the regular (or Fr\'{e}chet) subdifferential, the (limiting) subdifferential and the horizon subdifferential of $f$ at $\bm{x}\in{\rm dom}\,f$ are defined respectively as
\begin{equation*}
\begin{aligned}
\widehat{\partial} f(\bm{x}) &:= \left\{\bm{d} \in \mathbb{R}^{n} : \liminf\limits_{\bm{y} \rightarrow \bm{x}, \,\bm{y} \neq \bm{x}}\, \frac{f(\bm{y})-f(\bm{x})-\langle \bm{d}, \bm{y}-\bm{x}\rangle}{\|\bm{y}-\bm{x}\|} \geq 0\right\}, \\
\partial f(\bm{x}) &:= \left\{ \bm{d} \in \mathbb{R}^{n}: \exists \,\bm{x}^k \xrightarrow{f} \bm{x}, ~\bm{d}^k \rightarrow \bm{d} ~~\mathrm{with}~~\bm{d}^k\in\widehat{\partial} f(\bm{x}^k) \right\}, \\
\partial^{\infty} f(\bm{x}) &:= \left\{ \bm{d} \in \mathbb{R}^{n}: \exists \,\bm{x}^k \xrightarrow{f} \bm{x}, ~\lambda_k\bm{d}^k \rightarrow \bm{d},~\lambda_k\downarrow0 ~~\mathrm{with}~~\bm{d}^k\in\widehat{\partial} f(\bm{x}^k) \right\}.
\end{aligned}
\end{equation*}
It can be observed from the above definitions (or see \cite[Proposition 8.7]{rw1998variational}) that
\begin{equation}\label{robust}
\begin{aligned}
\left\{ \bm{d}\in\mathbb{R}^{n}: \exists \,\bm{x}^k \xrightarrow{f} \bm{x}, ~\bm{d}^k \rightarrow \bm{d} ~\mathrm{with}~\bm{d}^k \in \partial f(\bm{x}^k) \right\} &\subseteq \partial f(\bm{x}), \\
\left\{ \bm{d}\in\mathbb{R}^{n}: \exists \,\bm{x}^k \xrightarrow{f} \bm{x}, ~\lambda_k\bm{d}^k \rightarrow \bm{d},~\lambda_k\downarrow0 ~\mathrm{with}~\bm{d}^k \in \partial f(\bm{x}^k) \right\} &\subseteq \partial^{\infty} f(\bm{x}).
\end{aligned}
\end{equation}
When $f$ is convex, the above (limiting) subdifferential coincides with the classical subdifferential in convex analysis \cite[Proposition~8.12]{rw1998variational}. Moreover, if $f$ is continuously differentiable, we have $\partial f(\bm{x})=\{\nabla f(\bm{x})\}$, where $\nabla f(\bm{x})$ is the gradient of $f$ at $\bm{x}$ \cite[Exercise~8.8(b)]{rw1998variational}. For a closed set $\mathcal{X}\subseteq\mathbb{R}^{n}$, its indicator function $\delta_{\mathcal{X}}$ is defined by $\delta_{\mathcal{X}}(\bm{x})=0$ if $\bm{x}\in\mathcal{X}$ and $\delta_{\mathcal{X}}(\bm{x})=+\infty$ otherwise. In addition, we use $\mathcal{B}(\bm{y};\delta)$ to denote the closed ball of radius $\delta$ centered at $\bm{y}$, i.e., $\mathcal{B}(\bm{y};\delta):=\{\bm{x}\in\mathbb{R}^n: \|\bm{x} - \bm{y}\|_2 \leq \delta\}$, and ${\rm FEA}(A,\bm{b},\sigma,q):=\left\{\bm{x}\in \mathbb{R}^n : \|A\bm{x} - \bm{b}\|_q \leq \sigma\right\}$ to denote the feasible set of problem \eqref{lplqmodel}.

%%%%%%%%%%%%%%%%%%%%%%%%%%%%%%%%%%%%%%%%%%%%%%%%%%%%%%
\section{Properties of solutions of problem \eqref{lplqmodel}}\label{secpro}

In this section, we characterize the properties of the optimal solutions of problem \eqref{lplqmodel} with different choices of $p$ and $q$. We first give a supporting lemma.

\begin{lemma}\label{lem-lqineq}
Let $1 \leq q \leq \infty$. For any $\bm{x}\in\mathbb{R}^n$, we have
\begin{equation*}
n^{\min\left\{\frac{1}{q}-\frac{1}{2}, \,0\right\}}\|\bm{x}\|_2 \leq \|\bm{x}\|_q \leq n^{\max\left\{\frac{1}{q}-\frac{1}{2}, \,0\right\}}\|\bm{x}\|_2.
\end{equation*}
\end{lemma}
\beginproof
We consider the following two cases.
\begin{itemize}
\item $1\leq q \leq 2$. In this case, it is easy to see that $\|\bm{x}\|_q \geq \|\bm{x}\|_2$. On the other hand, since $2/q\geq1$, it then follows from the H\"{o}lder's inequality that
    \begin{equation*}
    \|\bm{x}\|^q_q = \sum_{i=1}^{n} |x_i|^q = \sum_{i=1}^{n} |x_i|^q\cdot1 \leq \left(\sum_{i=1}^{n} \left(|x_i|^q\right)^{\frac{2}{q}}\right)^{\frac{q}{2}} \left(\sum_{i=1}^{n} 1^{\frac{2}{2-q}}\right)^{1-\frac{q}{2}} = n^{1-\frac{q}{2}} \|\bm{x}\|_2^q,
    \end{equation*}
    which results in $\|\bm{x}\|_q \leq n^{\frac{1}{q}-\frac{1}{2}} \|\bm{x}\|_2$.

\item $q\geq2$. In this case, it is easy to see that $\|\bm{x}\|_q \leq \|\bm{x}\|_2$. On the other hand, since $q/2\geq1$, it then follows from the H\"{o}lder's inequality that
    \begin{equation*}
    \|\bm{x}\|^2_2 = \sum_{i=1}^{n} |x_i|^2 = \sum_{i=1}^{n} |x_i|^2\cdot1 \leq \left(\sum_{i=1}^{n} \left(|x_i|^2\right)^{\frac{q}{2}}\right)^{\frac{2}{q}} \left(\sum_{i=1}^{n} 1^{\frac{q}{q-2}}\right)^{1-\frac{2}{q}} = n^{1-\frac{2}{q}} \|\bm{x}\|_q^2,
    \end{equation*}
    which results in $\|\bm{x}\|_q \geq n^{\frac{1}{q}-\frac{1}{2}} \|\bm{x}\|_2$.
\end{itemize}
Combing the above results, we prove this lemma.
\endproof

The following theorem is given for $0 \leq p \leq 1$ and $1 \leq q \leq \infty$.

\begin{theorem}\label{thmeqvl0}
Let $1 \leq q \leq \infty$. For any $\bm{x}^* \in{\rm SOL}(A,\bm{b},\sigma,p,q)$, the following statements hold with ${\mathcal J}:=\mathrm{supp}(\bm{x}^*)$.
\begin{itemize}[leftmargin=0.8cm]
\item[{\rm (i)}] For $0<p\leq1$, $\|A\bm{x}^*-\bm{b}\|_q=\sigma$; and for $p=0$, there is a scalar $\alpha\in (0,1]$ such that $\|A(\alpha\bm{x}^*)-\bm{b}\|_q=\sigma$ and $\alpha' \bm{x}^* \in{\rm SOL}(A,\bm{b},\sigma,0,q)$ for any $\alpha'\in[\alpha,\,1]$.

\item[{\rm (ii)}] For $0\leq p<1$, $\|\bm{x}^*\|_0=|\mathcal J|=\mathrm{rank}(A_{\mathcal J})$.

\item[{\rm (iii)}] For $0\leq p<1$,
         \begin{equation*}
         \frac{(\|\bm{b}\|_q-\sigma)\,m^{\min\left\{\frac{1}{2}-\frac{1}{q}, \,0\right\}}}{\sqrt{|\mathcal{J}|\lambda_{\max}(A_{\mathcal{J}}^{\top}A_{\mathcal{J}})}}
         \leq \|\bm{x}^*\|_{\infty}
         \leq \frac{\sigma m^{\max\left\{\frac{1}{2}-\frac{1}{q}, \,0\right\}}+\|\bm{b}\|_2}{\sqrt{\lambda_{\min}(A_{\mathcal{J}}^{\top}A_{\mathcal{J}})}},
         \end{equation*}
         where $\lambda_{\max}(A_{\mathcal{J}}^{\top}A_{\mathcal{J}})$ and $\lambda_{\min}(A_{\mathcal{J}}^{\top}A_{\mathcal{J}})$ are the largest and smallest eigenvalues of $A_{\mathcal{J}}^{\top}A_{\mathcal{J}}$, respectively. Moreover, when $\sigma=0$, we have $\bm{x}_{\mathcal{J}}^*=(A_{\mathcal{J}}^{\top}A_{\mathcal{J}})^{-1}A_{\mathcal{J}}^{\top}\bm{b}$.
\end{itemize}
\end{theorem}
\beginproof
\textit{Statement (i)}. If $\|A\bm{x}^*-\bm{b}\|_q=\sigma$, the results hold trivially. Next, we assume that $\|A\bm{x}^*-\bm{b}\|_q<\sigma$.

Consider $0<p\leq1$. From $\|\bm{b}\|_q > \sigma$, we see that $A\bm{x}^* \neq 0$. Then, it is easy to verify that there exists a constant $0<c<1$ such that $\|A(c\bm{x}^*)-\bm{b}\|_q < \sigma$. Thus, $c\bm{x}^* \in {\rm FEA}(A,\bm{b},\sigma,q)$, but $\|c\bm{x}^*\|^p_p = c^p\|\bm{x}^*\|_p^p < \|\bm{x}^*\|_p^p$ for $0<p\leq1$. This leads to a contradiction. Hence, we have $\|A\bm{x}^*-\bm{b}\|_q=\sigma$.

Consider $p=0$. Let $f(t):=\|A(t\bm{x}^*)-\bm{b}\|_q$. Then, from the continuity of $f$, $f(0)=\|\bm{b}\|_q>\sigma$ and $f(1)=\|A\bm{x}^*-\bm{b}\|_q<\sigma$, there exists a scalar $\alpha\in (0,1)$ such that $f(\alpha)=\|A(\alpha\bm{x}^*)-\bm{b}\|_q=\sigma$. Moreover, it is easy to verify that $f$ is convex on $[0,\,1]$. Thus, for any $\alpha'\in[\alpha,\,1]$, there exists a $0\leq\lambda\leq1$ such that $\alpha'=\lambda\alpha+(1-\lambda)$ and $f(\alpha')\leq\lambda f(\alpha) + (1-\lambda)f(1) \leq \sigma$. Hence, $\alpha'\bm{x}^*$ is feasible. This together with $\|\alpha'\bm{x}^*\|_0=\|\bm{x}^*\|_0$ shows that $\alpha' \bm{x}^* \in{\rm SOL}(A,\bm{b},\sigma,0,q)$ for any $\alpha'\in[\alpha,\,1]$.

\textit{Statement (ii)}. Let $s:=\|\bm{x}^*\|_0=|\mathcal{J}|$ for simplicity. We then consider the following two cases.

{\bf Case 1}, $p=0$. First, it is not hard to see that $s \leq m$ since any set of $m+1$ vectors in $\mathbb{R}^m$ is linearly dependent. Thus, we have $\mathrm{rank}(A_{\mathcal{J}}) \leq \min\{m,\,s\}=s$. We next prove $\mathrm{rank}(A_{\mathcal{J}})=s$ by contradiction. Assume that $\mathrm{rank}(A_{\mathcal{J}}) < s$. Then, there exists a vector $\hat{\bm{h}}\in\mathbb{R}^{s}$ such that $\hat{\bm{h}}\neq0$ and $A_{\mathcal{J}}\hat{\bm{h}} = 0$. Let $\bm{h}\in\mathbb{R}^n$ be a vector such that $\bm{h}_{\mathcal{J}}=\hat{\bm{h}}$ and $\bm{h}_{\mathcal{J}^c}=0$. Thus, we have $A\bm{h}=0$. Now, let
\begin{equation*}
\tau := \min\limits_{h_i\neq0,\,i\in\mathcal{J}} \left\{\frac{x^*_i}{h_i}\right\}=\frac{x^*_{i_0}}{h_{i_0}}~~\mathrm{for}~\mathrm{some}~i_0.
\end{equation*}
Then, we see that $\tilde{\bm{x}}:=\bm{x}^* - \tau\bm{h} \in {\rm FEA}(A,\bm{b},\sigma,q)$ since $A\tilde{\bm{x}}=A(\bm{x}^* - \tau\bm{h})=A\bm{x}^*$. Moreover, from the definition of $\tau$, one can verify that $\tilde{x}_{i_0}=0$ and thus $\|\tilde{\bm{x}}\|_0 < \|\bm{x}^*\|_0$. This leads to a contradiction. Hence, we only have $\mathrm{rank}(A_{\mathcal{J}}) = s=\|\bm{x}^*\|_0$.

{\bf Case 2}, $0<p<1$. We first prove $s \leq m$ by contradiction. Assume that $s > m$. Thus, there exists a vector $\tilde{\bm{h}}\in\mathbb{R}^{s}$ such that $\tilde{\bm{h}}\neq0$ and $A_{\mathcal{J}}\tilde{\bm{h}} = 0$, since $\mathrm{rank}(A_{\mathcal{J}}) \leq \min \{m, \, s\} = m < s$. Let $\bm{h} \in \mathbb{R}^n$ be a vector such that $\bm{h}_{\mathcal{J}}=\tilde{\bm{h}}$ and $\bm{h}_{\mathcal{J}^c}=0$. Thus, we have that $A\bm{h}=0$ and hence $\bm{x}^* + t\bm{h} \in {\rm FEA}(A,\bm{b},\sigma,q)$ for any $t\in\mathbb{R}$. Moreover, we can choose a sufficiently small real positive number $t_0 > 0$ such that, for all $|t|\leq t_0$,
\begin{equation}\label{keepsgn}
\bm{x}^*_{\mathcal{J}} + t\bm{h}_{\mathcal{J}} \neq 0, \quad \mathrm{and} \quad \mathrm{sgn}(x^*_i) = \mathrm{sgn}(x^*_i + th_i) \quad \mathrm{for} \quad i \in \mathcal{J}.
\end{equation}
Let $f(t):=\sum_{i\in\mathcal{J}}\left[\mathrm{sgn}(x_i^*)(x^*_i + th_i)\right]^p$. Then, we have
\begin{equation*}%\label{add-eqt0}
\begin{aligned}
f(0)
&=\sum_{i\in\mathcal{J}}|x_i^*|^p = \|\bm{x}^*\|_p^p = \min\limits_{t \in[-t_0,\,t_0]} \|\bm{x}^* + t\bm{h}\|_p^p \\
&=\min\limits_{t \in[-t_0,\,t_0]} \,\sum\limits_{i\in\mathcal{J}}\left[\mathrm{sgn}(x_i^*+ th_i)(x^*_i + th_i)\right]^p = \min\limits_{t \in[-t_0,\,t_0]} f(t),
\end{aligned}
\end{equation*}
where the third equality follows because $\bm{x}^* \in{\rm SOL}(A,\bm{b},\sigma,p,q)$ and the last equality follows from \eqref{keepsgn}. However, for all $|t|\leq t_0$,
\begin{equation*}
f''(t) = p(p-1){\textstyle\sum_{i\in\mathcal{J}}}\left[\mathrm{sgn}(x_i^*)(x^*_i + th_i)\right]^{p-2}h_i^2 < 0.
\end{equation*}
This leads to a contradiction. Hence, we have $s \leq m$ and $\mathrm{rank}(A_{\mathcal{J}}) \leq \min \{m, \, s\} = s$. We further assume that $\mathrm{rank}(A_{\mathcal{J}}) < s$. Then, there also exists a vector $\hat{\bm{h}}\in\mathbb{R}^{s}$ such that $\hat{\bm{h}}\neq0$ and $A_{\mathcal{J}}\hat{\bm{h}} = 0$. Using the similar arguments as above, we can get a contradiction. Hence, we only have that $\mathrm{rank}(A_{\mathcal{J}})=s$.

\textit{Statement (iii)}. From statement (ii), $A_{\mathcal{J}}$ has full column rank and hence $\lambda_{\min}(A_{\mathcal{J}}^{\top}A_{\mathcal{J}})\neq0$. Then, we see that
\begin{equation*}
\begin{aligned}
\sigma
&\geq \|A\bm{x}^*-\bm{b}\|_q = \|A_{\mathcal{J}}\bm{x}_{\mathcal{J}}^*-\bm{b}\|_q
\geq m^{\min\left\{\frac{1}{q}-\frac{1}{2}, \,0\right\}} \|A_{\mathcal{J}}\bm{x}_{\mathcal{J}}^*-\bm{b}\|_2 \\
&\geq m^{\min\left\{\frac{1}{q}-\frac{1}{2}, \,0\right\}}(\|A_{\mathcal{J}}\bm{x}_{\mathcal{J}}^*\|_2 - \|\bm{b}\|_2)
\geq m^{\min\left\{\frac{1}{q}-\frac{1}{2}, \,0\right\}}\Big{(}\sqrt{\lambda_{\min}(A_{\mathcal{J}}^{\top}A_{\mathcal{J}})}\|\bm{x}_{\mathcal{J}}^*\|_2 - \|\bm{b}\|_2\Big{)},
\end{aligned}
\end{equation*}
where the second inequality follows from Lemma \ref{lem-lqineq} and the last inequality follows from $\|A_{\mathcal{J}}\bm{x}_{\mathcal{J}}^*\|_2^2 \geq \lambda_{\min}(A_{\mathcal{J}}^{\top}A_{\mathcal{J}})\|\bm{x}_{\mathcal{J}}^*\|_2^2$. Thus, the above relation implies that
\begin{equation*}
\|\bm{x}^*\|_{\infty} \leq \|\bm{x}^*\|_2 = \|\bm{x}_{\mathcal{J}}^*\|_2 \leq \frac{\sigma m^{\max\left\{\frac{1}{2}-\frac{1}{q}, \,0\right\}}+\|\bm{b}\|_2}{\sqrt{\lambda_{\min}(A_{\mathcal{J}}^{\top}A_{\mathcal{J}})}},
\end{equation*}
which gives the upper bound for $\|\bm{x}^*\|_{\infty}$. On the other hand, we have
\begin{equation*}
\begin{aligned}
\sigma
&\geq \|A\bm{x}^*-\bm{b}\|_q = \|A_{\mathcal{J}}\bm{x}_{\mathcal{J}}^*-\bm{b}\|_q
\geq \|\bm{b}\|_q - \|A_{\mathcal{J}}\bm{x}_{\mathcal{J}}^*\|_q   \\
&\geq \|\bm{b}\|_q - m^{\max\left\{\frac{1}{q}-\frac{1}{2}, \,0\right\}}\|A_{\mathcal{J}}\bm{x}_{\mathcal{J}}^*\|_2
\geq \|\bm{b}\|_q - m^{\max\left\{\frac{1}{q}-\frac{1}{2}, \,0\right\}}\sqrt{\lambda_{\max}(A_{\mathcal{J}}^{\top}A_{\mathcal{J}})}\|\bm{x}_{\mathcal{J}}^*\|_2,
\end{aligned}
\end{equation*}
where the third inequality follows from Lemma \ref{lem-lqineq} and the last inequality follows from $\|A_{\mathcal{J}}\bm{x}_{\mathcal{J}}^*\|_2^2 \leq \lambda_{\max}(A_{\mathcal{J}}^{\top}A_{\mathcal{J}})\|\bm{x}_{\mathcal{J}}^*\|_2^2$. This results in
\begin{equation*}
\|\bm{x}^*\|_{\infty} = \|\bm{x}^*_{\mathcal{J}}\|_{\infty} \geq \frac{\|\bm{x}_{\mathcal{J}}^*\|_2}{\sqrt{|\mathcal{J}|}} \geq \frac{(\|\bm{b}\|_q-\sigma)\,m^{\min\left\{\frac{1}{2}-\frac{1}{q}, \,0\right\}}}{\sqrt{|\mathcal{J}|\lambda_{\max}(A_{\mathcal{J}}^{\top}A_{\mathcal{J}})}},
\end{equation*}
which gives the lower bound for $\|\bm{x}^*\|_{\infty}$. Recall that $\|\bm{b}\|_q > \sigma$ (our blanket assumption). Thus, this lower bound is nontrivial. Moreover, when $\sigma=0$, we have $A\bm{x}^*=A_{\mathcal{J}}\bm{x}_{\mathcal{J}}^*=\bm{b}$ and hence $\bm{x}_{\mathcal{J}}^*=(A_{\mathcal{J}}^{\top}A_{\mathcal{J}})^{-1}A_{\mathcal{J}}^{\top}\bm{b}$. We then complete the proof.
\endproof

\begin{remark}[\textbf{The sparse solution of the $L_p$-$L_2$ problem}]
Theorem \ref{thmeqvl0}(ii) implies that without any condition on the sensing matrix $A$, $\|\bm{x}^*\|_0 \leq \min(m,n)$ for any $\bm{x}^* \in{\rm SOL}(A,\bm{b},\sigma,p,q)$ with $0<p<1$ and $1 \leq q \leq \infty$, while Shen and Mousavi show in \cite[Proposition 3.1]{sm2018} that $\|\bm{x}^*\|_0 \geq n-m+1$ for any $\bm{x}^* \in {\rm SOL}(A,\bm{b},\sigma,p,2)$ with $p>1$ and $n \geq m$ if every $m \times m$ submatrix of $A$ is invertible. Combining these results gives a formal confirmation that if $n\gg m$, all solutions of the $L_p$-$L_2$ problem with $0\leq p<1$ are sparse, but the $L_p$-$L_2$ problem with $p>1$ may not have sparse solutions.
\end{remark}

In the following, we shall derive more theoretical results for the optimal solution set of the $L_1$-constrained problem \eqref{lpl1model} with $0 \leq p < 1$. But we should point out that all results established later can be extended without much difficulty to the $L_{\infty}$-constrained case or other more general cases; see Remarks \ref{addremark1} and \ref{addremark2} for more details. As we shall see later, solving problem \eqref{lpl1model} with an arbitrarily sufficiently small $0<p<1$ actually gives an optimal solution of problem \eqref{lpl1model} with $p=0$. This nice result is obtained based on a simple observation that the feasible set ${\rm FEA}(A,\bm{b},\sigma,1)$ is indeed a convex polyhedron in $\mathbb{R}^n$ (see Lemma \ref{polyhd}). Moreover, observe that $\mathbb{R}^n$ can be represented as a union of $2^n$ orthants, denoted by $\mathbb{P}_j$ for $j=1,\cdots,2^n$, such that any two vectors $\bm{x}$ and $\bm{y}$ in each $\mathbb{P}_j$ have the same sign for each entry, i.e., for each $\mathbb{P}_j$, we have
\begin{equation}\label{defPorthant}
\forall\,\bm{x}, \,\bm{y} \in \mathbb{P}_j  \quad \Longrightarrow \quad x_iy_i \geq 0~~\mathrm{for}~~i = 1, \cdots, n.
\end{equation}
For example, when $n=2$, we have $\mathbb{R}^2=\bigcup^4_{j=1}\mathbb{P}_j$, where $\mathbb{P}_1=\{\bm{x}:x_1\geq0, x_2\geq0\}$, $\mathbb{P}_2=\{\bm{x}:x_1\geq0, x_2\leq0\}$, $\mathbb{P}_3=\{\bm{x}:x_1\leq0, x_2\geq0\}$ and $\mathbb{P}_4=\{\bm{x}:x_1\leq0, x_2\leq0\}$. Then, for each $j$, one can see that $\mathbb{P}_{j} \cap {\rm FEA}(A,\bm{b},\sigma,1)$ is empty or a polyhedron that has a finite number of extreme points because $\mathbb{P}_{j} \cap {\rm FEA}(A,\bm{b},\sigma,1)$ contains no lines; see
\cite[Corollary 18.5.3]{r1970convex} and \cite[Corollary 19.1.1]{r1970convex}.

\begin{lemma}\label{exoctant}
Let $0<p<1$. Suppose that $j \in \{1,\cdots,2^n\}$ is an arbitrary index such that $\mathbb{P}_{j}\cap {\rm FEA}(A,\bm{b},\sigma,1) \neq \emptyset$, where $\mathbb{P}_j$ is defined in \eqref{defPorthant}. Then, any optimal solution of the following problem
\begin{equation}\label{LPP}
\min\limits_{\bm{x}\in\mathbb{R}^n}\,\,\|\bm{x}\|_p^p \qquad \mbox{\rm s.t.} \qquad \bm{x} \in \mathbb{P}_{j}\cap {\rm FEA}(A,\bm{b},\sigma,1)
\end{equation}
is an extreme point of $\mathbb{P}_{j}\cap{\rm FEA}(A,\bm{b},\sigma,1)$.
\end{lemma}
\beginproof
Let $\bm{x}^*$ be an optimal solution of \eqref{LPP}. Suppose that there exist $\bm{y},\,\bm{z} \in \mathbb{P}_{j}\cap {\rm FEA}(A,\bm{b},\sigma,1)$ such that $\bm{x}^*=\lambda \bm{y} + (1-\lambda)\bm{z}$ for some $0< \lambda< 1$. Then, we have
\begin{equation*}
\begin{aligned}
\|\bm{x}^*\|_p^p
&=\|\lambda \bm{y} + (1-\lambda)\bm{z}\|_p^p
={\textstyle\sum_{i=1}^n}\left|\lambda y_i+(1-\lambda)z_i\right|^p
= {\textstyle\sum_{i=1}^n}\left(\lambda |y_i|+(1-\lambda)|z_i|\right)^p \\
&\geq {\textstyle\sum_{j=1}^n} \left(\lambda|y_j|^p + (1-\lambda)|z_j|^p\right)
= \lambda\|\bm{y}\|_p^p+(1-\lambda)\|\bm{z}\|_p^p
\geq \|\bm{x}^*\|_p^p,
\end{aligned}
\end{equation*}
where the third equality follows because any $\bm{y}, \,\bm{z} \in \mathbb{P}_{j}$ have the same sign for each entry, the first inequality follows because $f(t)=t^p$ is strictly concave for $t \geq 0$, and the last inequality follows because $\bm{x}^*$ is an optimal solution of \eqref{LPP}. Note that the above relation holds if and only if $\bm{y}=\bm{z}=\bm{x}^*$. This implies that $\bm{x}^*$ is an extreme point of $\mathbb{P}_{j}\cap {\rm FEA}(A,\bm{b},\sigma,1)$.
\endproof

Based on Lemma \ref{exoctant}, we are able to characterize the number of the optimal solutions of problem \eqref{lpl1model} with $0<p<1$. For notational simplicity, for $j=1,\cdots,2^n$, let
\begin{equation*}
{\rm EXT}\left(\mathbb{P}_{j}\cap{\rm FEA}(A,\bm{b},\sigma,1)\right) := \big{\{}\mathrm{all~extreme~points~of~} \mathbb{P}_{j}\cap{\rm FEA}(A,\bm{b},\sigma,1)\big{\}}.
\end{equation*}

\begin{proposition}\label{extrp}
For any $0<p<1$, the optimal solution set ${\rm SOL}(A,\bm{b},\sigma,p,1)$ of problem \eqref{lpl1model} is a finite set. Moreover, the set $\bigcup_{0<p<1}{\rm SOL}(A,\bm{b},\sigma,p,1)$ is a finite set.
\end{proposition}
\beginproof
For a given $0<p<1$, let $\bm{x}^*$ be an optimal solution of problem \eqref{lpl1model}, i.e., $\bm{x}^*\in{\rm SOL}(A,\bm{b},\sigma,p,1)$. Then, there must exist a $j^* \in \{1,\cdots,2^n\}$ such that $\bm{x}^* \in \mathbb{P}_{j^*}\cap{\rm FEA}(A,\bm{b},\sigma,1)$ and $\bm{x}^*$ is also an optimal solution of \eqref{LPP} with $j^*$ in place of $j$. Then, it follows from Lemma \ref{exoctant} that $\bm{x}^*$ is an extreme point of $\mathbb{P}_{j^*}\cap{\rm FEA}(A,\bm{b},\sigma,1)$. This implies that
\begin{equation}\label{boundbyex}
{\rm SOL}(A,\bm{b},\sigma,p,1) \subseteq \bigcup_{j \in \{1,\cdots,2^n\}} {\rm EXT}\left(\mathbb{P}_{j}\cap{\rm FEA}(A,\bm{b},\sigma,1)\right).
\end{equation}
Note that, for each $j$, $\mathbb{P}_{j} \cap {\rm FEA}(A,\bm{b},\sigma,1)$ is empty or a polyhedron that has a finite number of extreme points since $\mathbb{P}_{j} \cap {\rm FEA}(A,\bm{b},\sigma,1)$ contains no lines; see \cite[Corollary 18.5.3]{r1970convex} and \cite[Corollary 19.1.1]{r1970convex}. This together with \eqref{boundbyex} implies that ${\rm SOL}(A,\bm{b},\sigma,p,1)$ is a finite set.

Moreover, since \eqref{boundbyex} holds for any $0<p<1$, then we have
\begin{equation*}
\bigcup_{0<p<1}{\rm SOL}(A,\bm{b},\sigma,p,1) \subseteq \bigcup_{j \in \{1,\cdots,2^n\}}{\rm EXT}\left(\mathbb{P}_{j}\cap{\rm FEA}(A,\bm{b},\sigma,1)\right),
\end{equation*}
which implies $\bigcup_{0<p<1}{\rm SOL}(A,\bm{b},\sigma,p,1)$ is a finite set. This completes the proof.
\endproof

\begin{remark}[\textbf{Comments on Proposition \ref{extrp}}]\label{addremark1}
Proposition \ref{extrp} is obtained based on the observation that the feasible set ${\rm FEA}(A,\bm{b},\sigma,1)$ is  a convex polyhedron in $\mathbb{R}^n$. From this observation, we can extend Proposition \ref{extrp} to that for any $0<p<1$, the optimal solution set ${\rm SOL}(A,\bm{b},\sigma,p,q)$ of \eqref{lplqmodel} with $q=\infty$ is a finite set. However, it is not clear whether for any $0<p<1$, the optimal solution set ${\rm SOL}(A,\bm{b},\sigma,p,q)$ of problem \eqref{lplqmodel} with $q=2$ is a finite set. Thanks to Theorem \ref{thmeqvl0}, we can claim that if $A$ satisfies $\mathrm{rank}(A)=2$, the optimal solution set ${\rm SOL}(A,\bm{b},\sigma, \frac{1}{k}, 2)$ is a finite set, where $k\ge 2$ is a positive integer. Indeed, in this case, by Theorem \ref{thmeqvl0}(ii), any optimal solution $\bm{x}^*$ satisfies that $\|\bm{x}^*\|_0=|\mathcal{J}|=\mathrm{rank}(A_{\mathcal J})\leq\mathrm{rank}(A)=2$ and hence has at most two nonzero entries supported on $\mathcal{J}$. Then, there are only $\frac{n(n-1)}{2}$ different choices of the support set $\mathcal{J}$. Let $\nu^*$ be the optimal objective value and, without loss of generality, assume that $x^*_1\ge 0$, $x^*_2\ge 0$, $x^*_3=\cdots= x^*_n=0$. Then, $\sqrt[k]{x^*_1} + \sqrt[k]{x^*_2}=\nu^*$. Also, let $t:= \sqrt[k]{x^*_1}$ and $\sqrt[k]{x^*_2}=\nu^*-t$. We then see from Theorem \ref{thmeqvl0}(i) that $\|A_{\mathcal{J}}\bm{x}^*_{\mathcal{J}}-\bm{b}\|^2_2=\|A\bm{x}^*-\bm{b}\|^2_2=\sigma^2$ and this equation can be further written as a $2k$-th order polynomial equation $f(t)=0$, which has at most $2k$ real roots. This implies that, for each $\mathcal{J}$ satisfying $|\mathcal{J}|=2$, there are only $2k$ different choices of $x^*_1$ and $x^*_2$. Hence, the optimal solution set ${\rm SOL}(A,\bm{b},\sigma, \frac{1}{k}, 2)$ is a finite set and the number of solutions is at most $n(n-1)k$.
\end{remark}

We next give two supporting lemmas and relegate the proofs to Appendices \ref{prooflemmas1} and \ref{prooflemmas2}, respectively.

\begin{lemma}\label{inequ}
Suppose that $\bm{a}=(a_1,\cdots,a_n)^{\top}\in\mathbb{R}^n$ and $\bm{b}=(b_1,\cdots,b_n)^{\top}\in\mathbb{R}^n$ satisfy
\begin{equation*}
a_1 \leq a_2 \leq \cdots \leq a_n, \quad
b_1 \leq b_2 \leq \cdots \leq b_n, \quad
{\textstyle\sum_{j=1}^n} a_j^k = {\textstyle\sum_{j=1}^n} b_j^k, \quad
k = 1,\cdots,n,
\end{equation*}
then $\bm{a} = \bm{b}$.
\end{lemma}

\begin{lemma}\label{addlemma}
Given $\bm{a}$, $\bm{b}\in\mathbb{R}^n$ with $\|\bm{a}\|_0=\|\bm{b}\|_0=s$. Let $\{a_{i_1}, \cdots, a_{i_s}\}$ and $\{b_{t_1},\cdots,b_{t_s}\}$ be the nonzero entries in $\bm{a}$ and $\bm{b}$, respectively, and, without loss of generality, assume that $|a_{i_1}|\leq\cdots\leq|a_{i_s}|$ and $|b_{t_1}|\leq\cdots\leq|b_{t_s}|$. For $k = 1, \cdots, s$, define
\begin{equation}\label{defDeltak}
\Delta_k(\bm{a}, \bm{b}) := {\textstyle\sum_{j=1}^{s}}\left((\ln |a_{i_j}|)^{k}-(\ln |b_{t_j}|)^{k}\right).
\end{equation}
Then, the following statements hold.
\begin{itemize}[leftmargin=0.8cm]
\item[{\rm (i)}] If $\Delta_k(\bm{a}, \bm{b})=0$ for all $k=1,\cdots,s$, then $\|\bm{a}\|_p^p=\|\bm{b}\|_p^p$ holds for any $p>0$.
\item[{\rm (ii)}] Otherwise, there exists a sufficiently small $p'$ such that either $\|\bm{a}\|_p^p<\|\bm{b}\|_p^p$ or $\|\bm{a}\|_p^p>\|\bm{b}\|_p^p$ holds for any $p\in(0,\,p']$.
\end{itemize}
\end{lemma}

Now, we are ready to present our results concerning the optimal solution set ${\rm SOL}(A,\bm{b},\sigma,p,1)$ with different choices of $p$.

\begin{theorem}\label{mianthm1}
There exists a $p^*\in (0,\,1]$ such that ${\rm SOL}(A,\bm{b},\sigma,p,1)\subseteq {\rm SOL}(A,\bm{b},$ $\sigma,0,1)$ for any $p\in(0, \,p^*)$. Moreover, there exists a $\overline{p}\in (0,\,p^*)$ such that ${\rm SOL}(A,\bm{b},$ $\sigma,p,1)= {\rm SOL}(A,\bm{b},\sigma,\overline{p},1)$ for any $p\in(0,\,\overline{p}]$.
\end{theorem}
\beginproof
We prove the first result by contradiction. Assume that there does \textit{not} exist a number $p^*\in (0,\,1]$ such that, for any $p\in(0, \,p^*)$, ${\rm SOL}(A,\bm{b},\sigma,p,1)\subseteq {\rm SOL}(A,\bm{b},\sigma,0,1)$. Consider a sequence $\{p_k\}$ with $0<p_k<1$ and $p_k \to 0$ as $k \to \infty$. Thus, from the hypothesis, for each $p_k$, there exists a point $\bm{x}^{k}$ such that $\bm{x}^k \in {\rm SOL}(A,\bm{b},\sigma,p_k,1)$ and $\bm{x}^k \notin {\rm SOL}(A,\bm{b},\sigma,0,1)$. Now, we consider the sequence $\{\bm{x}^{k}\}$. Note that all elements in $\{\bm{x}^{k}\}$ come from the set $\bigcup_{0<p<1}{\rm SOL}(A,\bm{b},\sigma,p,1)$ but they are \textit{not} contained in ${\rm SOL}(A,\bm{b},\sigma,0,1)$. Since there are only finitely many points in $\bigcup_{0<p<1}{\rm SOL}(A,\bm{b},\sigma,p,1)$ (by Proposition \ref{extrp}), then there exists at least one point $\hat{\bm{x}} \in \bigcup_{0<p<1}{\rm SOL}(A,\bm{b},\sigma,p,1)$ such that $\{\bm{x}^{k}\}$ contains \textit{infinitely} many $\hat{\bm{x}}$, i.e., there exists a subsequence $\{\bm{x}^{k_j}\}$ so that $\bm{x}^{k_j} \equiv \hat{\bm{x}}$ for all $k_j$. Moreover, let $\bm{x}^*\in{\rm SOL}(A,\bm{b},\sigma,0,1)$. Then, for all $k_j$, we have $\|\bm{x}^{k_j}\|_{p_{k_j}}^{p_{k_j}} \leq \|\bm{x}^*\|_{p_{k_j}}^{p_{k_j}}$ since $\bm{x}^{k_j} \in {\rm SOL}(A,\bm{b},\sigma,p_{k_j},1)$. Then, we see that
\begin{equation*}
\|\bm{x}^*\|_0=\lim\limits_{k_j \to \infty}\|\bm{x}^*\|_{p_{k_j}}^{p_{k_j}} \geq \lim\limits_{k_j \to \infty} \|\bm{x}^{k_j}\|_{p_{k_j}}^{p_{k_j}} = \lim\limits_{k_j \to \infty} \|\hat{\bm{x}}\|_{p_{k_j}}^{p_{k_j}} = \|\hat{\bm{x}}\|_0,
\end{equation*}
which implies that $\hat{\bm{x}}\in{\rm SOL}(A,\bm{b},\sigma,0,1)$. This leads to a contradiction and completes the proof for the first result.

Next, we prove the second result. For notational simplicity, let $\mathcal{S}_{0\sim p^*}:=\bigcup_{0<p<p^*}{\rm SOL}(A,\bm{b},\sigma,p,1)$ and $s:=\|\bm{x}^*\|_0$, where $\bm{x}^*\in{\rm SOL}(A,\bm{b},\sigma,0,1)$. For any $\bm{x}\in\mathcal{S}_{0\sim p^*}$, we have $\|\bm{x}\|_0=s$ (by the first result) and define a set as $\mathcal{C}(\bm{x}):=\big\{\bm{z}\in\mathcal{S}_{0\sim p^*} : \Delta_k(\bm{x}, \bm{z})=0,\, \forall\,k=1,\cdots,s\big\}$, where $\Delta_k(\cdot,\cdot)$ is defined as \eqref{defDeltak}. Then, given $\bm{x}\in\mathcal{S}_{0\sim p^*}$ and $\bm{y}\in\mathcal{S}_{0\sim p^*}\setminus\mathcal{C}(\bm{x})$, it follows from Lemma \ref{addlemma}(ii) that there exists a sufficiently small $p^{(x,y)}\in(0, \,p^*)$ such that either $\|\bm{x}\|_p^p<\|\bm{y}\|_p^p$ or $\|\bm{x}\|_p^p>\|\bm{y}\|_p^p$ holds for any $p\in(0,\,p^{(x,y)}]$. Since $\mathcal{S}_{0\sim p^*}$ is contained in $\bigcup_{0<p<1}{\rm SOL}(A,\bm{b},\sigma,p,1)$, then the number of such a pair $(\bm{x},\,\bm{y})$ is finite. Therefore, we must have a sufficiently small $\tilde{p}\in(0, \,p^*)$ such that, for any $\bm{x}\in\mathcal{S}_{0\sim p^*}$ and $\bm{y}\in\mathcal{S}_{0\sim p^*}\setminus\mathcal{C}(\bm{x})$, either $\|\bm{x}\|_p^p<\|\bm{y}\|_p^p$ or $\|\bm{x}\|_p^p>\|\bm{y}\|_p^p$ holds for any $p\in(0,\,\tilde{p}]$. Now, for such $\tilde{p}$, consider any $p'\in(0,\,\tilde{p}]$ and let $\bm{x}'\in {\rm SOL}(A,\bm{b},\sigma,p',1)$. We must have $\|\bm{x}'\|_{p'}^{p'}<\|\bm{y}\|_{p'}^{p'}$ for any $\bm{y} \in \mathcal{S}_{0\sim p^*}\setminus\mathcal{C}(\bm{x}')$. This together with Lemma \ref{addlemma}(ii) implies that for any $0<p<p'$, $\|\bm{x}'\|_{p}^{p}<\|\bm{y}\|_{p}^{p}$ for any $\bm{y} \in \mathcal{S}_{0\sim p^*}\setminus\mathcal{C}(\bm{x}')$. Moreover, from Lemma \ref{addlemma}(i), for any $p>0$, $\|\bm{x}'\|_{p}^{p}=\|\bm{y}\|_{p}^{p}$ for any $\bm{y} \in \mathcal{C}(\bm{x}')$. These two facts show that for any $0<p<p'$, $\|\bm{x}'\|_{p}^{p}\leq\|\bm{y}\|_{p}^{p}$ for any $\bm{y} \in \mathcal{S}_{0\sim p^*}$. Hence, we have $\bm{x}'\in{\rm SOL}(A,\bm{b},\sigma,p,1)$ for any $0<p<p'\leq\tilde{p}$. Since $p'$ is arbitrary and $\bm{x}'\in{\rm SOL}(A,\bm{b},\sigma,p',1)$ is also arbitrary, we can conclude that ${\rm SOL}(A,\bm{b},\sigma,p',1)\subseteq{\rm SOL}(A,\bm{b},\sigma,p'',1)$ for any $0<p''<p'\leq\tilde{p}$.

We now prove by contradiction that there must exist a $\overline{p}\in (0,\,\tilde{p}]$ such that ${\rm SOL}(A,\bm{b},\sigma,p,1)= {\rm SOL}(A,\bm{b},\sigma,\overline{p},1)$ for any $p\in(0,\,\overline{p}]$. Assume this is not true. Then, for any $p'\in (0,\,\tilde{p}]$, there exists a $p''\in(0,\,p')$ such that ${\rm SOL}(A,\bm{b},\sigma,p'',1)\neq{\rm SOL}(A,\bm{b},\sigma,p',1)$. This together with the conclusion obtained above implies that ${\rm SOL}(A,\bm{b},\sigma,p',1)$ must be \textit{strictly} contained in ${\rm SOL}(A,\bm{b},\sigma,p'',1)$, i.e., ${\rm SOL}(A,\bm{b},$ $\sigma, p',1)\subset{\rm SOL}(A,\bm{b},\sigma,p'',1)$. With this fact, we generate a sequence $\{p^k\}$ as follows. Let $p^0=\tilde{p}$. Then, there exists a $p_1\in(0,\,p_0)$ such that ${\rm SOL}(A,\bm{b},\sigma,p_0,1)\subset{\rm SOL}(A,\bm{b},\sigma,p_1,1)$. For such $p_1$, there exists a $p_2\in(0,\,p_1)$ such that ${\rm SOL}(A,\bm{b},\sigma,p_1,1)$ $\subset{\rm SOL}(A,\bm{b},\sigma,p_2,1)$. Repeating this procedure, we can obtain a sequence $\{p^k\}$ such that $p_0>p_1>\cdots>0$ and ${\rm SOL}(A,\bm{b},\sigma,p_0,1)\subset{\rm SOL}(A,\bm{b},\sigma,p_1,1)\subset\cdots\subset{\rm SOL}(A,\bm{b},\sigma,0,1)$. Thus, along such sequence $\{p^k\}$, the number of elements of ${\rm SOL}(A,\bm{b},\sigma,p,1)$ will \textit{strictly} increase and hence $\bigcup_{\{p^k\}}{\rm SOL}(A,\bm{b},\sigma,p,1)$ must have infinitely many elements. This leads to a contradiction and completes the proof.
\endproof

\begin{remark}[\textbf{Comments on Theorem \ref{mianthm1}}]\label{addremark2}
Theorem \ref{mianthm1} is established based on the observation that the feasible set ${\rm FEA}(A,\bm{b},\sigma,1)$ of problem \eqref{lpl1model} is a polyhedron, and then, for each $j$, $\mathbb{P}_{j}\cap{\rm FEA}(A,\bm{b},\sigma,1)$ has at most a finite number of extreme points. Thus, one can also consider minimizing $\|\bm{x}\|_0$ under many other polyhedral constraints, for example, $\left\{\bm{x}\in \mathbb{R}^n :\|A\bm{x}-\bm{b}\|_{\infty}\leq\sigma\right\}$ and $\{\bm{x}\in \mathbb{R}^n : \|A\bm{x}-\bm{b}\|_1\leq\sigma, \,\bm{l}\leq\bm{x}\leq\bm{u}\}$ with $\bm{l}\in\mathbb{R}^n\cup\{-\infty\}^n$, $\bm{u}\in\mathbb{R}^n\cup\{\infty\}^n$ and $\bm{l}<\bm{u}$, to fit different scenarios in practice. Following the similar arguments presented in this paper, one can obtain the similar results in Theorem \ref{mianthm1} as well as Theorem \ref{pstarthm} under these polyhedral constraints. Moreover, it is also possible to extend our smoothing penalty method presented in the next section to solve problems in these cases. Here, we will omit more details to avoid overcomplicating the presentation. In addition, we are aware that the first result in Theorem \ref{mianthm1} has also been discussed in \cite{yhw2019sparsest}. However, the analysis there is much more tedious.
\end{remark}

Based on Theorem \ref{mianthm1}, it is easy to give the following corollary for $\sigma=0$ (namely, the noiseless case), which has also been discussed in \cite[Theorem 1]{pyl2015np}.

\begin{corollary}
There exists a $p^* \in (0,\,1]$ such that, for any $p\in(0,\,p^*)$, every optimal solution of problem $\min\big{\{}\|\bm{x}\|_p^p : A\bm{x} = \bm{b}\big{\}}$ is an optimal solution of problem $\min\big{\{}\|\bm{x}\|_0 : A\bm{x} = \bm{b}\big{\}}$.
\end{corollary}

Theorem \ref{mianthm1} says that there exists a $p^* \in (0, \,1]$ such that solving problem \eqref{lpl1model} with any $p\in(0, \,p^*)$ also solves problem \eqref{lpl1model} with $p=0$. Therefore, the constant $p^*$ is obviously the key for such nice relation and we are interested in estimating such $p^*$ in the next theorem. Our analysis is motivated by that of \cite[Theorem 1]{pyl2015np}, but makes use of results developed in Theorem \ref{thmeqvl0} and Lemma \ref{exoctant} for the more general feasible set. Before proceeding, we define two constants as follows:
\begin{eqnarray}
r &:=& \frac{\sigma+\|\bm{b}\|_2}{\sqrt{\lambda^*}}, ~\mathrm{where ~\lambda^* ~is ~the ~smallest ~nonzero ~eigenvalue ~of~} A^{\top}A, \label{defr1}  \\ [2pt]
\tilde{r} &:=& \min\left\{\,|x_i|~:~
\bm{x}\in{\textstyle\bigcup_{j \in \{1,\cdots,2^n\}}} {\rm EXT}\left(\mathbb{P}_{j}\cap{\rm FEA}(A,\bm{b},\sigma,1)\right), ~|x_i|\neq0, ~1 \leq i \leq n\right\}.  \label{defrt}
\end{eqnarray}
Note that for any subset $\mathcal{I}\subseteq \{1,\cdots,n\}$ such that $A_{\mathcal{I}}$ has full column rank, $A_{\mathcal{I}}^{\top}A_{\mathcal{I}}$ is a principal submatrix of $A^{\top}A$. Then, it follows from \cite[Theorem 1.4.10]{hj2012matrix} that $\lambda_{\min}(A_{\mathcal{I}}^{\top}A_{\mathcal{I}})>0$ is an eigenvalue of $A^{\top}A$ and hence $\lambda_{\min}(A_{\mathcal{I}}^{\top}A_{\mathcal{I}})\geq\lambda^*$. This together with Theorem \ref{thmeqvl0}(iii) implies that
\begin{equation}\label{ineqadd3}
\|\bm{x}^*\|_\infty \leq r, ~\mathrm{for~any~} \bm{x}^*\in{\rm SOL}(A,\bm{b},\sigma,p,1)~\mathrm{with}~0\leq p<1.
\end{equation}
From \eqref{defrt}, \eqref{ineqadd3} and Lemma \ref{exoctant}, one can also see that $r\geq\tilde{r}$.

\begin{theorem}\label{pstarthm}
Let $s$ be the optimal objective value of problem \eqref{lpl1model} with $p=0$ and
\begin{equation}\label{defpstar}
p^*:=\min \left\{1, ~\frac{\ln(s+1)-\ln s}{\ln r - \ln \tilde{r}} \right\}.
\end{equation}
Then, for any $p\in(0, \,p^*)$, ${\rm SOL}(A,\bm{b},\sigma,p,1)\subseteq {\rm SOL}(A,\bm{b},\sigma,0,1)$.
\end{theorem}
\beginproof
First, we show that
\begin{equation}\label{ineqadd1}
\left(\frac{r}{\tilde{r}}\right)^{p}s < s + 1
\end{equation}
holds for any $p\in(0, \,p^*)$. Since $\frac{r}{\tilde{r}} \geq 1$ and \eqref{ineqadd1} holds trivially when $\frac{r}{\tilde{r}} = 1$, then we only consider $\frac{r}{\tilde{r}} > 1$ in the following two cases.
\begin{itemize}[leftmargin=0.8cm]
\item $\frac{r}{\tilde{r}} \leq \frac{s+1}{s}$. In this case, $p^*=1$. Since $\frac{r}{\tilde{r}} > 1$, then $\left(\frac{r}{\tilde{r}}\right)^p < \frac{r}{\tilde{r}} \leq \frac{s+1}{s}$ for any $p\in(0, \,1)$.

\item $\frac{r}{\tilde{r}} > \frac{s+1}{s}$. In this case, $p^*=\frac{\ln(s+1)-\ln s}{\ln r-\ln\tilde{r}}$. Since $\frac{r}{\tilde{r}} > 1$, then $\left(\frac{r}{\tilde{r}}\right)^p < \left(\frac{r}{\tilde{r}}\right)^{p^*} = \frac{s+1}{s}$ for any $p\in(0, \,p^*)$.
\end{itemize}
Hence, \eqref{ineqadd1} holds for any $p\in(0, \,p^*)$.

Next, let $\bm{x}^*$ be an arbitrary optimal solution of problem \eqref{lpl1model} with $p\in(0, \,p^*)$, i.e., $\bm{x}^*\in{\rm SOL}(A,\bm{b},\sigma,p,1)$. It then follows from Lemma \ref{exoctant} that $\bm{x}^*$ is an extreme point of $\mathbb{P}_{j^*}\cap{\rm FEA}(A,\bm{b},\sigma,1)$ for some $j^* \in \{1,\cdots,2^n\}$. Thus, we have $\frac{|x^*_i|}{\tilde{r}}\geq 1$ for any $|x^*_i|\neq0$. Moreover, we see that
\begin{equation*}
\begin{aligned}
s
&\leq \|\bm{x}^*\|_0
= \lim\limits_{p'\downarrow0}\,\sum^n_{i=1}|x_i^*|^{p'}
= \lim\limits_{p'\downarrow0}\,\sum^n_{i=1}\left(\frac{|x_i^*|}{\tilde{r}}\right)^{p'}
\leq \sum^n_{i=1}\left(\frac{|x_i^*|}{\tilde{r}}\right)^{p}
= \tilde{r}^{-p} \,\|\bm{x}^*\|^p_p  \\
&= \tilde{r}^{-p} \min\left\{\|\bm{x}\|_p^p: \|A\bm{x} - \bm{b}\|_1 \leq \sigma\right\}
\stackrel{{\rm(i)}}{=} \tilde{r}^{-p} \min\left\{\|\bm{x}\|_p^p: \|A\bm{x} - \bm{b}\|_1 \leq \sigma, \,\|\bm{x}\|_{\infty} \leq r\right\} \\
&= \left(\frac{r}{\tilde{r}}\right)^{p} \min\left\{\|r^{-1}\bm{x}\|_p^p: \|A\bm{x} - \bm{b}\|_1 \leq \sigma, \,\|\bm{x}\|_{\infty} \leq r\right\} \\
&\leq \left(\frac{r}{\tilde{r}}\right)^{p} \min\left\{\|r^{-1}\bm{x}\|_0 : \|A\bm{x} - \bm{b}\|_1 \leq \sigma, \,\|\bm{x}\|_{\infty} \leq r\right\}   \\
&= \left(\frac{r}{\tilde{r}}\right)^{p} \min\left\{\|\bm{x}\|_0 : \|A\bm{x} - \bm{b}\|_1 \leq \sigma, \,\|\bm{x}\|_{\infty} \leq r\right\} \\
&\stackrel{{\rm(ii)}}{=} \left(\frac{r}{\tilde{r}}\right)^{p} \min\left\{\|\bm{x}\|_0 : \|A\bm{x} - b\|_1 \leq \sigma\right\}
= \left(\frac{r}{\tilde{r}}\right)^{p}s  < s + 1,
\end{aligned}
\end{equation*}
where the second inequality follows because for any $t\geq1$, the function $p \mapsto t^p$ is non-decreasing on $[0, 1)$, the equality (i) follows from \eqref{ineqadd3}, the third inequality follows because for any $0 \leq t \leq 1$, the function $p \mapsto t^p$ is non-increasing on $[0, 1)$, the equality (ii) follows again from \eqref{ineqadd3}, and the last inequality follows from \eqref{ineqadd1}. Then, from the above relation, we have that $\|\bm{x}^*\|_0=s$ and hence $\bm{x}^*$ is an optimal solution of problem \eqref{lpl1model} with $p=0$. This completes the proof.
\endproof

Before closing this section, we present a simple example to illustrate our previous theoretical results.

\begin{example}
Let $A = \begin{bmatrix*}[r] 1 & 1 & 1 \\ 1 & 1 & -1\end{bmatrix*}$, $\bm{b} = \begin{bmatrix} 3 \\ 3 \end{bmatrix}$ and $\sigma = 1$. Then, we consider
\begin{equation}\label{exam1}
\min\limits_{\bm{x}\in\mathbb{R}^3}~\|\bm{x}\|_p^p \quad~~ \mbox{\rm s.t.} \quad~~ \|A\bm{x} - \bm{b}\|_q \leq 1
\end{equation}
with $0 \leq p \leq 1$ and $q=1,\,2,\,\infty$. Next, for each $q$, we discuss the optimal solution sets of problem \eqref{exam1} with different choices of $p$.

\textbf{For} $q=1$, the feasible set of \eqref{exam1} is
\begin{equation*}
\mathcal{S}_1:=\big{\{}\bm{x}\in\mathbb{R}^3:|x_1+x_2+x_3-3|+|x_1+x_2-x_3-3|\leq1\big{\}}.
\end{equation*}
Then,
\begin{equation}\label{egresults}
\hspace{6mm}\begin{aligned}
{\rm Arg}\min\limits_{\bm{x}\in\mathcal{S}_1}\{\|\bm{x}\|_0\} &= {\textstyle\left\{(x_1, 0, 0)^{\top} : \frac{5}{2} \leq x_1 \leq \frac{7}{2}\right\} \cup \left\{(0, x_2, 0)^{\top}: \frac{5}{2} \leq x_2 \leq \frac{7}{2}\right\}},  \\
{\rm Arg}\min\limits_{\bm{x}\in\mathcal{S}_1}\{\|\bm{x}\|_p^p\} &= {\textstyle\left\{(\frac{5}{2}, \,0, \,0)^{\top}, \,(0, \,\frac{5}{2}, \,0)^{\top}\right\}} \quad \mathrm{for}~\mathrm{any}~0<p<1,  \\
{\rm Arg}\min\limits_{\bm{x}\in\mathcal{S}_1}\{\|\bm{x}\|_1\} &= {\textstyle\left\{(x_1, x_2, 0)^{\top}: x_1 + x_2 = \frac{5}{2}, \,x_1\geq 0, \,x_2\geq 0 \right\}}.
\end{aligned}
\end{equation}

\textbf{For} $q=2$, the feasible set of \eqref{exam1} is
\begin{equation*}
\mathcal{S}_2:=\big{\{}\bm{x}\in\mathbb{R}^3:\sqrt{(x_1+x_2+x_3-3)^2+(x_1+x_2-x_3-3)^2}\leq1\big{\}}.
\end{equation*}
Then,
\begin{equation*}
\begin{aligned}
{\rm Arg}\min\limits_{\bm{x}\in\mathcal{S}_2}\{\|\bm{x}\|_0\} &= {\textstyle\left\{(x_1, 0, 0)^{\top} : 3-\frac{\sqrt{2}}{2} \leq x_1 \leq 3+\frac{\sqrt{2}}{2}\right\} \cup \left\{(0, x_2, 0)^{\top} : 3-\frac{\sqrt{2}}{2} \leq x_2 \leq 3+\frac{\sqrt{2}}{2}\right\}},  \\
{\rm Arg}\min\limits_{\bm{x}\in\mathcal{S}_2}\{\|\bm{x}\|_p^p\} &= {\textstyle\left\{(3-\frac{\sqrt{2}}{2}, \,0, \,0)^{\top}, ~(0, \,3-\frac{\sqrt{2}}{2}, \,0)^{\top}\right\}} \quad \mathrm{for}~\mathrm{any}~0<p<1,  \\
{\rm Arg}\min\limits_{\bm{x}\in\mathcal{S}_2}\{\|\bm{x}\|_1\} &= {\textstyle\left\{(x_1, x_2, 0)^{\top}: x_1 + x_2 = 3-\frac{\sqrt{2}}{2}, \,x_1\geq 0, \,x_2\geq 0 \right\}}.
\end{aligned}
\end{equation*}

\textbf{For} $q=\infty$, the feasible set of \eqref{exam1} is
\begin{equation*}
\mathcal{S}_{\infty}:=\left\{\bm{x}\in\mathbb{R}^3:\max\{|x_1+x_2+x_3-3|,\, |x_1+x_2-x_3-3|\}\leq1\right\}.
\end{equation*}
Then,
\begin{equation*}
\begin{aligned}
{\rm Arg}\min\limits_{\bm{x}\in\mathcal{S}_{\infty}}\{\|\bm{x}\|_0\} &= \left\{(x_1, 0, 0)^{\top} : 2 \leq x_1 \leq 4\right\} \cup \left\{(0, x_2, 0)^{\top} : 2 \leq x_2 \leq 4\right\},  \\
{\rm Arg}\min\limits_{\bm{x}\in\mathcal{S}_{\infty}}\{\|\bm{x}\|_p^p\} &= \left\{(2, \,0, \,0)^{\top}, \,(0, \,2, \,0)^{\top}\right\} \quad \mathrm{for}~\mathrm{any}~0<p<1,  \\
{\rm Arg}\min\limits_{\bm{x}\in\mathcal{S}_{\infty}}\{\|\bm{x}\|_1\} &= \left\{(x_1, x_2, 0)^{\top} : x_1 + x_2 = 2, \,x_1\geq 0, \,x_2\geq 0 \right\}.
\end{aligned}
\end{equation*}

\end{example}

From this example, one can easily see that every optimal solution $\bm{x}^*$ of \eqref{exam1} is at the boundary of the feasible set for $0<p\leq1$ and there is a $\alpha \in (0,1]$ such that $\alpha\bm{x}^*$ is at the boundary of the feasible set for $p=0$, as claimed in Theorem \ref{thmeqvl0}(i). Moreover, every optimal solution of \eqref{exam1} with $0<p<1$ is exactly a sparsest solution over $\{\bm{x}\in\mathbb{R}^3:\|A\bm{x}-\bm{b}\|_q\leq1\}$ for $q=1,\,2,\,\infty$, while an optimal solution of \eqref{exam1} with $p=1$ may not be a sparest one. This shows the potential advantage of using the $L_p$-norm ($0<p<1$) to approximate the $L_0$-norm. In particular, when $q=1$, one can further estimate $p^*$ by \eqref{defpstar} for this example. Indeed, it is easy to see that $s=1$. Then, from \eqref{defr1}, we compute that $ r=\frac{1+3\sqrt{2}}{\sqrt{2}}$. Moreover, one can verify that
\begin{equation*}
\bigcup_{j \in \{1,\cdots,8\}}{\rm EXT}\left(\mathbb{P}_{j}\cap\mathcal{S}_1\right) = \left\{
\begin{array}{llll}
\left(\frac{7}{2}, \,0, \,\frac{1}{2}\right), &\left(0, \,\frac{7}{2}, \,\frac{1}{2}\right), &\left(\frac{7}{2}, \,0, \,-\frac{1}{2}\right), &\left(0, \,\frac{7}{2}, \,-\frac{1}{2}\right) \vspace{1mm}\\
\left(\frac{5}{2}, \,0, \,\frac{1}{2}\right), &\left(0, \,\frac{5}{2}, \,\frac{1}{2}\right), &\left(\frac{5}{2}, \,0, \,-\frac{1}{2}\right), &\left(0, \,\frac{5}{2}, \,-\frac{1}{2}\right) \vspace{1mm}\\
\left(\frac{7}{2}, \,0, \,0\right), &\left(0, \,\frac{7}{2}, \,0\right), &\left(\frac{5}{2}, \,0, \,0\right), &\left(0, \,\frac{5}{2}, \,0\right)
\end{array}
\right\}.
\end{equation*}
Thus, it follows from \eqref{defrt} that $\tilde{r}=\frac{1}{2}$. Now, using \eqref{defpstar}, since $\frac{r}{\tilde{r}}=6+\sqrt{2} > \frac{s+1}{s} = 2$, we have
\begin{equation*}
p^* = \frac{\ln(s+1)-\ln s}{\ln r - \ln \tilde{r}} = \frac{\ln 2}{\ln (6+\sqrt{2})} \approx 0.346.
\end{equation*}
Recalling Theorem \ref{pstarthm}, we know that every optimal solution of \eqref{exam1} with $p\in\left(0, \,\frac{\ln 2}{\ln (6+\sqrt{2})}\right)$ shall be an optimal solution of \eqref{exam1} with $p=0$. This is clearly evident in \eqref{egresults}. In fact, for this example, every optimal solution of \eqref{exam1} with $0<p<1$ is an optimal solution of \eqref{exam1} with $p=0$. This shows that $p^*$ given in \eqref{defpstar} may not be the optimal upper bound of $p$ such that ${\rm SOL}(A,\bm{b},\sigma,p,1)\subseteq {\rm SOL}(A,\bm{b},\sigma,0,1)$ for any $p\in(0, \,p^*)$. In addition, our current estimate $p^*$ in \eqref{defpstar} depends on the knowledge on the optimal value $s$, which may be unknown or difficult to find in practice. Fortunately, we observe that $p^*$, viewed as a function of $s$, is actually decreasing when $\frac{\ln(s+1)-\ln s}{\ln r - \ln \tilde{r}}\leq1$. Thus, one may estimate a proper upper bound $\tilde{s}$ for the true optimal value $s$ (i.e., $\tilde{s} \geq s$) and compute $\tilde{p}^*=\frac{\ln(\tilde{s}+1)-\ln \tilde{s}}{\ln r - \ln \tilde{r}}$ satisfying $\tilde{p}^* \leq p^*$. It then follows from Theorem \ref{pstarthm} that ${\rm SOL}(A,\bm{b},\sigma,p,1)\subseteq {\rm SOL}(A,\bm{b},\sigma,0,1)$ for any $p\in(0, \,\tilde{p}^*)$. But it should be noticed that such $\tilde{p}^*$ can be more conservative. Improving estimations of $p^*$ and $\tilde{p}^*$ will be an interesting research topic in the future.

%%%%%%%%%%%%%%%%%%%%%%%%%%%%%%%%%%%%%%%%%%
\section{A smoothing penalty method}\label{secpen}

In this section, we propose a smoothing penalty method for solving the $L_1$-constrained problem \eqref{lpl1model} with $0<p<1$. Before proceeding, we would like to point out that the smoothing penalty method presented in this paper can be extended without much difficulty to solve the $L_{\infty}$-constrained problem, namely, problem \eqref{lplqmodel} with $0<p<1$ and $q=\infty$. Because the $L_{\infty}$-constrained problem is similar to the $L_1$-constrained problem in the sense that both constraints $\|A\bm{x}-\bm{b}\|_1 \leq \sigma$ and $\|A\bm{x}-\bm{b}\|_{\infty} \leq \sigma$ are polyhedral constraints, and the functions $\bm{x}\mapsto\|A\bm{x}-\bm{b}\|_1$ and $\bm{x}\mapsto\|A\bm{x}-\bm{b}\|_\infty$ are piecewise linear. On the other hand, for $1<q<\infty$, the function $\bm{x}\mapsto\|A\bm{x}-\bm{b}\|^q_q$ is continuously differentiable. Then, one can readily extend the smoothing penalty method proposed in \cite{clp2016penalty} to solve problem \eqref{lplqmodel} with $0<p<1$ and $1<q<\infty$. However, the approach in \cite{clp2016penalty} cannot be directly adapted for $q\in\{1,\infty\}$ due to the nonsmoothness of the function $\bm{x} \mapsto \|A\bm{x}-\bm{b}\|_q$ in these two cases. In view of the above, in this paper, we consider an alternative smoothing penalty method for solving the $L_1$-constrained problem and omit the discussions on solving the $L_{\infty}$-constrained problem to save space.

We first study the first-order optimality conditions for problem \eqref{lpl1model} with $0<p<1$. For simplicity, from now on, let $\Phi(\bm{x}):=\|\bm{x}\|_p^p$. Then, problem \eqref{lpl1model} with $0<p<1$ can be equivalently written as follows:
\begin{equation}\label{lpl1modelre}
\min\limits_{\bm{x}\in\mathbb{R}^n}\,\,\Phi(\bm{x}) + \delta_{{\rm FEA}(A,\bm{b},\sigma,1)}(\bm{x}).
\end{equation}
It is known from the generalized Fermat's rule
\cite[Theorem 10.1]{rw1998variational} that, at any local minimizer $\bar{\bm{x}}$ of \eqref{lpl1modelre} (hence \eqref{lpl1model}), the following first-order necessary condition holds:
\begin{equation}\label{optcondlpl1}
0 \in \partial \big{(} \Phi + \delta_{{\rm FEA}(A,\bm{b},\sigma,1)} \big{)}(\bar{\bm{x}}).
\end{equation}
This motivates the following definition.

\begin{definition}[\textbf{Stationary point of problem \eqref{lpl1model} with $0<p<1$}]\label{defKKTpoint}
A point $\bm{x}^*$ is said to be a stationary point of problem \eqref{lpl1model} with $0<p<1$ if $\bm{x}^*\in{\rm FEA}(A,\bm{b},\sigma,1)$ and \eqref{optcondlpl1} is satisfied with $\bm{x}^*$ in place of $\bar{\bm{x}}$.
\end{definition}

\begin{figure}[ht]
\centering
\includegraphics[width=7cm]{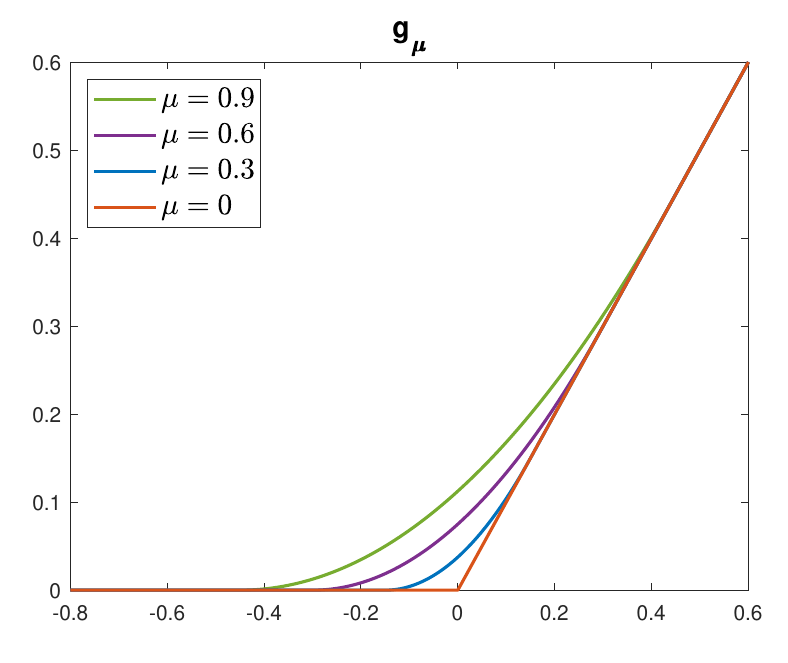}~~~~
\includegraphics[width=7cm]{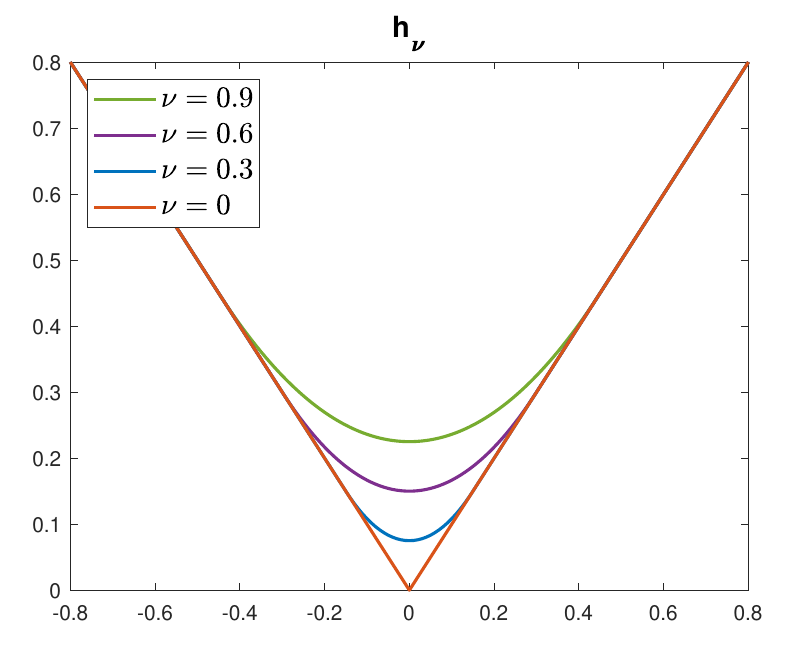}
\caption{Graphs of the smoothing functions $g_{\mu}$ and $h_{\nu}$ with different values of $\mu$ and $\nu$.}\label{figsmfuncs}
\end{figure}

Note that finding an optimal solution of problem \eqref{lpl1model} with $0<p<1$ is NP-hard \cite{cgwy2014complexity,gjy2011note}. Therefore, we shall focus on finding a stationary point of this problem. To this end, we introduce the following auxiliary penalty problem:
\begin{equation}\label{penprob}
\min\limits_{\bm{x}\in\mathbb{R}^n}~F_{\lambda}(\bm{x}) := \Phi(\bm{x}) + \lambda (\|A\bm{x} - \bm{b}\|_1 - \sigma)_+,
\end{equation}
where $\lambda>0$ is the penalty parameter and $(\cdot)_+:=\max\{\cdot, \,0\}$. This problem is indeed an exact penalty problem for problem \eqref{lpl1model} with $0<p<1$. The detailed analysis for the exact penalization results regarding global and local minimizers is given in Appendix \ref{appexact}. However, problem \eqref{penprob} is still not conceivably solvable because both parts in \eqref{penprob} are nonsmooth, and moreover, $\Phi$ is nonconvex and non-Lipschitz. We then consider a partially smoothing problem of \eqref{penprob} as follows:
\begin{equation}\label{penprobsmooth}
\min\limits_{\bm{x}\in\mathbb{R}^n}~F_{\lambda,\mu,\nu}(\bm{x}):=\Phi(\bm{x}) + f_{\lambda,\mu,\nu}(\bm{x}),
\end{equation}
where $\mu$, $\nu>0$ are smoothing parameters and
\begin{equation*}
f_{\lambda,\mu,\nu}(\bm{x}):=\lambda\,g_{\mu}\big{(} H_{\nu}(A\bm{x}-\bm{b})-\sigma\big{)}
\end{equation*}
with
\begin{equation*}
\begin{aligned}
g_{\mu}(s) := \left\{\begin{array}{ll}
               \!\!(s)_+, &\mathrm{if}~|s| \geq \frac{\mu}{2},  \\[2pt]
               \!\!\frac{s^2}{2\mu}+\frac{s}{2} + \frac{\mu}{8},  &\mathrm{if}~|s| < \frac{\mu}{2},
               \end{array}\right.
\quad
H_{\nu}(\bm{z}) := \sum^m_{i=1}h_{\nu}(z_i),
\quad
h_{\nu}(t) := \left\{\begin{array}{ll}
               \!\!|t|, &\mathrm{if}~|t| \geq \frac{\nu}{2},   \\[2pt]
               \!\!\frac{t^2}{\nu}+\frac{\nu}{4},  &\mathrm{if}~|t| < \frac{\nu}{2}.
               \end{array}\right.
\end{aligned}
\end{equation*}
Note that $g_{\mu}(s)$ and $h_{\nu}(t)$ are the smoothing functions of $(s)_+$ and $|t|$, respectively (see Figure \ref{figsmfuncs}), and they have the following nice properties:
\begin{eqnarray}
&&0 \leq g_{\mu}(s) - (s)_+ \leq {\textstyle\frac{\mu}{8}}, \quad \forall\,s\in\mathbb{R}, \label{smooth1} \\ [4pt]
&&0 \leq h_{\nu}(t) - |t| \leq {\textstyle\frac{\nu}{4}}, \quad \forall\,t\in\mathbb{R},  \nonumber \\ [4pt]
&&0 \leq H_{\nu}(A\bm{x}-\bm{b}) - \|A\bm{x} - \bm{b}\|_1 \leq {\textstyle\frac{m}{4}}\nu, \quad \forall\,\bm{x}\in\mathbb{R}^n. \label{Hsmoothbd}
\end{eqnarray}
More details on these smoothing functions can be found in
\cite[Section 3]{c2012smoothing} and references therein. Thus, the composite function $f_{\lambda,\mu,\nu}(\bm{x})$ is indeed obtained by applying the smoothing technique \textit{twice}. Hence, it is continuously differentiable and can be viewed as a smoothing function of $\lambda (\|A\bm{x} - \bm{b}\|_1 - \sigma)_+$. One can also show that
\begin{equation}\label{smoothfun}
0 \leq f_{\lambda,\mu,\nu}(\bm{x}) - \lambda (\|A\bm{x} - \bm{b}\|_1 - \sigma)_+ \leq {\textstyle\frac{m}{4}}\lambda\nu + {\textstyle\frac{1}{8}}\lambda\mu.
\end{equation}
Moreover, it is worth mentioning that when $\sigma=0$, the auxiliary penalty problem \eqref{penprob} reduces to $\min\limits_{\bm{x}\in\mathbb{R}^n}\left\{\Phi(\bm{x}) + \lambda\|A\bm{x} - \bm{b}\|_1\right\}$. Then, the smoothing function $g_{\mu}$ of $(\cdot)_+$ is no longer needed and the subsequent analysis can also be simplified in this special case. Now, based on \eqref{penprobsmooth}, we are ready to present a smoothing penalty method as Algorithm \ref{algpen} for solving problem \eqref{lpl1model} with $0<p<1$. We call it SPeL1 for short in the rest of this paper.

\begin{algorithm}[ht]
\caption{A smoothing penalty method for solving \eqref{lpl1model} with $0<p<1$ (SPeL1)}\label{algpen}
\textbf{Input:} a feasible point $\bm{x}^{\mathrm{feas}}\in{\rm FEA}(A,\bm{b},\sigma,1)$, $\lambda_0>0$, $\mu_0>0$, $\nu_0>0$, $\epsilon_0>0$, $\rho>1$ and $0<\theta<1$. Set $k=0$ and $\bm{x}^{0}=\bm{x}^{\mathrm{feas}}$. \\
\textbf{while} a termination criterion is not met, \textbf{do}
\begin{itemize}[leftmargin=2cm]
\item[\textbf{Step 1}.] If $F_{\lambda_k,\mu_k,\nu_k}(\bm{x}^{k}) \!\leq\! F_{\lambda_k,\mu_k,\nu_k}(\bm{x}^{\mathrm{feas}})$, set $\bm{x}^{k,0}=\bm{x}^{k}$; otherwise, set $\bm{x}^{k,0}=\bm{x}^{\mathrm{feas}}$.

\item[\textbf{Step 2}.] Apply certain method with $\bm{x}^{k,0}$ as the initial point to find an approximate first-order stationary point $\bm{x}^{k,l_k}$ of \eqref{penprobsmooth} with $(\lambda_k, \mu_k, \nu_k)$ such that
    \begin{eqnarray}
    \mathrm{dist}\left(0,\,\partial \Phi(\bm{x}^{k,l_k+1}) + \nabla f_{\lambda_k,\mu_k,\nu_k}(\bm{x}^{k,l_k})\right)
    &\leq& \epsilon_k, \label{suboptcond} \\ [4pt]
    \left\|\bm{x}^{k,l_k+1}-\bm{x}^{k,l_k}\right\|_2 &\leq& \epsilon_k, \label{suboptsucc} \\ [4pt]
    F_{\lambda_k,\mu_k,\nu_k}(\bm{x}^{k,l_k}) - F_{\lambda_k,\mu_k,\nu_k}(\bm{x}^{k,0}) &\leq& 0. \label{suboptcondadd}
    \end{eqnarray}

\item[\textbf{Step 3}.] Set $\bm{x}^{k+1}=\bm{x}^{k,l_k}$, $\lambda_{k+1}=\rho\lambda_k$, $\mu_{k+1}=\theta\mu_k$, $\nu_{k+1}=\theta\nu_k$ and $\epsilon_{k+1}=\theta\epsilon_k$.

\item [\textbf{Step 4}.] Set $k = k+1$ and go to \textbf{Step 1}.
\end{itemize}
\textbf{end while}  \\
\textbf{Output}: $\bm{x}^k$ \vspace{0.5mm}
\end{algorithm}

The reader may have observed that, since problem \eqref{penprob} is a penalty counterpart of problem \eqref{lpl1model} and problem \eqref{penprobsmooth} is a partially smoothing counterpart of problem \eqref{penprob}, our method actually adapts the \textit{penalty} strategy and the \textit{smoothing} strategy at the same time for solving the nonconvex nonsmooth non-Lipschitz constrained problem \eqref{lpl1model} with $0<p<1$. Specifically, in our method, at each iteration, we solve problem \eqref{penprobsmooth} approximately with given $(\lambda,\mu,\nu)$, and then update $\bm{x}$ and $(\lambda,\mu,\nu)$. The cooperation of these two strategies indeed provides an efficient practical way to solve problem \eqref{lpl1model} with $0<p<1$. This circumvents the potential disadvantages of the traditional penalty approach that directly solves the penalty problem \eqref{penprob} with an exact penalty parameter $\lambda^*$, because (i) it is still not easy to solve problem \eqref{penprob} efficiently; (ii) it is, in general, hard to estimate the exact penalty parameter $\lambda^*$ and the overestimation may make the penalty problem \eqref{penprob} ill-conditioned. The convergence result that characterizes a cluster point of the sequence generated by the SPeL1 in Algorithm \ref{algpen} is shown in the next theorem. We should note that, though the proofs are motivated by those in
\cite[Theorem 4.2]{clp2016penalty} and \cite[Theorem 2]{lpt2017a}, the technical details become much more involved since our smoothing function $f_{\lambda,\mu,\nu}$ is obtained by a composition of two smoothing functions $g_{\mu}$ and $H_{\nu}$.

For the ease of future reference, we write down the gradients of $f_{\lambda,\mu,\nu}$ and $H_{\nu}$ as well as the derivatives of $g_{\mu}$ and $h_{\nu}$ as follows:
\begin{eqnarray}
\nabla f_{\lambda,\mu,\nu}(\bm{x}) &=& \lambda\,g'_{\mu}\left( H_{\nu}(A\bm{x}-\bm{b})-\sigma\right)A^{\top}\nabla H_{\nu}(A\bm{x}-\bm{b}), \label{gradf} \\ [4pt]
g'_{\mu}(s) &=& \min\big{\{}\max\big{\{}{\textstyle\frac{s}{\mu}+\frac{1}{2}}, \,0\big{\}}, \,1\big{\}},  \label{defgde}\\ [4pt]
\nabla H_{\nu}(\bm{z}) &=& \big{(}h'_{\nu}(z_1), \cdots, h'_{\nu}(z_m)\big{)}^{\top},  \nonumber\\ [4pt]
h'_{\nu}(t) &=& \min\big{\{}\max\big{\{}{\textstyle\frac{2}{\nu}}\,t, \,-1\big{\}}, \,1\big{\}}.  \label{defhde}
\end{eqnarray}
Moreover, we claim that $\Phi$ is regular at any $\bm{x}\in\mathbb{R}^n$ as follows. Let $\phi(t)=|t|^p$ for any $t\in\mathbb{R}$. It is easy to see that $\phi(t)$ is regular at any $t\neq0$, because $\phi(t)$ is smooth in a neighborhood of any $t\neq0$; see \cite[Exercise 8.8]{rw1998variational} and \cite[Corollary 8.11]{rw1998variational}. For $t=0$, it follows from \cite[Lemma 2.5]{clp2016penalty} and its proof that $\widehat{\partial} \phi(0)=\partial \phi(0)=\partial^{\infty} \phi(0) = \mathbb{R}$. Moreover, from the definition of the horizon cone (see \cite[Definition 3.3]{rw1998variational}), we have that $\widehat{\partial} \phi(0)^{\infty}=\mathbb{R}$. Using these facts and \cite[Corollary 8.11]{rw1998variational}, we see that $\phi(t)$ is also regular at $t=0$. Therefore, it follows from \cite[Proposition 10.5]{rw1998variational} that $\Phi$ is regular at any $\bm{x}\in\mathbb{R}^n$.

\begin{theorem}\label{thmconv}
Suppose that $\rho>1$ and $0<\theta<1$ are chosen such that $\theta\rho\leq1$. Let $\{\bm{x}^k\}^{\infty}_{k=0}$ be the sequence generated by the SPeL1 in Algorithm \ref{algpen}. Then, the following statements hold.
\begin{itemize}[leftmargin=0.8cm]
\item[(i)] $\{\bm{x}^k\}$ is bounded.

\item[(ii)] Any cluster point $\bm{x}^*$ of $\{\bm{x}^k\}$ is a feasible point of problem \eqref{lpl1model} with $0<p<1$.

\item[(iii)] Suppose that $\bm{x}^*$ is a cluster point of $\{\bm{x}^k\}$ and it holds at $\bm{x}^*$ that
    \begin{equation}\label{CQcond}
    -\,\partial^{\infty}\Phi(\bm{x}^*)\,\cap\,\mathcal{N}_{{\rm FEA}(A,\bm{b},\sigma,1)}(\bm{x}^*) = \{0\}.
    \end{equation}
    Then, $\bm{x}^*$ is a stationary point of problem \eqref{lpl1model} with $0<p<1$.
\end{itemize}
\end{theorem}
\beginproof
\textit{Statement (i)}. First, we see that
\begin{equation*}
\begin{aligned}
\Phi(\bm{x}^{k+1})
&\leq F_{\lambda_k,\mu_k,\nu_k}(\bm{x}^{k+1})
= F_{\lambda_k,\mu_k,\nu_k}(\bm{x}^{k,l_k})
\leq F_{\lambda_k,\mu_k,\nu_k}(\bm{x}^{k,0})
\leq F_{\lambda_k,\mu_k,\nu_k}(\bm{x}^{\mathrm{feas}}) \\[2pt]
&
\leq \Phi(\bm{x}^{\mathrm{feas}}) + {\textstyle\frac{1}{8}}\lambda_k(2m\nu_k + \mu_k)
= \Phi(\bm{x}^{\mathrm{feas}}) + {\textstyle\frac{1}{8}}\lambda_0(2m\nu_0 + \mu_0)(\theta\rho)^k  \\[2pt]
&\leq \Phi(\bm{x}^{\mathrm{feas}}) + {\textstyle\frac{1}{8}}\lambda_0(2m\nu_0 + \mu_0),
\end{aligned}
\end{equation*}
where the first inequality follows from the nonnegativity of $f_{\lambda_k,\mu_k,\nu_k}(\bm{x}^{k+1})$ (since $g_{\mu}(s)\geq0$ for all $s$), the second inequality follows from \eqref{suboptcondadd}, the third inequality follows from Step 1 in Algorithm \ref{algpen}, the fourth inequality follows from \eqref{smoothfun} and the last inequality follows from $\theta\rho\leq1$. This together with the level-boundedness of $\Phi$ (recall that $\Phi(\bm{x}):=\|\bm{x}\|_p^p$) implies that $\{\bm{x}^k\}$ is bounded.

\textit{Statement (ii)}. Since $\{\bm{x}^k\}$ is bounded, there exists at least one cluster point. Suppose that $\bm{x}^*$ is a cluster point of $\{\bm{x}^k\}$ and let $\{\bm{x}^{k_i}\}$ be a convergent subsequence such that $\lim\limits_{i\to\infty} \bm{x}^{k_i} = \bm{x}^*$. Note that
\begin{equation*}
\begin{aligned}
\lambda_{k-1} \left(\|A\bm{x}^{k} - \bm{b}\|_1 - \sigma\right)_+
&\leq f_{\lambda_{k-1},\mu_{k-1},\nu_{k-1}}(\bm{x}^{k})
\leq F_{\lambda_{k-1},\mu_{k-1},\nu_{k-1}}(\bm{x}^{k})
\leq F_{\lambda_{k-1},\mu_{k-1},\nu_{k-1}}(\bm{x}^{\mathrm{feas}}) \\[2pt]
&
\leq \Phi(\bm{x}^{\mathrm{feas}}) + {\textstyle\frac{1}{8}}\lambda_{k-1}(2m\nu_{k-1} + \mu_{k-1}), 
\end{aligned}
\end{equation*}
where the first inequality follows from \eqref{smooth1}, \eqref{Hsmoothbd} and the fact that $g_{\mu}$ is non-decreasing, and the last inequality follows from \eqref{smoothfun}. Then,
\begin{equation*}
(\|A\bm{x}^{k} - \bm{b}\|_1 - \sigma)_+ \leq \frac{\Phi(\bm{x}^{\mathrm{feas}})}{\lambda_{k-1}} + \frac{m}{4}\nu_{k-1} + \frac{1}{8}\mu_{k-1}.
\end{equation*}
Taking limit in above inequality along $\{\bm{x}^{k_i}\}$ and recalling that $\lambda_{k_i-1}\rightarrow\infty$, $\mu_{k_i-1}\to0$, $\nu_{k_i-1}\to0$ (see Step 3 in Algorithm \ref{algpen}), we see that $\|A\bm{x}^* - \bm{b}\|_1 \leq \sigma$. Hence, $\bm{x}^*$ is a feasible point of \eqref{lpl1model} with $0<p<1$.

\textit{Statement (iii)}. We next show that $\bm{x}^*$ is a stationary point of problem \eqref{lpl1model} with $0<p<1$. For simplicity, let $\bm{a}_j\in\mathbb{R}^n$ ($j=1,\cdots,m$) be the column vector formed from the $j$th row of $A$, i.e., $A=[\bm{a}_1, \cdots, \bm{a}_m]^{\top}\in\mathbb{R}^{m\times n}$. Moreover, let $\bm{y}^{k+1}:=\bm{x}^{k,l_k+1}$. Then, $\lim\limits_{i\to\infty} \bm{y}^{k_i} = \bm{x}^*$ thanks to $\bm{x}^{k_i}\to\bm{x}^*$ and \eqref{suboptsucc} with $\epsilon_k\to0$. Thus, from \eqref{suboptcond} and \eqref{gradf}, we see that for any $k\geq1$, there exists a $\bm{\xi}^k \in \partial \Phi(\bm{y}^{k})$ such that
\begin{equation}\label{suboptcond2}
\begin{aligned}
&\quad\big{\|}\bm{\xi}^{k} + \nabla f_{\lambda_{k-1},\mu_{k-1},\nu_{k-1}}(\bm{x}^{k})\big{\|} \\
&=\big{\|}\bm{\xi}^{k} + \lambda_{k-1}g'_{\mu_{k-1}}\left( H_{\nu_{k-1}}(A\bm{x}^k-\bm{b})-\sigma\right)A^{\top}\nabla H_{\nu_{k-1}}(A\bm{x}^k-\bm{b})\big{\|}   \\
&=\big{\|}\bm{\xi}^{k} + \lambda_{k-1}g'_{\mu_{k-1}}\big{(} H_{\nu_{k-1}}(A\bm{x}^k-\bm{b})-\sigma\big{)}{\textstyle\sum^{m}_{j=1}}h'_{\nu_{k-1}}
\big{(}[A\bm{x}^k-\bm{b}]_j\big{)}\bm{a}_j\big{\|}  \\
&\leq \epsilon_{k-1}.
\end{aligned}
\end{equation}
In the following, we consider two cases: $\|A\bm{x}^* - \bm{b}\|_1 < \sigma$ and $\|A\bm{x}^* - \bm{b}\|_1 = \sigma$.

\textbf{Case 1.} In this case, we suppose that $\|A\bm{x}^* - \bm{b}\|_1 < \sigma$. Since $\|A\bm{x}^{k_i} - \bm{b}\|_1 \to \|A\bm{x}^* - \bm{b}\|_1$, then, for any $0<\gamma<\sigma-\|A\bm{x}^* - \bm{b}\|_1$, there exists a sufficiently large $K_{\gamma}>0$ such that $\big{|}\|A\bm{x}^{k_i} - \bm{b}\|_1-\|A\bm{x}^* - \bm{b}\|_1\big{|}\leq\gamma$ for all $k_i \geq K_{\gamma}$. Note that
\begin{equation*}
\begin{aligned}
&~~\frac{H_{\nu_{k_i-1}}(A\bm{x}^{k_i}-\bm{b}) - \sigma}{\mu_{k_i-1}}+\frac{1}{2}
\leq \frac{\|A\bm{x}^{k_i} - \bm{b}\|_1 - \sigma + \frac{m}{4}\nu_{k_i-1}}{\mu_{k_i-1}}+\frac{1}{2}  \\
&=\frac{\|A\bm{x}^{k_i} - \bm{b}\|_1 - \sigma}{\mu_{k_i-1}} + \frac{m}{4}\frac{\nu_{0}}{\mu_{0}}+\frac{1}{2}
\leq \frac{\|A\bm{x}^* - \bm{b}\|_1-\sigma+\gamma}{\mu_{k_i-1}} + \frac{m}{4}\frac{\nu_{0}}{\mu_{0}}+\frac{1}{2}
< 0,
\end{aligned}
\end{equation*}
where the first inequality follows from \eqref{Hsmoothbd}, the equality follows from $\frac{\nu_{k}}{\mu_{k}}=\frac{\theta\nu_{k-1}}{\theta\mu_{k-1}}=\cdots=\frac{\nu_{0}}{\mu_{0}}$, the second inequality holds for all $k_i \geq K_{\gamma}$ and the last inequality follows whenever $k_i\geq\widetilde{K}_{\gamma}$ for some $\widetilde{K}_{\gamma}\geq K_{\gamma}$ because $\mu_{k_i}\to0$ and $\|A\bm{x}^* - \bm{b}\|_1-\sigma+\gamma<0$. This together with \eqref{defgde} implies that $g'_{\mu_{k_i-1}}\big{(}H_{\nu_{k_i-1}}(A\bm{x}^{k_i}-\bm{b})-\sigma\big{)} = 0$ for all sufficiently large $k_i$. Hence, \eqref{suboptcond2} reduces to $\|\bm{\xi}^{k_i}\| \leq \epsilon_{k_i-1}$ for all sufficiently large $k_i$. Then, we have from \eqref{robust} that $0=\bm{\xi}^*\in\partial\Phi(\bm{x}^*)$. This together with $\mathcal{N}_{{\rm FEA}(A,\bm{b},\sigma,1)}(\bm{x}^*)=\{0\}$ (since $\|A\bm{x}^* - \bm{b}\|_1 < \sigma$) implies that
\begin{equation*}
0 \in \partial\Phi(\bm{x}^*) + \mathcal{N}_{{\rm FEA}(A,\bm{b},\sigma,1)}(\bm{x}^*).
\end{equation*}
Moreover, since $\Phi$ and $\delta_{{\rm FEA}(A,\bm{b},\sigma,1)}$ are regular, then it follows from \cite[Corollary 8.11]{rw1998variational} and \cite[Exercise 8.14]{rw1998variational} that $\partial\Phi(\bm{x}^*)=\widehat{\partial}\Phi(\bm{x}^*)$ and $\mathcal{N}_{{\rm FEA}(A,\bm{b},\sigma,1)}(\bm{x}^*)=\partial \delta_{{\rm FEA}(A,\bm{b},\sigma,1)}(\bm{x}^*)=\widehat{\partial}\delta_{{\rm FEA}(A,\bm{b},\sigma,1)}(\bm{x}^*)$. Using these facts and recalling \cite[Theorem 8.6]{rw1998variational}, \cite[Corollary 10.9]{rw1998variational}, we have
\begin{equation*}
\begin{aligned}
0 &\in \partial\Phi(\bm{x}^*) + \mathcal{N}_{{\rm FEA}(A,\bm{b},\sigma,1)}(\bm{x}^*)
=\widehat{\partial}\Phi(\bm{x}^*) + \widehat{\partial}\delta_{{\rm FEA}(A,\bm{b},\sigma,1)}(\bm{x}^*) \\
&\subseteq \widehat{\partial} \left(\Phi+\delta_{{\rm FEA}(A,\bm{b},\sigma,1)}\right)(\bm{x}^*)
\subseteq \partial \left(\Phi+\delta_{{\rm FEA}(A,\bm{b},\sigma,1)}\right)(\bm{x}^*),
\end{aligned}
\end{equation*}
which implies that $\bm{x}^*$ is a stationary point of problem \eqref{lpl1model} with $0<p<1$.

\textbf{Case 2.} In this case, we suppose that $\|A\bm{x}^* - \bm{b}\|_1 = \sigma$. For such $\bm{x}^*$, one can follow
\cite[Theorem 1.3.5 in Section D]{hl2001fundamentals} to compute that
\begin{equation}\label{normalcone}
\begin{aligned}
\mathcal{N}_{{\rm FEA}(A,\bm{b},\sigma,1)}(\bm{x}^*)
&= \left\{ c A^{\top}\bm{d} \,:\, \bm{d} \in \partial\|\cdot\|_1(A\bm{x}^*-\bm{b}), \,c \geq 0\right\} \\[2pt]
&= \left\{ \sum_{j=1}^m z_j\bm{a}_j\,:\,c\geq0, ~\begin{array}{ll}
z_j=c,   &~\mathrm{if}~~ [A\bm{x}^*-\bm{b}]_j > 0, \\
z_j\in[-c,\,c], &~\mathrm{if}~~ [A\bm{x}^*-\bm{b}]_j = 0, \\
z_j=-c,  &~\mathrm{if}~~ [A\bm{x}^*-\bm{b}]_j < 0, \\
\end{array}
~\mathrm{for}~j=1,\cdots,m\right\}.
\end{aligned}
\end{equation}
For simplicity, let $\tilde{t}_k:=\lambda_{k-1}g'_{\mu_{k-1}}( H_{\nu_{k-1}}(A\bm{x}^k-\bm{b})-\sigma)$ and $t^j_k:=\tilde{t}_k h'_{\nu_{k-1}}([A\bm{x}^k-\bm{b}]_j)$ for $j=1,\cdots,m$. Also, let $\mathcal{J}^0:=\{j:[A\bm{x}^*-\bm{b}]_j=0\}$, $\mathcal{J}^+:=\{j:[A\bm{x}^*-\bm{b}]_j>0\}$ and $\mathcal{J}^-:=\{j:[A\bm{x}^*-\bm{b}]_j<0\}$. Then, \eqref{suboptcond2} is equivalent to
\begin{equation}\label{suboptcond2re}
\big{\|}\bm{\xi}^k + \textstyle{\sum_{j=1}^m} t^j_k \bm{a}_j\big{\|} \leq \epsilon_{k-1}.
\end{equation}
Since $\bm{x}^{k_i}\to\bm{x}^*$ and $\nu_{k_i-1}\to0$, there exists a sufficiently large $K>0$ such that for all $k_i\geq K$, we have $[A\bm{x}^{k_i} - \bm{b}]_j>0$ and $\frac{2}{\nu_{k_i-1}}[A\bm{x}^{k_i} - \bm{b}]_j\geq1$ for all $j\in\mathcal{J}^+$, and have $[A\bm{x}^{k_i} - \bm{b}]_j<0$ and $\frac{2}{\nu_{k_i-1}}[A\bm{x}^{k_i} - \bm{b}]_j\leq-1$ for all $j\in\mathcal{J}^-$. Thus, it follows from \eqref{defhde} that for all $k_i\geq K$, we have $h'_{\nu_{k_i-1}}([A\bm{x}^{k_i}-\bm{b}]_j)=1$ for all $j\in\mathcal{J}^+$ and $h'_{\nu_{k_i-1}}([A\bm{x}^{k_i}-\bm{b}]_j)=-1$ for all $j\in\mathcal{J}^-$. Moreover, for all $k_i\geq1$, we see from \eqref{defgde} and \eqref{defhde} that $g'_{\mu_{k_i-1}}( H_{\nu_{k_i-1}}(A\bm{x}^{k_i}-\bm{b})-\sigma) \geq 0$ and  $h'_{\nu_{k_i-1}}([A\bm{x}^{k_i}-\bm{b}]_j)\in[-1,1]$ for all $j$. Then, for all $k_i\geq K$, we have that $t^j_{k_i}=\tilde{t}_{k_i}\geq0$ for all $j\in\mathcal{J}^+$, $t^j_{k_i}=-\tilde{t}_{k_i}\leq0$ for all $j\in\mathcal{J}^-$ and $t^j_{k_i}\in[-\tilde{t}_{k_i}, \,\tilde{t}_{k_i}]$ for all $j\in\mathcal{J}^0$.

We next prove by contradiction that $\{\bm{\xi}^{k_i}\}$ is bounded. Suppose that $\{\bm{\xi}^{k_i}\}$ is unbounded. Without loss of generality, we assume that $\|\bm{\xi}^{k_i}\|\to\infty$ and that $\frac{1}{\|\bm{\xi}^{k_i}\|}\bm{\xi}^{k_i}\to\bm{\xi}^*$ for some $\bm{\xi}^*$. Then, it follows from \eqref{suboptcond2re} that
\begin{equation}\label{suboptcond3}
\left\|\frac{1}{\|\bm{\xi}^{k_i}\|}\bm{\xi}^{k_i} + \sum^m_{j=1}\frac{t_{k_i}^j}{\|\bm{\xi}^{k_i}\|}\bm{a}_j\right\| \leq \frac{\epsilon_{k_i-1}}{\|\bm{\xi}^{k_i}\|}.
\end{equation}
Moreover, from the discussions in the last paragraph, for all $k_i\geq K$, we have that $t_{k_i}^j/\|\bm{\xi}^{k_i}\|=\tilde{t}_{k_i}/\|\bm{\xi}^{k_i}\|\geq0$ for all $j\in\mathcal{J}^+$, $t_{k_i}^j/\|\bm{\xi}^{k_i}\|=-\tilde{t}_{k_i}/\|\bm{\xi}^{k_i}\|\leq0$ for all $j\in\mathcal{J}^-$ and $t_{k_i}^j/\|\bm{\xi}^{k_i}\|\in\left[-\tilde{t}_{k_i}/\|\bm{\xi}^{k_i}\|, \,\tilde{t}_{k_i}/\|\bm{\xi}^{k_i}\|\right]$ for all $j\in\mathcal{J}^0$. Then, it follows from \eqref{normalcone} that
\begin{equation*}
\sum^m_{j=1}\frac{t_{k_i}^j}{\|\bm{\xi}^{k_i}\|}\bm{a}_j \in \mathcal{N}_{{\rm FEA}(A,\bm{b},\sigma,1)}(\bm{x}^*)
\end{equation*}
for all $k_i\geq K$. Then, passing to the limit in \eqref{suboptcond3} along $\{\bm{x}^{k_i}\}$, together with $\frac{\epsilon_{k_i-1}}{\|\bm{\xi}^{k_i}\|} \to 0$ and the closeness of $\mathcal{N}_{{\rm FEA}(A,\bm{b},\sigma,1)}(\bm{x}^*)$, it is not hard to see that
\begin{equation*}
\bm{\xi}^* \in \partial^{\infty} \Phi(\bm{x}^*) \quad \mathrm{and} \quad
-\bm{\xi}^* \in \mathcal{N}_{{\rm FEA}(A,\bm{b},\sigma,1)}(\bm{x}^*).
\end{equation*}
Since $\bm{\xi}^*\neq0$ due to $\|\bm{\xi}^*\|=1$, this is in contradiction to \eqref{CQcond}. Hence, $\{\bm{\xi}^{k_i}\}$ is bounded. Without loss of generality, assume that $\bm{\xi}^{k_i}\to\bm{\xi}^*$. Then, passing to the limit in \eqref{suboptcond2re} along $\{\bm{x}^{k_i}\}$ and $\{\bm{y}^{k_i}\}$, making use of \eqref{normalcone} and the closeness of $\mathcal{N}_{{\rm FEA}(A,\bm{b},\sigma,1)}(\bm{x}^*)$, recalling \eqref{robust}, we obtain that
\begin{equation*}
0 \in \partial\Phi(\bm{x}^*) + \mathcal{N}_{{\rm FEA}(A,\bm{b},\sigma,1)}(\bm{x}^*).
\end{equation*}
Thus, following the similar arguments in \textbf{Case 1}, one can show that $\bm{x}^*$ is a stationary point of problem \eqref{lpl1model} with $0<p<1$. This completes the proof.
\endproof

\begin{remark}[\textbf{Comments on condition \eqref{CQcond}}]
Condition \eqref{CQcond} used for Theorem \ref{thmconv}(iii) is actually a classic constraint qualification for nonconvex nonsmooth optimization problems; see \cite[Theorem 8.15]{rw1998variational}. Note that, for any $\bm{x}^*\in{\rm FEA}(A,\bm{b},\sigma,1)$, we have
\begin{equation*}
\begin{aligned}
\mathcal{N}_{{\rm FEA}(A,\bm{b},\sigma,1)}(\bm{x}^*)
&=\left\{\begin{aligned}
&\big{\{} c A^{\top}\bm{d} \,:\, \bm{d} \in \partial\|\cdot\|_1(A\bm{x}^*-\bm{b}), \,c \geq 0\big{\}}\neq\{0\},  &&\mathrm{if}~~\|A\bm{x}^* - \bm{b}\|_1 = \sigma,  \\[2pt]
&\big{\{}0\big{\}},  &&\mathrm{if}~~\|A\bm{x}^* - \bm{b}\|_1 < \sigma.
\end{aligned}\right.
\end{aligned}
\end{equation*}
Moreover, recall from \cite[Lemma 2.5(ii)]{clp2016penalty} that
\begin{equation*}
\partial^{\infty}\Phi(\bm{x}^*) = \big{\{}\bm{v}\in\mathbb{R}^n : v_i=0~\mathrm{for}~i\in\mathrm{supp}(\bm{x}^*)\big{\}}.
\end{equation*}
Thus, condition \eqref{CQcond} obviously holds at a point $\bm{x}^*$ satisfying $\|A\bm{x}^*-\bm{b}\|_1<\sigma$. For a point $\bm{x}^*$ satisfying $\|A\bm{x}^*-\bm{b}\|_1=\sigma$, one sufficient condition for \eqref{CQcond} is that, for some $i\in\mathrm{supp}(\bm{x}^*)$, $[A^{\top}\bm{d}]_i\neq0$ holds for any $\bm{d}\in\partial\|\cdot\|_1(A\bm{x}^*-\bm{b})$, i.e., $\mathrm{Diag}(\bm{x}^*)A^{\top}\bm{d} \neq 0$ for any $\bm{d}\in\partial\|\cdot\|_1(A\bm{x}^*-\bm{b})$.
\end{remark}

To end this section, we briefly discuss the method for approximately solving the smoothing penalty problem \eqref{penprobsmooth} such that conditions \eqref{suboptcond}--\eqref{suboptcondadd} hold. Note that, for any given $(\lambda,\mu,\nu)$, $F_{\lambda,\mu,\nu}$ is a continuous function that consists of a nonconvex nonsmooth non-Lipschitz function $\Phi$ and a smooth function $f_{\lambda,\mu,\nu}$. It is also not hard to verify that the gradient of $f_{\lambda,\mu,\nu}$ is Lipschitz continuous. Moreover, $F_{\lambda,\mu,\nu}$ is level-bounded because $\Phi$ is level-bounded and $f_{\lambda,\mu,\nu}$ is nonnegative since $g_{\mu}$ is nonnegative. Hence, the well-known proximal gradient method and its variants are suitably applied for solving \eqref{penprobsmooth} with convergence guarantee; see, for example, \cite{abs2013convergence,beck2017first,clp2016penalty,y2017proximal}. In our numerical experiments, we follow \cite{clp2016penalty} to adapt the nonmonotone proximal gradient (NPG) method. The NPG method is basically the proximal gradient method with a non-monotone line search technique and allows the occasional increases in objective. By incorporating this technique, the NPG has been shown to have more favorable numerical performance over the monotone version in many applications; see, for example, \cite{gzlhy2013a,wnf2009sparse,y2017proximal}. The iterative scheme of the NPG for solving \eqref{penprobsmooth} with $(\lambda_k,\mu_k,\nu_k)$ is given as follows:
\vspace{2mm}
\begin{center}
\fbox{\parbox{0.96\linewidth}{
Choose $L_k^{\max} \geq L_k^{\min} > 0$, $\tau>1$, $c>0$, and an integer $N\geq0$. At the $l$-th ($l\geq0$) iteration, choose $L_{k,l}^0\in[L_k^{\min},L_k^{\max}]$ and find the smallest nonnegative integer $i_l$ such that
\begin{equation}\label{npgiter}
\left\{\begin{aligned}
&\bm{w} \in \arg\min\limits_{\bm{x}\in\mathbb{R}^n} \Big{\{} \Phi(\bm{x}) + \langle \nabla f_{\lambda_k,\mu_k,\nu_k}(\bm{x}^{k,l}), \,\bm{x}\rangle + \frac{\tau^{i_l}L_{k,l}^0}{2}\|\bm{x}-\bm{x}^{k,l}\|^2 \Big{\}}, \\
&F_{\lambda_k,\mu_k,\nu_k}(\bm{w}) - \max\limits_{[l-N]_{+}\leq i\leq l}F_{\lambda_k,\mu_k,\nu_k}(\bm{x}^{k,i}) \leq -\frac{c}{2}\|\bm{w}-\bm{x}^{k,l}\|^2.
\end{aligned}\right.
\end{equation}
Then, set $\bm{x}^{k,l+1}=\bm{w}$ and $\bar{L}_{k,l}=\tau^{i_l}L_{k,l}^0$.
}}
\end{center}
\vspace{2mm}
One can also show that, for any given $(\lambda_k,\mu_k,\nu_k)$ and $\epsilon_k$, a point $\bm{x}^{k,l_k}$ satisfying conditions \eqref{suboptcond}--\eqref{suboptcondadd} can be found by the NPG within a finite number of iterations. Indeed, it follows from \cite[Proposition A.1(i)]{clp2016penalty} that \eqref{suboptcondadd} holds for all $l\geq0$. Moreover, from the optimality condition of \eqref{npgiter}, we see that
\begin{equation*}
\begin{aligned}
&0 \in \partial \Phi(\bm{x}^{k,l+1}) + \nabla f_{\lambda_k,\mu_k,\nu_k}(\bm{x}^{k,l}) + \bar{L}_{k,l}(\bm{x}^{k,l+1}-\bm{x}^{k,l}), \\[3pt]
\Longrightarrow ~~ & -\bar{L}_{k,l}(\bm{x}^{k,l+1}-\bm{x}^{k,l}) \in \partial \Phi(\bm{x}^{k,l+1}) + \nabla f_{\lambda_k,\mu_k,\nu_k}(\bm{x}^{k,l}),
\end{aligned}
\end{equation*}
which implies that
\begin{equation}\label{optdistbd}
\begin{aligned}
\mathrm{dist}\left(0,\,\partial \Phi(\bm{x}^{k,l+1}) + \nabla f_{\lambda_k,\mu_k,\nu_k}(\bm{x}^{k,l})\right)
\leq \bar{L}_{k,l}\|\bm{x}^{k,l+1}-\bm{x}^{k,l}\|.
\end{aligned}
\end{equation}
This together with the boundedness of $\{\bar{L}_{k,l}\}_{l\geq0}$ (see \cite[Proposition A.1(ii)]{clp2016penalty}) and $\|\bm{x}^{k,l+1}-\bm{x}^{k,l}\|\to0$ as $l\to\infty$ (see \cite[Theorem A.1]{clp2016penalty}) implies that \eqref{suboptcond} and \eqref{suboptsucc} hold when $l$ is sufficiently large. In view of the above, the sequence $\{\bm{x}^k\}$ generated by the SPeL1 in Algorithm \ref{algpen} is well-defined.

%%%%%%%%%%%%%%%%%%%%%%%%%%%%%%%%%%%%%%%%%%%%%%%%%%%%%
\section{Numerical simulations}\label{secnum}

In this section, we conduct some numerical experiments for problem \eqref{lpl1model} with $0<p<1$ on finding sparse solutions to implicitly illustrate the theoretical results established in Section \ref{secpro} and show the efficiency of our SPeL1 in Algorithm \ref{algpen}. All experiments are run in {\sc Matlab} R2016a on a workstation with Intel(R) Xeon(R) Processor E-2176G@3.70GHz and 64GB of RAM, equipped with 64-bit Windows 10 OS.

For the SPeL1, we set $\lambda_0=\mu_0=\nu_0=1$ and $\bm{x}^0=\bm{x}^{\mathrm{feas}}=A^{\dagger}\bm{b}$, where the computation of $A^{\dagger}\bm{b}$ is not counted in the CPU time below. At the $k$th outer iteration, we compute
\begin{equation}\label{defeta}
\eta_1^k := \frac{\|\bm{x}^{k+1}-\bm{x}^{k}\|}{1+\|\bm{x}^{k+1}\|}, \quad
\eta_2^k := \frac{|\Phi(\bm{x}^{k+1})-\Phi(\bm{x}^{k})|}{1+\Phi(\bm{x}^{k+1})}, \quad
\eta_3^k := \max\big{\{}\|A\bm{x}^{k+1}-\bm{b}\|_1-\sigma, \,0\big{\}}.
\end{equation}
Then, based on these quantities, we set
\begin{equation*}
\theta=1/\rho \quad \mathrm{and} \quad
\rho=\left\{\begin{array}{ll}
1.2, ~&\mathrm{if}~\max\left\{\eta_1^k, \,\eta_2^k, \,\eta_3^k\right\}<10^{-2}, \\
2,   ~&\mathrm{otherwise}.
\end{array}\right.
\end{equation*}
The initial tolerance for the subproblem is set to $\epsilon_0=10^{-3}$ and $\epsilon_{k+1}$ is updated as $\max\{\theta\epsilon_k, 10^{-8}\}$ (instead of $\theta\epsilon_k$) in our implementation. Finally, we terminate the SPeL1 when
\begin{equation*}
\max\left\{\eta_1^k, \,\eta_2^k, \,\eta_3^k\right\}<10^{-8}.
\end{equation*}
Once the SPeL1 is terminated and returns an approximate solution $\bm{x}^*$, we also perform a \textit{refinement} step by setting $x_i^*=0$ if $|x_i^*|/\|\bm{x}^*\|_{\infty}<10^{-8}$ to improve the quality of the approximate solution.

For solving each subproblem \eqref{penprobsmooth} with $(\lambda_k,\mu_k,\nu_k)$ in the SPeL1, we adapt the NPG described in \eqref{npgiter} with $L_k^{\min}=10^{-6}$, $L_k^{\max}=\big{(}\frac{m}{\mu_k}+\frac{2}{\nu_k}\big{)}\lambda_k\|A\|^2$, $\tau=2$, $c=10^{-4}$ and $N=2$. Moreover, we set $L_{k,0}^0=1$ and, for any $l\geq1$,
\begin{equation*}
L_{k,l}^0 = \min\left\{\max\left\{\max\left\{\widetilde{\Delta}_k, \,0.5\bar{L}_{k,l-1}\right\}, \,L_k^{\min}\right\}, \,L_k^{\max}\right\}
\end{equation*}
with $\bm{x}^{k,-1}=\bm{x}^{k,0}$, where
\begin{equation*}
\begin{aligned}
\widetilde{\Delta}_k & := \frac{\Delta_k(\bm{x}^{k,l}, \bm{x}^{k,l-1})+\Delta_k(\bm{x}^{k,l}, \bm{x}^{k,l-2})+\Delta_k(\bm{x}^{k,l-1}, \bm{x}^{k,l-2})}{3}, \\[3pt]
\Delta_k(\bm{y}, \tilde{\bm{y}}) &:=
\left\{\begin{aligned}
&\frac{\langle\,\bm{y}-\tilde{\bm{y}}, \,\nabla f_{\lambda_k,\mu_k,\nu_k}(\bm{y}) - \nabla f_{\lambda_k,\mu_k,\nu_k}(\tilde{\bm{y}})\,\rangle}{\|\bm{y}-\tilde{\bm{y}}\|^2}, && \mathrm{if}~~\bm{y}\neq\tilde{\bm{y}}, \\
&0, && \mathrm{otherwise}.
\end{aligned}\right.
\end{aligned}
\end{equation*}
The NPG method is terminated when the number of iterations exceeds 1000 \textit{or}
\begin{equation*}
\frac{\bar{L}_{k,l}\|\bm{x}^{k,l+1}-\bm{x}^{k,l}\|}{1+\|\bm{x}^{k,l+1}\|} < \epsilon_k
\quad \textit{or} \quad
\frac{\left|F_{\lambda_k,\mu_k,\nu_k}(\bm{x}^{k,l+1})-F_{\lambda_k,\mu_k,\nu_k}(\bm{x}^{k,l})\right|}
{1+\left|F_{\lambda_k,\mu_k,\nu_k}(\bm{x}^{k,l+1})\right|} < \epsilon_k^{1.2}.
\end{equation*}
Note from \eqref{optdistbd} that if the first inequality above holds, condition \eqref{suboptcond} is then approximately satisfied.

In the following experiments, we consider randomly generated instances. Given a dimensional triple $(m, n, s)$, we randomly generate an instance as follows. First, we generate a matrix $A \in \mathbb{R}^{m\times n}$ with i.i.d. standard Gaussian entries and then normalize $A$ so that each column of $A$ has unit norm. We next choose a subset $\mathcal{S}\subset\{1, \cdots, n\}$ of size $s$ uniformly at random and generate an $s$-sparse vector $\hat{\bm{x}}\in\mathbb{R}^{n}$, which has i.i.d. standard Gaussian entries on $\mathcal{S}$ and zeros on $\mathcal{S}^c$. Then, we generate the vector $\bm{b} \in \mathbb{R}^{m}$ by setting $\bm{b} = A\hat{\bm{x}}+\delta\bm{\xi}$, where $\delta>0$ is a scaling parameter and $\bm{\xi}\in\mathbb{R}^m$ is the noisy vector with each entry $\xi_i$ independently following certain distribution. We shall consider two cases:
\begin{itemize}[leftmargin=0.8cm]
\item \textbf{Case 1.} We use the standard Gaussian distribution via the {\sc Matlab} command: \texttt{xi = randn(m,1)}.
\item \textbf{Case 2.} We use the Student's $t(2)$ distribution via the {\sc Matlab} command: \texttt{xi = trnd(2,m,1)}.
\end{itemize}
Finally, we set $\sigma=\delta\|\bm{\xi}\|_1$ so that $\hat{\bm{x}}\in{\rm FEA}(A,\bm{b},\sigma,1)$. In particular, for such $\sigma$, we have observed from our simulations that all random instances satisfy $\|\bm{b}\|_1>\sigma$ and hence $0\notin{\rm FEA}(A,\bm{b},\sigma,1)$.

Table \ref{ResSPeL1} presents the numerical results of the SPeL1 for solving problem \eqref{lpl1model} with $0<p<1$, where we use $\delta=10^{-3}$ and consider different choices of $(m, n, s)$ and $p$ under different noisy cases. In this table, ``\textbf{nnz}" denotes the number of nonzero entries in the refined terminating solution $\bm{x}^*$; ``\textbf{rank}" denotes the rank of $A_{\mathcal{J}}$ with ${\mathcal J}=\mathrm{supp}(\bm{x}^*)$; $\mathbf{err}_1 := \max\big{\{}\|\bm{x}^*\|_\infty - (\lambda_{\min}(A_{\mathcal{J}}^{\top}A_{\mathcal{J}}))^{-\frac{1}{2}}(\sigma+\|\bm{b}\|_2), \, (|\mathcal{J}|\lambda_{\max}(A_{\mathcal{J}}^{\top}A_{\mathcal{J}}))^{-\frac{1}{2}}(\|\bm{b}\|_1-\sigma)-\|\bm{x}^*\|_{\infty}, \,0\big{\}}$; and $\mathbf{err}_2 := \sigma - \|A\bm{x}^*-\bm{b}\|_1$. All results presented are the average of 10 independent instances for each $(m, n, s)$ and we display the rounding numbers for ``\textbf{nnz}" and ``\textbf{rank}". From Table \ref{ResSPeL1}, one can see that \textbf{nnz} $=$ \textbf{rank}, $\mathbf{err}_1=0$ and $\mathbf{err}_2\approx0$ always hold, clearly matching Theorem \ref{thmeqvl0} established for an optimal solution of problem \eqref{lpl1model} with $0<p<1$. This implies that our SPeL1 is able to find a `good' stationary point of problem \eqref{lpl1model} with $0<p<1$, which has important properties of an optimal solution.

\begin{figure}[ht]
\centering
\includegraphics[width=7.5cm]{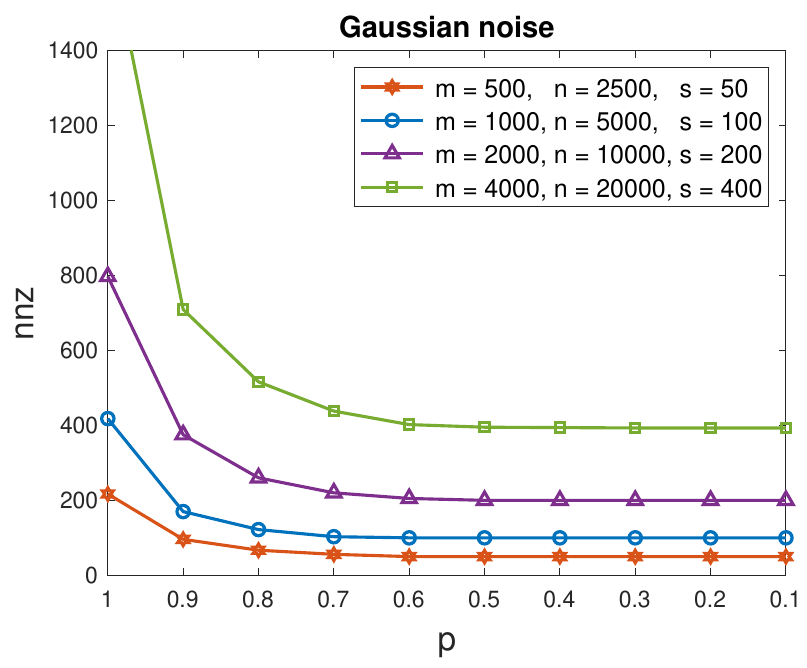}~~~~
\includegraphics[width=7.5cm]{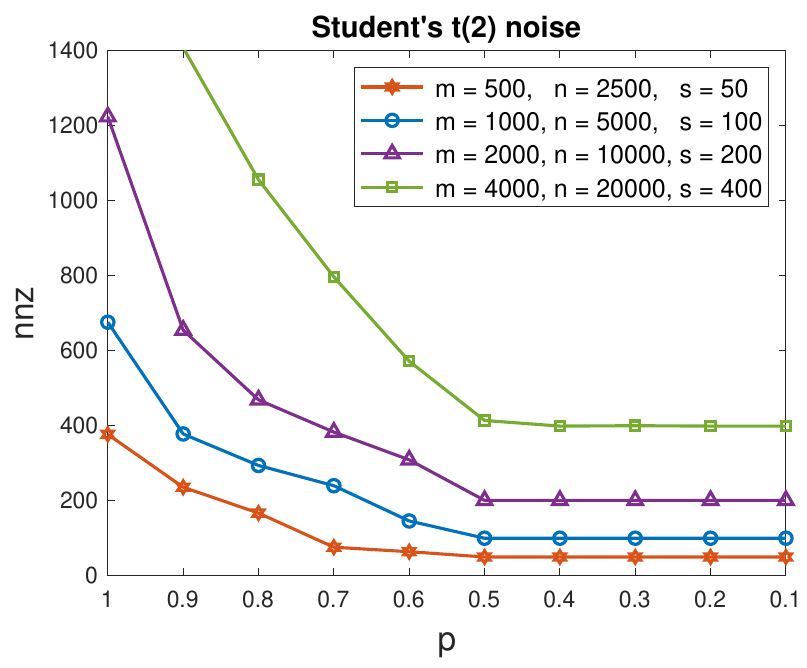}
\caption{The number of nonzero entries in the approximate solution for different $p$.}\label{figdiffp}
\end{figure}

We further generate one random instance for each $(m, n, s)$ under different noisy cases, and then apply our SPeL1 to solve problem \eqref{lpl1model} with different $p$. The number of nonzero entries in the approximate solution obtained for different $p$ are presented in Figure \ref{figdiffp}. From this figure, we see that solving problem \eqref{lpl1model} with a smaller $p$ always gives a sparser approximate solution, and the sparsity is almost unchanged and is close to the sparsity of $\hat{\bm{x}}$ when $p$ is smaller than a certain threshold. This observation implicitly matches Theorem \ref{mianthm1}, which says that ${\rm SOL}(A,\bm{b},\sigma,p,1)\subseteq {\rm SOL}(A,\bm{b},\sigma,0,1)$ and ${\rm SOL}(A,\bm{b},\sigma,p,1)$ remains unchanged for any sufficiently small $p$, and shows the potential advantage of solving problem \eqref{lpl1model} with a small $p$ for finding a sparse solution. Moreover, in practice, such $p$ may not be necessarily too small. From our experiments, we observe that $p=0.5$ is small enough for problem \eqref{lpl1model} to give a sparse solution.

\begin{table}[ht]
\TABLE
{Numerical results of the SPeL1 for solving \eqref{lpl1model} with $0<p<1$.\label{ResSPeL1}}
{\begin{tabular}{|ccc|c|cccc|cccc|}
\hline
\multicolumn{4}{|c|}{} & \multicolumn{4}{c|}{\footnotesize{Gaussian noise}} & \multicolumn{4}{c|}{\footnotesize{Student's $t(2)$ noise}} \\
\hline
$m$ & $n$ & $s$ & $p$ & \textbf{nnz} & \textbf{rank} & $\mathbf{err}_1$ & $\mathbf{err}_2$ & \textbf{nnz} & \textbf{rank} & $\mathbf{err}_1$ & $\mathbf{err}_2$  \\
\hline
\multirow{5}{*}{500} & \multirow{5}{*}{2500} & \multirow{5}{*}{50}
&   0.9 &  92 &  92 & 0 & 2.54e-7 & 163 & 163 & 0 & 1.82e-7 \\
&&& 0.7 &  53 &  53 & 0 & 1.15e-7 &  87 &  87 & 0 & 6.65e-8 \\
&&& 0.5 &  50 &  50 & 0 & 2.24e-7 &  50 &  50 & 0 & 2.17e-7 \\
&&& 0.3 &  50 &  50 & 0 & 4.02e-7 &  50 &  50 & 0 & 2.47e-7 \\
&&& 0.1 &  50 &  50 & 0 & 3.43e-7 &  50 &  50 & 0 & 3.00e-7 \\
\hline
\multirow{5}{*}{1000} & \multirow{5}{*}{5000} & \multirow{5}{*}{100}
&   0.9 & 164 & 164 & 0 & 4.71e-7 & 366 & 366 & 0 & 3.69e-7 \\
&&& 0.7 & 105 & 105 & 0 & 1.87e-7 & 210 & 210 & 0 & 1.35e-7 \\
&&& 0.5 &  99 &  99 & 0 & 3.53e-7 & 100 & 100 & 0 & 4.00e-7 \\
&&& 0.3 &  99 &  99 & 0 & 4.08e-7 &  99 &  99 & 0 & 4.43e-7 \\
&&& 0.1 &  99 &  99 & 0 & 6.97e-7 &  99 &  99 & 0 & 5.11e-7 \\
\hline
\multirow{5}{*}{2000} & \multirow{5}{*}{10000} & \multirow{5}{*}{200}
&   0.9 & 337 & 337 & 0 & 8.03e-7 & 706 & 706 & 0 & 6.04e-7 \\
&&& 0.7 & 214 & 214 & 0 & 3.72e-7 & 426 & 426 & 0 & 2.34e-7 \\
&&& 0.5 & 199 & 199 & 0 & 4.90e-7 & 199 & 199 & 0 & 3.97e-7 \\
&&& 0.3 & 198 & 198 & 0 & 6.10e-7 & 198 & 198 & 0 & 6.71e-7 \\
&&& 0.1 & 198 & 198 & 0 & 7.01e-7 & 198 & 198 & 0 & 6.86e-7 \\
\hline
\multirow{5}{*}{4000} & \multirow{5}{*}{20000} & \multirow{5}{*}{400}
&   0.9 & 703 & 703 & 0 & 1.01e-6 & 1611 & 1611 & 0 & 9.52e-7 \\
&&& 0.7 & 433 & 433 & 0 & 5.24e-7 & 873 & 873 & 0   & 3.34e-7 \\
&&& 0.5 & 398 & 398 & 0 & 6.72e-7 & 418 & 418 & 0   & 6.52e-7 \\
&&& 0.3 & 397 & 397 & 0 & 6.79e-7 & 396 & 396 & 0   & 1.16e-6 \\
&&& 0.1 & 396 & 396 & 0 & 8.68e-7 & 396 & 396 & 0   & 1.16e-6 \\
\hline
\end{tabular}}
{}
\end{table}

Next, we consider using model \eqref{lpl1model} to recover a sparse solution of an underdetermined linear system from noisy measurements, and compare its performance with that of using the widely-studied $L_2$-constrained problem (see, for example, \cite{bde2009sparse,clp2016penalty,cw2018spherical,vf2008probing}):
\begin{equation}\label{lpl2model}
\min\limits_{\bm{x}\in\mathbb{R}^n}~~\|\bm{x}\|_p^p \quad~~ \mbox{\rm s.t.} \quad~~ \|A\bm{x}-\bm{b}\|_2 \leq \sigma.
\end{equation}
We will solve problem \eqref{lpl2model} with $0<p<1$ by the smoothing penalty method\footnote{The {\sc Matlab} codes implemented by the authors in \cite{clp2016penalty} are available at \url{http://www.mypolyuweb.hk/~tkpong/Exact_lp_codes/}} proposed in \cite{clp2016penalty} and call it SPeL2 for short. All parameters in the SPeL2 are chosen as the default settings, except that we terminate its subroutine NPG when the inner iteration number exceeds 1000 to save the cost for solving the subproblem, while maintaining the quality of the eventual solution. Moreover, we initialize the SPeL2 at the same point as the SPeL1 and terminate the SPeL2 at the $k$th iteration when $\max\left\{\eta_1^k, \,\eta_2^k, \,\eta_4^k\right\}<10^{-8}$, where $\eta_1^k$, $\eta_2^k$ are defined in \eqref{defeta} and $\eta_4^k := \max\big{\{}\|A\bm{x}^{k+1}-\bm{b}\|_2-\sigma, \,0\big{\}}$. We also adapt the refinement step for the approximate solution obtained by the SPeL2 to improve the quality of the approximate solution.

In comparisons below, we use $p=0.5$ and consider different $(m,n,s)$ and $\delta$ under different noisy cases. For each $(m,n,s)$ and $\delta$, we randomly generate $A$, $\hat{\bm{x}}$, $\bm{b}$, $\bm{\xi}$ as described above, but set $\sigma=\delta\|\bm{\xi}\|_1$ for \eqref{lpl1model} and set $\sigma=\delta\|\bm{\xi}\|$ for \eqref{lpl2model} so that both resulting feasible sets of \eqref{lpl1model} and \eqref{lpl2model} will contain the sparse vector $\hat{\bm{x}}$ as a boundary point. The computational results are reported in Table \ref{ResComp1}, where ``\textbf{nnz}" denotes the number of nonzero entries in the refined terminating solution $\bm{x}^*$; ``\textbf{feas}" denotes the deviation of $\bm{x}^*$ from the constraint, which is given by $\eta_3^k$ for \eqref{lpl1model} and $\eta_4^k$ for \eqref{lpl2model}; ``\textbf{recerr}" denotes the relative recovery error $\|\bm{x}^*-\hat{\bm{x}}\|_2/\|\hat{\bm{x}}\|_2$; ``\textbf{time}" denotes the computational time (in seconds). All results reported are the average of 10 independent instances for each $(m, n, s)$ and $\delta$. One can observe from this table that for the Gaussian noisy case, the performance of our SPeL1 is comparable with that of the SPeL2 with respect to the relative recovery error, while for the Student's $t(2)$ noisy case, our SPeL1 gives sparse solutions with smaller relative recovery errors for all instances. It is worth noting that, for the problem of recovering sparse solutions, even marginal improvements on recovery error could be very hard. Moreover, all approximate solutions obtained by the SPeL1 are exactly the feasible points of \eqref{lpl1model} and the sparsity of each solution is closer to that of the true sparse vector for most cases.

\begin{table}[ht]
\TABLE
{Comparisons between SPeL1 and SPeL2.\label{ResComp1}}
{\renewcommand\arraystretch{1.15}
\begin{tabular}{|c|ccc|c|cccc|cccc|}
\hline
\multicolumn{5}{|c|}{Problem Setting} & \multicolumn{4}{c|}{\footnotesize{SPeL1}} & \multicolumn{4}{c|}{\footnotesize{SPeL2}} \\
\hline
\textbf{noise}  &$m$ & $n$ & $s$ & $\delta$ & \textbf{nnz} & \textbf{feas} & \textbf{recerr} & \textbf{time}
& \textbf{nnz} & \textbf{feas} & \textbf{recerr} & \textbf{time}  \\
\hline
\multirow{12}{*}{Gaussian}
&\multirow{3}{*}{500} & \multirow{3}{*}{2500} & \multirow{3}{*}{50}
&    $10^{-1}$ &  44 & 0 & 2.29e-1 & 0.41 &  37 & 1.07e-9 & 2.13e-1 & 1.88 \\
&&&& $10^{-2}$ &  49 & 0 & 1.86e-2 & 0.61 &  48 & 2.47e-9 & 1.82e-2 & 1.41 \\
&&&& $10^{-3}$ &  50 & 0 & 1.79e-3 & 0.71 &  50 & 3.22e-9 & 1.79e-3 & 0.99 \\
\cline{2-13}
&\multirow{3}{*}{1000} & \multirow{3}{*}{5000} & \multirow{3}{*}{100}
&    $10^{-1}$ &  91 & 0 & 2.11e-1 & 2.79 &  73 & 7.02e-10 & 1.89e-1 & 16.31 \\
&&&& $10^{-2}$ &  97 & 0 & 1.59e-2 & 4.08 &  94 & 2.53e-9 & 1.51e-2 & 10.35 \\
&&&& $10^{-3}$ &  99 & 0 & 1.52e-3 & 5.11 &  99 & 3.32e-9 & 1.54e-3 & 10.30 \\
\cline{2-13}
&\multirow{3}{*}{2000} & \multirow{3}{*}{10000} & \multirow{3}{*}{200}
&    $10^{-1}$ & 184 & 0 & 1.94e-1 & 11.41 & 150 & 4.65e-10 & 1.73e-1 & 79.54 \\
&&&& $10^{-2}$ & 196 & 0 & 1.49e-2 & 16.46 & 190 & 2.46e-9  & 1.43e-2 & 51.62 \\
&&&& $10^{-3}$ & 199 & 0 & 1.44e-3 & 21.82 & 198 & 5.58e-9  & 1.43e-3 & 26.58 \\
\cline{2-13}
&\multirow{3}{*}{4000} & \multirow{3}{*}{20000} & \multirow{3}{*}{400}
&    $10^{-1}$ & 374 & 0 & 2.03e-1 & 46.68 & 294 & 7.15e-10 & 1.81e-1 & 438.38 \\
&&&& $10^{-2}$ & 398 & 0 & 1.54e-2 & 60.79 & 382 & 2.16e-9 & 1.47e-2 & 213.82 \\
&&&& $10^{-3}$ & 399 & 0 & 1.49e-3 & 99.13 & 397 & 5.24e-9 & 1.48e-3 & 140.70 \\
\hline
\multirow{12}{*}{Student's $t(2)$}
&\multirow{3}{*}{500} & \multirow{3}{*}{2500} & \multirow{3}{*}{50}
&    $10^{-1}$ &  45 & 0 & 3.61e-1 & 0.41 &  22 & 8.63e-10 & 6.51e-1 & 5.02 \\
&&&& $10^{-2}$ &  51 & 0 & 2.89e-2 & 0.71 &  45 & 9.02e-10 & 5.74e-2 & 2.12 \\
&&&& $10^{-3}$ &  50 & 0 & 2.43e-3 & 1.18 &  50 & 2.89e-9 & 5.42e-3 & 1.11 \\
\cline{2-13}
&\multirow{3}{*}{1000} & \multirow{3}{*}{5000} & \multirow{3}{*}{100}
&    $10^{-1}$ &  93 & 0 & 3.45e-1 & 2.64 &  44 & 4.02e-10 & 6.49e-1 & 31.04 \\
&&&& $10^{-2}$ & 104 & 0 & 2.59e-2 & 5.05 &  88 & 8.77e-10 & 5.23e-2 & 15.53 \\
&&&& $10^{-3}$ &  99 & 0 & 2.13e-3 & 9.24 &  97 & 2.67e-9  & 5.03e-3 & 11.38 \\
\cline{2-13}
&\multirow{3}{*}{2000} & \multirow{3}{*}{10000} & \multirow{3}{*}{200}
&    $10^{-1}$ & 196 & 0 & 3.33e-1 & 11.65 &  87 & 2.41e-10 & 6.68e-1 & 175.28 \\
&&&& $10^{-2}$ & 219 & 0 & 2.63e-2 & 22.64 & 178 & 3.65e-9 & 5.55e-2 & 50.26 \\
&&&& $10^{-3}$ & 206 & 0 & 2.23e-3 & 45.44 & 195 & 3.55e-9 & 5.47e-3 & 36.93 \\
\cline{2-13}
&\multirow{3}{*}{4000} & \multirow{3}{*}{20000} & \multirow{3}{*}{400}
&    $10^{-1}$ & 383 & 0 & 3.40e-1 & 48.74 & 179 & 3.53e-10 & 7.06e-1 & 1404.59 \\
&&&& $10^{-2}$ & 474 & 0 & 3.05e-2 & 84.90 & 351 & 2.57e-9 & 6.21e-2 & 384.40 \\
&&&& $10^{-3}$ & 410 & 0 & 2.27e-3 & 225.96& 390 & 5.17e-9 & 5.93e-3 & 177.51 \\
\hline
\end{tabular}}
{}
\end{table}

\begin{figure}[ht]
\centering
\subfigure[~Gaussian noise]{
\includegraphics[width=7.5cm]{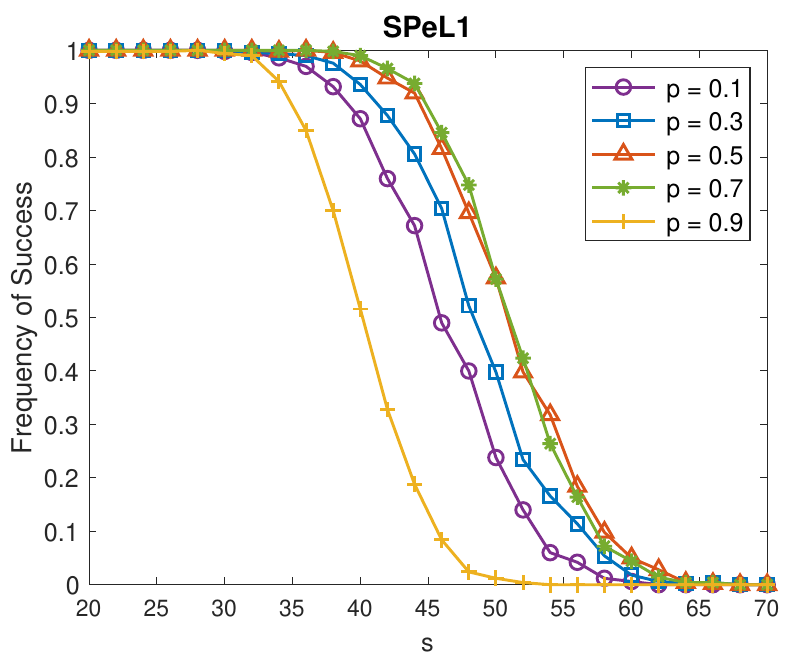}~~~~
\includegraphics[width=7.5cm]{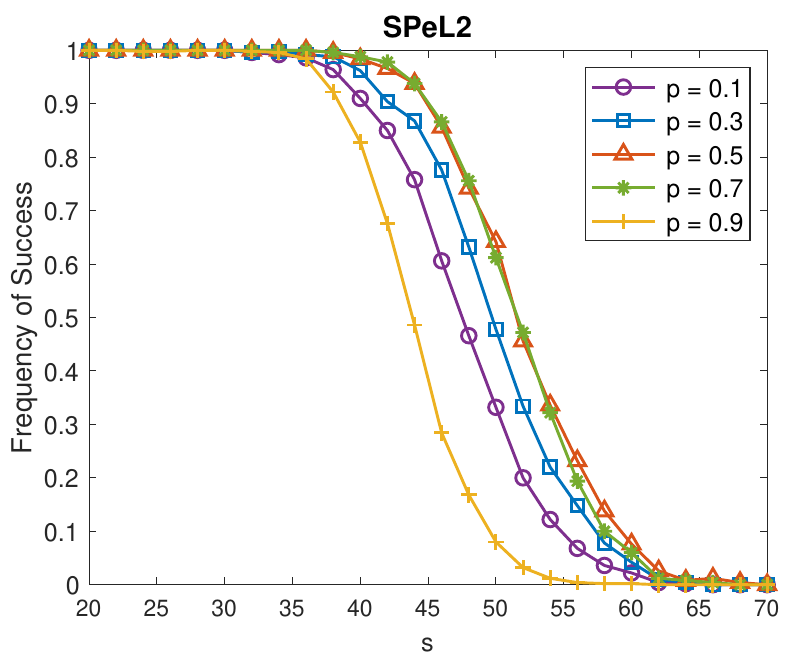}}
\subfigure[~Student's $t(2)$ noise]{
\includegraphics[width=7.5cm]{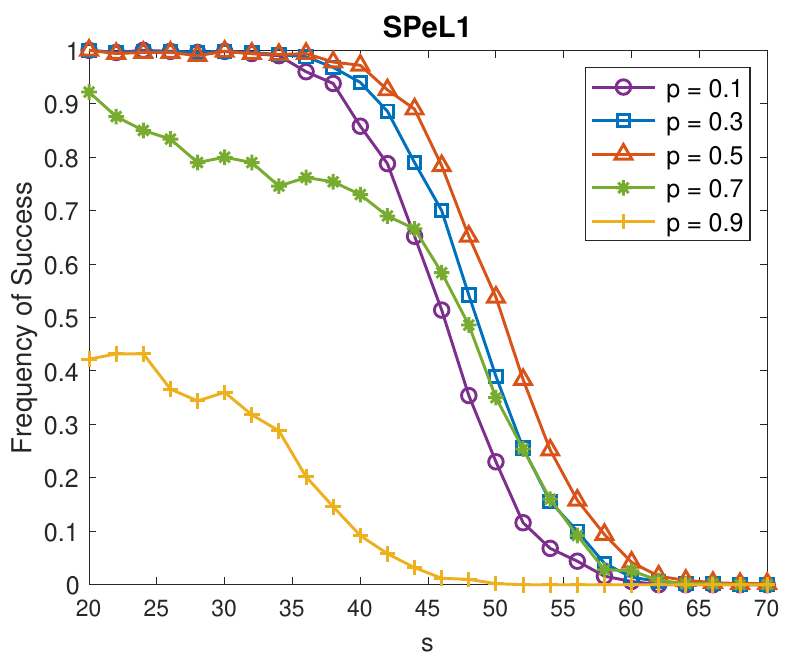}~~~~
\includegraphics[width=7.5cm]{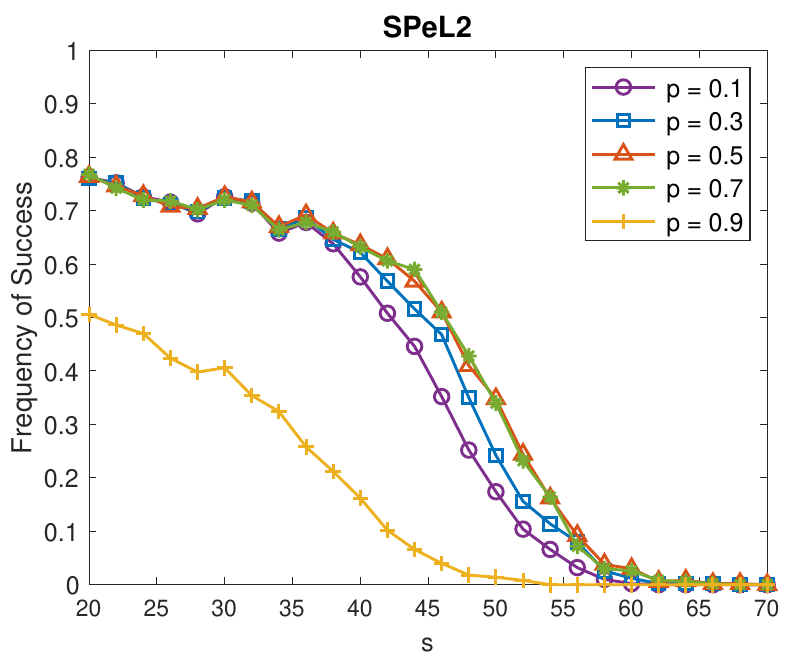}}
\caption{Comparisons between SPeL1 and SPeL2 with $m=128$, $n=512$ and different $s$.}\label{figRec}
\end{figure}

To better visualize the recovery performances of SPeL1 and SPeL2, we generate more instances to test and plot the ``frequency of success" for each method with different $p$. Specifically, we fix $m=128$, $n=512$ and vary $s$ from 20 to 70. The noisy level is set to $\delta=10^{-3}$. For each $(m,n,s)$, we generate 500 independent instances, and for each instance, we run each method to obtain an approximate solution $\bm{x}^*$ and consider the recovery successful if $\|\bm{x}^*-\hat{\bm{x}}\|_2/\|\hat{\bm{x}}\|_2 < 5\times 10^{-3}$. The results of the experiments are presented in Figure \ref{figRec}. Note that when the number of measurements is fixed, a larger $s$ generally leads to a more difficult recovery problem and thus the successful rate would be decayed, as shown in the figure. Moreover, one can see that for the Gaussian noisy case, the successful rate of our SPeL1 is comparable with that of the SPeL2, while for the Student's t(2) noisy case, our SPeL1 can give better successful rates especially when $p$ is small. This highlights the potential advantage of our approach for recovering a sparse solution under non-Gaussian noisy cases. One may also observe that when $s$ becomes larger and $p\leq0.5$, the successful rates of both methods appear to become lower as $p$ becomes smaller. The possible reason is that when $s$ is large and $p$ is too small, finding a solution of problem \eqref{lpl1model} or \eqref{lpl2model} can be rather difficult and hence it is less likely for a stationary point to be a good candidate. Therefore, both SPeL1 and SPeL2 may still need some improvements for the hard cases ($p$ is small and $s$ is large). We will leave this interesting research topic in the future.

\section{Concluding remarks}\label{secconc}

In this paper, we consider a unified $L_p$-$L_q$ sparse optimization problem \eqref{lplqmodel} and study various properties of its optimal solutions. Specifically, without any condition on the sensing matrix $A$, we provide upper bounds in cardinality and infinity norm for the optimal solutions, and show that all optimal solutions must be at the boundary of the feasible set when $0<p\leq 1$; see Theorem \ref{thmeqvl0}. Moreover, for $q \in \{1,\infty\}$, we show that the $L_q$-constrained problem with $0<p<1$ has finitely many optimal solutions; see Proposition \ref{extrp} and Remark \ref{addremark1}. We further show that, for $q \in \{1,\infty\}$, there exists $0<p^*<1$ such that the solution set of the problem with any $0<p<p^*$ is contained in the solution set of the problem with $p=0$ and there also exists $0<\overline{p}<p^*$ such that the solution set of the problem with any $0<p\leq\overline{p}$ remains unchanged; see Theorem \ref{mianthm1} and Remark \ref{addremark2}. An estimation of such $p^*$ is also provided in Theorem \ref{pstarthm}. A convergent smoothing penalty method is also proposed to solve the $L_1$-constrained problem with $0<p<1$. Some numerical examples are presented to implicitly illustrate the theoretical results and show the efficiency of the proposed method for solving the constrained $L_p$-$L_1$ problem under different noises.

%%%%%%%%%%%%%%%%%%%%%%%%%%%%%%%%%%%%%%%%%%%%%%%%%%%%%%
% Appendix here
% Options are (1) APPENDIX (with or without general title) or
%             (2) APPENDICES (if it has more than one unrelated sections)
% Outcomment the appropriate case if necessary
%
% \begin{APPENDIX}{<Title of the Appendix>}
% \end{APPENDIX}
%
% or

\begin{APPENDICES}

\section{Proof of Lemma \ref{inequ}}\label{prooflemmas1}

First, for $k = 1,\cdots,n$, we define $p_k(\bm{a}) := {\textstyle\sum_{j=1}^n}a_j^k$, $p_k(\bm{b}) := {\textstyle\sum_{j=1}^n}b_j^k$,
\begin{equation*}
\Lambda_k(\bm{a}) := \sum_{1\le i_1<i_2<\cdots<i_k\le n}a_{i_1}a_{i_2}\cdots  a_{i_k}
\quad \mathrm{and} \quad
\Lambda_k(\bm{b}) := \sum_{1\le i_1<i_2<\cdots<i_k\le n}b_{i_1}b_{i_2}\cdots  b_{i_k}.
\end{equation*}
Then, from Vi\`{e}te's formula \cite{Vieta}, we see that $a_1,\cdots,a_n$ and $b_1,\cdots,b_n$ are the roots of $q_n(t)$ and $r_n(t)$, respectively, where
\begin{equation*}
\begin{aligned}
q_n(t) &:= t^n - \Lambda_1(\bm{a})t^{n-1} + \Lambda_2(\bm{a})t^{n-2} + \cdots + (-1)^{n-1}\Lambda_{n-1}(\bm{a})t^{1} + (-1)^n\Lambda_{n}(\bm{a}) = 0,   \\
r_n(t) &:= t^n - \Lambda_1(\bm{b})t^{n-1} + \Lambda_2(\bm{b})t^{n-2} + \cdots + (-1)^{n-1}\Lambda_{n-1}(\bm{b})t^{1} + (-1)^{n}\Lambda_{n}(\bm{b}) = 0.
\end{aligned}
\end{equation*}
Moreover, from \cite[Eq.\,($2.11'$)]{MacDonald} and the discussions that follow, we have that, for $k=1,\cdots,n$,
\begin{equation}\label{Lam}
k\Lambda_k(\bm{a}) = {\textstyle\sum_{j=1}^k}(-1)^{j-1}p_{j}(\bm{a})\Lambda_{k-j}(\bm{a}),
\quad
k\Lambda_k(\bm{b}) = {\textstyle\sum_{j=1}^k}(-1)^{j-1}p_{j}(\bm{b})\Lambda_{k-j}(\bm{b})
\end{equation}
with $\Lambda_0(\bm{a})=\Lambda_0(\bm{b})=1$. Notice that $\Lambda_1(\bm{a})=\Lambda_1(\bm{b})=p_1(\bm{a})=p_1(\bm{b})$ and $p_k(\bm{a})=p_k(\bm{b})$ for $k=1,\cdots,n$. Thus, from \eqref{Lam}, it is not hard to show by induction that $\Lambda_k(\bm{a}) = \Lambda_k(\bm{b})$ holds for $k = 2, \cdots, n$. This implies that $q_n(t)$ and $r_n(t)$ have the same roots and hence $\bm{a} = \bm{b}$.
\endproof

\section{Proof of Lemma \ref{addlemma}}\label{prooflemmas2}

First, from the Taylor expansion (with Lagrange remainder), for any $0<p<1$, $c>0$ and $k\geq0$, we have
\begin{equation*}
c^p = e^{p\ln c} = 1 + p\ln c + \frac{(\ln c)^2}{2!}p^2 + \cdots + \frac{(\ln c)^k}{k!}p^k + \frac{e^{\xi_{k+1}}(\ln c)^{k+1}}{(k+1)!}p^{k+1},
\end{equation*}
where $\xi_{k+1}$ is a number between 0 and $p\ln c$. Then, for any $0<p<1$ and $k\geq0$, we have
\begin{equation}\label{TaylorE}
\begin{aligned}
\|\bm{a}\|_p^p &= \sum^s_{j=1} |a_{i_j}|^p = s + \sum_{l=1}^{k}\frac{\sum_{j=1}^{s}(\ln |a_{i_j}|)^l}{l!}p^l + \frac{\sum_{j=1}^{s}e^{\xi_{i_j,k+1}}(\ln |a_{i_j}|)^{k+1}}{(k+1)!}p^{k+1},  \\
\|\bm{b}\|_p^p &= \sum^s_{j=1} |b_{t_j}|^p = s + \sum_{l=1}^{k}\frac{\sum_{j=1}^{s}(\ln |b_{t_j}|)^l}{l!}p^l + \frac{\sum_{j=1}^{s}e^{\eta_{t_j,k+1}}(\ln |b_{t_j}|)^{k+1}}{(k+1)!}p^{k+1},
\end{aligned}
\end{equation}
where, for $j=1,\cdots,s$, $\xi_{i_j,k+1}$ is a number between 0 and $p\ln |a_{i_j}|$, and $\eta_{t_j,k+1}$ is a number between 0 and $p\ln |b_{t_j}|$. In the following, we consider two cases.

\textbf{Case 1}: $\Delta_k(\bm{a},\bm{b})=0$ for all $k=1,\cdots,s$, where $\Delta_k(\bm{a},\bm{b})$ is defined as \eqref{defDeltak}. In this case, we have $\sum_{j=1}^{s}(\ln |a_{i_j}|)^k = \sum_{j=1}^{s}(\ln |b_{t_j}|)^k$ for all $k=1,\cdots,s$. This together with Lemma \ref{inequ} further implies that $(\ln|a_{i_1}|,\cdots,\ln|a_{i_s}|) = (\ln|b_{t_1}|, \cdots, \ln|b_{t_s}|)$ and hence $(|a_{i_1}|,\cdots,|a_{i_s}|) = (|b_{t_1}|, \cdots, |b_{t_s}|)$. Then, we have $\|\bm{a}\|_p^p = \|\bm{b}\|_p^p$ for any $p>0$. This proves statement (i).

\textbf{Case 2}: Case 1 does not hold. In this case, there must exist some $1\leq\tilde{k}\leq s$ so that $\Delta_{\tilde{k}}(\bm{a}, \bm{b})\neq0$ and $\Delta_k(\bm{a}, \bm{b})=0$ for $k = 1, \cdots, \tilde{k}-1$. Then, we have from \eqref{TaylorE} and \eqref{defDeltak} that
\begin{equation}\label{diffab}
\|\bm{a}\|_{p}^{p} - \|\bm{b}\|_{p}^{p}
= \frac{p^{\tilde{k}}}{{\tilde{k}}!}\left(\Delta_{\tilde{k}}(\bm{a}, \bm{b})+\frac{p}{\tilde{k}+1}\Xi_{\tilde{k}+1}^p(\bm{a}, \bm{b})\right),
\end{equation}
where $\Xi_{\tilde{k}+1}^p(\bm{a}, \bm{b}) := {\textstyle\sum_{j=1}^{s}}\big{(}e^{\xi_{i_j,\tilde{k}+1}}(\ln |a_{i_j}|)^{\tilde{k}+1}-e^{\eta_{t_j,\tilde{k}+1}}(\ln |b_{t_j}|)^{\tilde{k}+1}\big{)}$. Note also that $\Delta_{\tilde{k}}(\bm{a}, \bm{b})\neq0$ and $\frac{p}{\tilde{k}+1}\Xi_{\tilde{k}+1}^p(\bm{a}, \bm{b})\to0$ as $p\to0$. Thus, there must exist a sufficiently small $p'$ such that
\begin{equation}\label{ineq-pdb}
{\textstyle \left|\frac{p}{\tilde{k}+1}\Xi_{\tilde{k}+1}^p(\bm{a}, \bm{b})\right|} \leq \frac{1}{2}|\Delta_{\tilde{k}}(\bm{a}, \bm{b})|, \quad \forall\,p\in(0,\,p'].
\end{equation}
We now consider the following two cases.
\begin{itemize}
\item $\Delta_{\tilde{k}}(\bm{a}, \bm{b})<0$: in this case, using \eqref{diffab} and \eqref{ineq-pdb}, we obtain that
    \begin{equation*}
    \|\bm{a}\|_{p}^{p} - \|\bm{b}\|_{p}^{p} \leq \frac{p^{\tilde{k}}}{{\tilde{k}}!}\left(\Delta_{\tilde{k}}(\bm{a}, \bm{b}) + \frac{1}{2}|\Delta_{\tilde{k}}(\bm{a}, \bm{b})|\right) = \frac{p^{\tilde{k}}}{2\tilde{k}!}\Delta_{\tilde{k}}(\bm{a}, \bm{b}) < 0, \quad \forall\,p\in(0,\,p'].
    \end{equation*}
    This implies that $\|\bm{a}\|_p^p<\|\bm{b}\|_p^p$ for any $p\in(0,\,p']$.

\item $\Delta_{\tilde{k}}(\bm{a}, \bm{b})>0$: in this case, using \eqref{diffab} and \eqref{ineq-pdb}, we obtain that
    \begin{equation*}
    \|\bm{a}\|_{p}^{p} - \|\bm{b}\|_{p}^{p} \geq \frac{p^{\tilde{k}}}{{\tilde{k}}!}\left(\Delta_{\tilde{k}}(\bm{a}, \bm{b}) - \frac{1}{2}|\Delta_{\tilde{k}}(\bm{a}, \bm{b})|\right) = \frac{p^{\tilde{k}}}{2\tilde{k}!}\Delta_{\tilde{k}}(\bm{a}, \bm{b}) > 0, \quad \forall\,p\in(0,\,p'].
    \end{equation*}
    This implies that $\|\bm{a}\|_p^p>\|\bm{b}\|_p^p$ for any $p\in(0,\,p']$.
\end{itemize}
Combing the above results, we complete the proof for statement (ii).
\endproof

\section{Exact penalization}\label{appexact}

In this section, we show that problem \eqref{penprob} is actually an exact penalization for problem \eqref{lpl1model} with $0<p<1$. For notational simplicity, we define a set $\mathcal{U}$ and a matrix $U$ as follows:
\begin{equation}\label{defu}
\mathcal{U}:=\big{\{}\bm{u}_1, \cdots, \bm{u}_{2^m}\big{\}} \quad \mathrm{and} \quad U:=[\bm{u}_1, \cdots, \bm{u}_{2^m}]^{\top}\in\mathbb{R}^{2^m \times m},
\end{equation}
where $\bm{u}_i\in\{-1,\,1\}^m$ and $\bm{u}_i\neq\bm{u}_j$ for any $i \neq j$. Since each entry of $\bm{u}_i$ is either $1$ \textit{or} $-1$ and the dimension of $\bm{u}_i$ is $m$, then one can have $2^m$ different choices of $\bm{u}_i$ and hence such $\mathcal{U}$ and $U$ are well-defined. Moreover, it is easy to see that if $\bm{u}_i \in \mathcal{U}$, then $-\bm{u}_i\in\mathcal{U}$. A simple example is given as follows: let $m=2$, then
\begin{equation*}
\mathcal{U}={\textstyle\left\{
\begin{bmatrix*}[r]1\\1\end{bmatrix*},\,
\begin{bmatrix*}[r]1\\-1\end{bmatrix*},\,
\begin{bmatrix*}[r]-1\\1\end{bmatrix*},\,
\begin{bmatrix*}[r]-1\\-1\end{bmatrix*}
\right\}}
\quad \mathrm{and} \quad
U=
\begin{bmatrix*}[r]
1  & 1  & -1 & -1 \\
1  & -1 & 1  & -1 \\
\end{bmatrix*}^{\top}.
\end{equation*}
We next present some auxiliary lemmas, which will be useful in our analysis.

\begin{lemma}\label{polyhd}
Let $A\in\mathbb{R}^{m\times n}$, $\bm{b}\in\mathbb{R}^m$ and $\sigma>0$. Then, ${\rm FEA}(A,\bm{b},\sigma,1)$ can be equivalently rewritten as $\{\bm{x}\in\mathbb{R}^n : UA \bm{x} \leq U\bm{b} + \sigma\mathbf{1}\}$, where $U$ is defined in \eqref{defu} and $\mathbf{1}:=(1,\cdots,1)^{\top}\in\mathbb{R}^{2^m}$.
\end{lemma}
\beginproof
Observe that
\begin{equation*}
\begin{aligned}
&~~\big{\{}\bm{x}\in\mathbb{R}^n : \|A\bm{x} - \bm{b}\|_1 \leq \sigma\big{\}}
= \big{\{}\bm{x}\in\mathbb{R}^n : \max\limits_{\|\bm{u}\|_{\infty} \leq 1} \langle \bm{u}, \,A\bm{x}-\bm{b} \rangle \leq \sigma\big{\}}  \\
&= \big{\{}\bm{x}\in\mathbb{R}^n : \max\limits_{\bm{u}\in\mathcal{U}}\,\langle \bm{u}, \,A\bm{x}-\bm{b} \rangle \leq \sigma\big{\}}
= \big{\{}\bm{x}\in\mathbb{R}^n : \bm{u}_i^{\top}(A\bm{x}-\bm{b}) \leq \sigma, ~\bm{u}_i\in\mathcal{U}, ~i = 1, \cdots, 2^m\big{\}}  \\
&= \big{\{}\bm{x}\in\mathbb{R}^n : U(A\bm{x}-\bm{b}) \leq \sigma\mathbf{1} \big{\}}
= \big{\{}\bm{x}\in\mathbb{R}^n : UA\bm{x} \leq U\bm{b} + \sigma\mathbf{1} \big{\}},
\end{aligned}
\end{equation*}
where the first equality follows from $\|A\bm{x} - \bm{b}\|_1 = \max\limits_{\|\bm{u}\|_{\infty} \leq 1} \langle \bm{u}, \,A\bm{x}-\bm{b} \rangle$, the second equality follows because the maximizer of $\max\limits_{\|\bm{u}\|_{\infty} \leq 1} \langle \bm{u}, \,A\bm{x}-\bm{b} \rangle$ must be an extreme point of $\{\bm{u}:\|\bm{u}\|_{\infty}\leq1\}$ (see \cite[Corollary 32.3.4]{r1970convex}) and $\mathcal{U}$ is the set of all extreme points of $\{\bm{u}:\|\bm{u}\|_{\infty}\leq1\}$. This completes the proof.
\endproof

From Lemma \ref{polyhd}, it is easy to see that the feasible set ${\rm FEA}(A,\bm{b},\sigma,1)$ is a convex polyhedron. This together with the Hoffman error bound theorem \cite{h1952on} gives the following lemma.

\begin{lemma}\label{hoffmanlem}
There exists a constant $\tilde{c} > 0$ such that
\begin{equation*}
\mathrm{dist}\left(\bm{x}, \,{\rm FEA}(A,\bm{b},\sigma,1)\right) \leq \tilde{c} \,\|(\tilde{A} \bm{x} - \tilde{\bm{b}})_+\|_1
\end{equation*}
holds for any $\bm{x}\in\mathbb{R}^n$, where $\tilde{A}=UA$, $\tilde{\bm{b}}=U\bm{b}+\sigma\mathbf{1}$ and $U$ is defined in \eqref{defu}.
\end{lemma}

Based on this error bound result, we further give the following lemma.
\begin{lemma}\label{errorbdlem}
There exists a constant $c > 0$ such that, for any $\bm{x}\in\mathbb{R}^n$, we have
\begin{equation*}
\mathrm{dist}\left(\bm{x}, \,{\rm FEA}(A,\bm{b},\sigma,1)\right) \leq c \,(\|A\bm{x}-\bm{b}\|_1 - \sigma)_+.
\end{equation*}
\end{lemma}
\beginproof
We first show that, for any $\bm{x}\in\mathbb{R}^n$, it holds that
\begin{equation}\label{eqresidual}
2^{1-m}\|(\tilde{A} \bm{x} - \tilde{\bm{b}})_+\|_1 \leq (\|A\bm{x}-\bm{b}\|_1 - \sigma)_+ \leq \|(\tilde{A} \bm{x} - \tilde{\bm{b}})_+\|_1,
\end{equation}
where $\tilde{A}$ and $\tilde{\bm{b}}$ are defined in Lemma \ref{hoffmanlem}. Indeed, for any $\bm{x}\in\mathbb{R}^n$, there exists some $\tilde{\bm{u}} \in \{-1,\,1\}^m$ such that $\|A\bm{x}-\bm{b}\|_1=\tilde{\bm{u}}^{\top}(A\bm{x}-\bm{b})$. Then, we have
\begin{equation*}
\begin{aligned}
&~~\|(\tilde{A} \bm{x} - \tilde{\bm{b}})_+\|_1
= \|(UA\bm{x} - U\bm{b} - \sigma\mathbf{1})_+\|_1 = \|(U(A\bm{x} - \bm{b}) - \sigma\mathbf{1})_+\|_1 \\
&= {\textstyle\sum_{\bm{u}_i\in\mathcal{U}}} (\bm{u}_i^{\top}(A\bm{x}-\bm{b}) - \sigma)_+
\geq (\tilde{\bm{u}}^{\top}(A\bm{x}-\bm{b})-\sigma)_+
= (\|A\bm{x}-\bm{b}\|_1-\sigma)_+.
\end{aligned}
\end{equation*}
On the other hand, from $\|A\bm{x} - \bm{b}\|_1 = \max\limits_{\|\bm{u}\|_{\infty} \leq 1} \langle \bm{u}, \,A\bm{x}-\bm{b} \rangle$, we have
\begin{equation}\label{ineqadd2}
\|A\bm{x} - \bm{b}\|_1 \geq \bm{u}_i^{\top}(A\bm{x}-\bm{b}), \quad  \forall\,\bm{u}_i\in\mathcal{U}, ~i = 1, \cdots, 2^m.
\end{equation}
Then, we see that
\begin{equation*}
\begin{aligned}
\|(\tilde{A} \bm{x} - \tilde{\bm{b}})_+\|_1
&= \sum_{i=1}^{2^{m}} (\bm{u}_i^{\top}(A\bm{x}-\bm{b}) - \sigma)_+
= \sum^{|\mathcal{K}|=2^{m-1}}_{j\in\mathcal{K}\subset\{1,\cdots,2^{m}\}}(\bm{u}_j^{\top}(A\bm{x}-\bm{b})-\sigma)_+  \\
& \leq \sum^{|\mathcal{K}|=2^{m-1}}_{j\in\mathcal{K}\subset\{1,\cdots,2^{m}\}} (\|A\bm{x} - \bm{b}\|_1-\sigma)_+
= 2^{m-1} (\|A\bm{x} - \bm{b}\|_1-\sigma)_+,
\end{aligned}
\end{equation*}
where the second equality follows because if $\bm{u}_i\in\mathcal{U}$ and $\bm{u}_i^{\top}(A\bm{x}-\bm{b})-\sigma > 0$ , then $-\bm{u}_i\in\mathcal{U}$ and $-\bm{u}_i^{\top}(A\bm{x}-\bm{b})-\sigma < -\left(\bm{u}_i^{\top}(A\bm{x}-\bm{b})-\sigma\right) < 0$, and the inequality follows from \eqref{ineqadd2}. From the above, we obtain \eqref{eqresidual}. This together with Lemma \ref{hoffmanlem} completes the proof.
\endproof

Now, we are ready to present our exact penalization results. Our first theorem concerns local minimizers of problems \eqref{lpl1model} and \eqref{penprob}. The other two theorems concern $\epsilon$-minimizers of problems \eqref{lpl1model} and \eqref{penprob} (see definitions later).

\begin{theorem}\label{localthm}
Suppose that $\bm{x}^*$ is a local minimizer of \eqref{lpl1model}. Then, there exists a $\lambda^* > 0$ such that $\bm{x}^*$ is a local minimizer of \eqref{penprob} whenever $\lambda \geq \lambda^*$.
\end{theorem}
\beginproof
We first assume that $\bm{x}^*=0$ and consider any bounded neighborhood $\mathcal{N}$ of 0 and $\lambda > 0$. Let $L$ denote a Lipschitz constant of the function $\bm{x} \mapsto \lambda (\|A\bm{x}-\bm{b}\|_1 - \sigma)_+$ on $\mathcal{N}$. For this $L$, one can verify that there exists a neighborhood $\widetilde{\mathcal{N}} \subseteq \mathcal{N}$ of 0 such that $\|\bm{x}\|_p^p \geq L \|\bm{x}\|$ for all $\bm{x} \in \widetilde{\mathcal{N}}$. Then, for any $\bm{x} \in \widetilde{\mathcal{N}}$, we have
\begin{equation*}
F_{\lambda}(\bm{x})
= \|\bm{x}\|^p_p + \lambda (\|A\bm{x}-\bm{b}\|_1 - \sigma)_+ \geq L \|\bm{x}\| + \lambda (\|A\bm{x}-\bm{b}\|_1 - \sigma)_+
\geq \lambda (\|\bm{b}\|_1 - \sigma)_+ = F_{\lambda}(0),
\end{equation*}
where the last inequality follows from the definition of $L$ being a Lipschitz constant. This shows that $\bm{x}^*=0$ is a local minimizer of \eqref{penprob} for any $\lambda > 0$.

From now on, we assume that $\bm{x}^* \neq 0$. Let $\mathcal{J}:=\mathrm{supp}(\bm{x}^*)$ for simplicity. Then, $\mathcal{J} \neq \emptyset$ since $\bm{x}^* \neq 0$. Since $\bm{x}^*$ is a local minimizer of \eqref{lpl1model}, one can verify that $\bm{x}^*_{\mathcal{J}}$ is a local minimizer of the following problem:
\begin{equation}\label{subspacepro}
\min\limits_{\bm{x}_{\mathcal{J}}} ~ \|\bm{x}_{\mathcal{J}}\|_p^p  \qquad \mathrm{s.t.} \qquad \bm{x}_{\mathcal{J}} \in \Omega_{\mathcal{J}}:=\left\{\bm{x}_{\mathcal{J}}:\|A_{\mathcal{J}}\bm{x}_{\mathcal{J}} - \bm{b}\|_1 \leq \sigma\right\}.
\end{equation}
Let $\tilde{\epsilon} = \frac{1}{2}\min\big{\{}|x^*_i| : i \in \mathcal{J}\big{\}}>0$. Thus, there exists a small $\delta>0$ such that $\bm{x}^*_{\mathcal{J}}$ is a local minimizer of \eqref{subspacepro} and $\min\big{\{}|x_i|:i \in \mathcal{J}\big{\}} > \tilde{\epsilon}$ for all $\bm{x}_{\mathcal{J}} \in \mathcal{B}(\bm{x}^*_{\mathcal{J}}; \delta)$. Moreover, note that $\bm{x}_{\mathcal{J}}\mapsto\|\bm{x}_{\mathcal{J}}\|_p^p$ is Lipschitz continuous on $\mathcal{B}(\bm{x}^*_{\mathcal{J}}; \delta)$ and there exists a constant $c'>0$ such that $\mathrm{dist}(\bm{x}_{\mathcal{J}}, \,\Omega_{\mathcal{J}}) \leq c' \,(\|A_{\mathcal{J}}\bm{x}_{\mathcal{J}}-\bm{b}\|_1 - \sigma)_+$ for all $\bm{x}_{\mathcal{J}} \in \mathcal{B}(\bm{x}^*_{\mathcal{J}}; \delta)$ (see Lemma \ref{errorbdlem}). Then, from \cite[Lemma 3.1]{clp2016penalty} (or \cite[Proposition 4]{ldn2012exact}), there exists a $\lambda^*>0$ such that, for any $\lambda \geq \lambda^*$, $\bm{x}^*_{\mathcal{J}}$ is a local minimizer of the following problem:
\begin{equation*}
\min\limits_{\bm{x}_{\mathcal{J}}}~F^{\mathcal{J}}_{\lambda}(\bm{x}_{\mathcal{J}}):=\|\bm{x}_{\mathcal{J}}\|^p_p + \lambda (\|A_{\mathcal{J}}\bm{x}_{\mathcal{J}} - b\|_1 - \sigma)_+,
\end{equation*}
i.e., there exists a neighborhood $\mathcal{N}_{\mathcal{J}}$ of 0 with $\mathcal{N}_{\mathcal{J}} \subseteq \mathcal{B}(0; \frac{\delta}{2})$ such that
\begin{equation}\label{sublocalmin}
F^{\mathcal{J}}_{\lambda}(\bm{x}^*_{\mathcal{J}} + \bm{v}_{\mathcal{J}}) \geq F^{\mathcal{J}}_{\lambda}(\bm{x}^*_{\mathcal{J}}), \quad \forall\,\bm{v}_{\mathcal{J}} \in \mathcal{N}_{\mathcal{J}}.
\end{equation}

We now show that $\bm{x}^*$ is a local minimizer of \eqref{penprob} for any $\lambda \geq \lambda^*$. Fix any $\epsilon > 0$ and any $\lambda \geq \lambda^*$. Consider the bounded neighborhood $\mathcal{V}:=\mathcal{N}_{\mathcal{J}}\times(-\epsilon,\,\epsilon)^{n-|\mathcal{J}|}$ of 0 and let $\widetilde{L}$ be a Lipschitz constant of the function $g_{\lambda}(\bm{x}):=\lambda (\|A\bm{x}-\bm{b}\|_1 - \sigma)_+$ on $\bm{x}^* + \mathcal{V}$. For this $\widetilde{L}$, there exists an $\tilde{\epsilon}\in(0,\epsilon)$ such that $\|\bm{v}_{\mathcal{J}^{c}}\|_p^p \geq \widetilde{L} \|\bm{v}_{\mathcal{J}^{c}}\|$ for all $\bm{v}_{\mathcal{J}^{c}} \in (-\tilde{\epsilon},\,\tilde{\epsilon})^{n-|\mathcal{J}|}$. Then, for any $\bm{v} \in \widetilde{\mathcal{V}}:=\mathcal{N}_{\mathcal{J}}\times(-\tilde{\epsilon},\,\tilde{\epsilon})^{n-|\mathcal{J}|}$, we have
\begin{equation*}
\begin{aligned}
&~~F_{\lambda}(\bm{x}^* + \bm{v})
= \|\bm{x}^* + \bm{v}\|^p_p + g_{\lambda}(\bm{x}^* + \bm{v})
= \|\bm{x}^*_{\mathcal{J}} + \bm{v}_{\mathcal{J}}\|_p^p + \|\bm{v}_{\mathcal{J}^{c}}\|_p^p + g_{\lambda}(\bm{x}^* + \bm{v})  \\
&\geq \|\bm{x}^*_{\mathcal{J}} + \bm{v}_{\mathcal{J}}\|_p^p + \|\bm{v}_{\mathcal{J}^{c}}\|_p^p + g_{\lambda}\begin{pmatrix}\bm{x}_{\mathcal{J}}^* + \bm{v}_{\mathcal{J}}\\0\end{pmatrix}
- \widetilde{L}\|\bm{x}^*_{\mathcal{J}^c} + \bm{v}_{\mathcal{J}^{c}}\| \\
&\geq \|\bm{x}^*_{\mathcal{J}} + \bm{v}_{\mathcal{J}}\|_p^p + \widetilde{L}\|\bm{v}_{\mathcal{J}^{c}}\| + \lambda \left(\|A_{\mathcal{J}}(\bm{x}_{\mathcal{J}}^* + \bm{v}_{\mathcal{J}}) - \bm{b}\|_1 - \sigma\right)_+
- \widetilde{L}\|\bm{v}_{\mathcal{J}^{c}}\| \\
&= F^{\mathcal{J}}_{\lambda}(\bm{x}^*_{\mathcal{J}} + \bm{v}_{\mathcal{J}})
\geq F^{\mathcal{J}}_{\lambda}(\bm{x}^*_{\mathcal{J}}) = F_{\lambda}(\bm{x}^*_{\mathcal{J}}),
\end{aligned}
\end{equation*}
where the first inequality follows from the Lipschitz continuity of $g_{\lambda}$ with Lipschtiz constant $\widetilde{L}$ and the last inequality follows from \eqref{sublocalmin}. This shows that $\bm{x}^*$ is a local minimizer of \eqref{penprob} for any $\lambda \geq \lambda^*$ and completes the proof.
\endproof

We next study $\epsilon$-minimizers of \eqref{lpl1model} and \eqref{penprob}, which are defined as follows.

\begin{definition}[\textbf{$\epsilon$-minimizer}]
Let $\epsilon>0$.
\begin{itemize}[leftmargin=0.8cm]
\item[(i)] $\bm{x}_{\epsilon}$ is said to be an $\epsilon$-minimizer of problem \eqref{lpl1model} if $\bm{x}_{\epsilon}\in{\rm FEA}(A,\bm{b},\sigma,1)$ and $\|\bm{x}_{\epsilon}\|_p^p \leq \min\big{\{}\|\bm{x}\|_p^p : \bm{x}\in{\rm FEA}(A,\bm{b},\sigma,1)\big{\}}+\epsilon$.

\item[(ii)] $\bm{x}_{\epsilon}$ is said to be an $\epsilon$-minimizer of problem \eqref{penprob} if $F_{\lambda}(\bm{x}_{\epsilon}) \leq \min\limits_{\bm{x}\in\mathbb{R}^n}F_{\lambda}(\bm{x}) + \epsilon$.
\end{itemize}
\end{definition}

We also introduce the following function:
\begin{equation}\label{defPsimu}
\Psi_{\mu}(\bm{x}) = \sum^n_{i=1} \big{(}\psi_{\mu}(x_i)\big{)}^p \quad \mathrm{with} \quad
\psi_{\mu}(t) = \left\{\begin{aligned}
&|t|,  &&~~|t| \geq \mu, \\
&{\textstyle\frac{t^2}{2\mu} + \frac{\mu}{2}},  &&~~\mathrm{otherwise},
\end{aligned}\right.
\end{equation}
where $\mu>0$ is a constant. Note that $\Psi_{\mu}$ is continuously differentiable. Moreover, from the discussions in \cite[Section 3.3]{clp2016penalty}, we have that
\begin{eqnarray}
&&0 \leq \Psi_{\mu}(\bm{x}) - \|\bm{x}\|_p^p \leq n{\textstyle\left(\mu/2\right)}^p,  \label{approxbd} \\ [4pt]
&&|\Psi_{\mu}(\bm{x}) - \Psi_{\mu}(\bm{y})| \leq \sqrt{n}p\mu^{p-1} \|\bm{x} - \bm{y}\|.  \label{lippsi}
\end{eqnarray}
Then, we characterize the relation between the global minimizer of problem \eqref{lpl1model} and the $\epsilon$-minimizer of problem \eqref{penprob} in the next theorem.

\begin{theorem}\label{epsithm1}
Suppose that $\bm{x}^*$ is a global minimizer of problem \eqref{lpl1model}. Then, for any $\epsilon>0$, there exists a $\lambda_{\epsilon}^* > 0$ such that $\bm{x}^*$ is an $\epsilon$-minimizer of problem \eqref{penprob} whenever $\lambda \geq \lambda_{\epsilon}^*$.
\end{theorem}
\beginproof
First, for any $\epsilon>0$, we consider $\mu=2\left(\epsilon/n\right)^{\frac{1}{p}}$ and $\Psi_{\mu}$ defined in \eqref{defPsimu}. Then, we see from \eqref{approxbd} and \eqref{lippsi} that
\begin{equation}\label{approxbd1}
0 \leq \Psi_{\mu}(\bm{x}) - \|\bm{x}\|_p^p \leq n\left(\frac{\mu}{2}\right)^p = \epsilon, \quad \forall\,\bm{x}\in\mathbb{R}^n,
\end{equation}
and $\Psi_{\mu}$ is globally Lipschitz continuous with Lipschitz constant $L_{\mu}:=\sqrt{n}p\mu^{p-1}$. Now, let $\lambda_{\epsilon}^*:=cL_{\mu}$, where $c>0$ is chosen as in Lemma \ref{errorbdlem}. For any $\bm{x}\in\mathbb{R}^n$, we also use $\mathcal{P}_{{\rm FEA}(A,\bm{b},\sigma,1)}(\bm{x})$ to denote the projection of $\bm{x}$ on ${\rm FEA}(A,\bm{b},\sigma,1)$. Then, for $\lambda \geq \lambda_{\epsilon}^*$ and any $\bm{x}\in\mathbb{R}^n$,
\begin{equation*}
\begin{aligned}
&~~F_{\lambda}(\bm{x})
=\|\bm{x}\|^p_p + \lambda (\|A\bm{x} - b\|_1 - \sigma)_+
\geq \Psi_{\mu}(\bm{x}) - \epsilon + \lambda (\|A\bm{x} - b\|_1 - \sigma)_+  \\
&\geq \Psi_{\mu}(\bm{x}) - \epsilon + \frac{\lambda}{c}\,\mathrm{dist}\left(\bm{x}, \,{\rm FEA}(A,\bm{b},\sigma,1)\right)
\geq \Psi_{\mu}(\bm{x}) + L_{\mu}\,\|\bm{x}-\mathcal{P}_{{\rm FEA}(A,\bm{b},\sigma,1)}(\bm{x})\| - \epsilon \\
&\geq \Psi_{\mu}(\mathcal{P}_{{\rm FEA}(A,\bm{b},\sigma,1)}(\bm{x})) - \epsilon
\geq \|\mathcal{P}_{{\rm FEA}(A,\bm{b},\sigma,1)}(\bm{x})\|_p^p - \epsilon
\geq \|\bm{x}^*\|_p^p - \epsilon
= F_{\lambda}(\bm{x}^*) - \epsilon,
\end{aligned}
\end{equation*}
where the first inequality follows from \eqref{approxbd1}, the second inequality follows from Lemma \ref{errorbdlem}, the third inequality follows from $\lambda\geq\lambda_{\epsilon}^*=cL_{\mu}$, the fourth inequality follows the Lipschitz continuity of $\Psi_{\mu}$ with Lipschtiz constant $L_{\mu}$, and the last two inequalities follows from \eqref{approxbd1} and the definition of $\bm{x}^*$ as a minimizer of problem \eqref{lpl1model}. This shows that $\bm{x}^*$ is an $\epsilon$-minimizer of problem \eqref{penprob} and completes the proof.
\endproof

From Theorems \ref{localthm} and \ref{epsithm1}, we see that if $\bm{x}^*$ is a local minimizer \textit{or} global minimizer of problem \eqref{lpl1model}, then it is also a local minimizer \textit{or} $\epsilon$-minimizer of problem \eqref{penprob}. Conversely, it is easy to see that if $\bm{x}^*$ is a local minimizer \textit{or} $\epsilon$-minimizer of problem \eqref{penprob} for some $\lambda>0$ and $\bm{x}^*\in{\rm FEA}(A,\bm{b},\sigma,1)$, then it is also a local minimizer \textit{or} $\epsilon$-minimizer of problem \eqref{lpl1model}. Finally, we shall study the case when $\bm{x}^*$ is a global minimizer of problem \eqref{penprob} for some $\lambda>0$ but $\bm{x}^*\notin{\rm FEA}(A,\bm{b},\sigma,1)$.

\begin{theorem}\label{epsithm2}
Suppose that $\tilde{\bm{x}}$ is an arbitrary feasible point of problem \eqref{lpl1model}, i.e., $\tilde{\bm{x}}\in{\rm FEA}(A,\bm{b},\sigma,1)$. Take any $\epsilon>0$ and consider any $\lambda \geq c\left(n^{\frac{p}{2}-1}\epsilon\right)^{-\frac{1}{p}}\|\tilde{\bm{x}}\|_p^p$, where $c>0$ is chosen as in Lemma \ref{errorbdlem}. Then, for any global minimizer $\bm{x}_{\lambda}^*$ of problem \eqref{penprob}, the projection $\mathcal{P}_{{\rm FEA}(A,\bm{b},\sigma,1)}(\bm{x}_{\lambda}^*)$ is an $\epsilon$-minimizer of problem \eqref{lpl1model}.
\end{theorem}
\beginproof
First, from the definition of $F_{\lambda}$ and the global optimality of $\bm{x}_{\lambda}^*$, we have
\begin{eqnarray}
&&\|\bm{x}_{\lambda}^*\|_p^p \leq F_{\lambda}(\bm{x}_{\lambda}^*) \leq F_{\lambda}(\bm{x}) = \|\bm{x}\|_p^p, \quad \forall\,\bm{x} \in {\rm FEA}(A,\bm{b},\sigma,1), \label{addbd} \\ [4pt]
&&(\|A\bm{x}_{\lambda}^* - b\|_1 - \sigma)_+ \leq \lambda^{-1}F_{\lambda}(\bm{x}_{\lambda}^*) \leq \lambda^{-1}F_{\lambda}(\tilde{\bm{x}}) = \lambda^{-1}\|\tilde{\bm{x}}\|_p^p.  \label{bdbyfea}
\end{eqnarray}
Then, for any $\bm{x} \in {\rm FEA}(A,\bm{b},\sigma,1)$, we have
\begin{equation*}
\begin{aligned}
&~~\|\mathcal{P}_{{\rm FEA}(A,\bm{b},\sigma,1)}(\bm{x}_{\lambda}^*)\|_p^p - \|\bm{x}\|_p^p
\leq \|\mathcal{P}_{{\rm FEA}(A,\bm{b},\sigma,1)}(\bm{x}_{\lambda}^*)\|_p^p - \|\bm{x}_{\lambda}^*\|_p^p  \\
&\leq \|\mathcal{P}_{{\rm FEA}(A,\bm{b},\sigma,1)}(\bm{x}_{\lambda}^*) - \bm{x}_{\lambda}^*\|_p^p
= n \cdot {\textstyle\frac{1}{n}\sum^n_{i=1}}\left(\left|[\mathcal{P}_{{\rm FEA}(A,\bm{b},\sigma,1)}(\bm{x}_{\lambda}^*)]_i - [\bm{x}_{\lambda}^*]_i\right|^2\right)^{\frac{p}{2}} \\
&\leq n \left({\textstyle\frac{1}{n}\sum^n_{i=1}}\left|[\mathcal{P}_{{\rm FEA}(A,\bm{b},\sigma,1)}(\bm{x}_{\lambda}^*)]_i - [\bm{x}_{\lambda}^*]_i\right|^2\right)^{\frac{p}{2}}
= n^{1-\frac{p}{2}} \|\mathcal{P}_{{\rm FEA}(A,\bm{b},\sigma,1)}(\bm{x}_{\lambda}^*) - \bm{x}_{\lambda}^*\|^p \\
& = n^{1-\frac{p}{2}} \big{[}\mathrm{dist}\left(\bm{x}_{\lambda}^*, \,{\rm FEA}(A,\bm{b},\sigma,1)\right)\big{]}^p
\leq n^{1-\frac{p}{2}}\big{[}c\,(\|A\bm{x}_{\lambda}^* - b\|_1 - \sigma)_+\big{]}^p \\
&\leq n^{1-\frac{p}{2}}\left[c\lambda^{-1}\|\tilde{\bm{x}}\|_p^p\right]^p
\leq \epsilon,
\end{aligned}
\end{equation*}
where the first inequality follows from \eqref{addbd}, the second inequality follows from \cite[Lemma 2.4]{clp2016penalty}, the third inequality follows from the concavity of the function $t \mapsto t^{\frac{p}{2}}$ for nonnegative $t$, the fourth inequality follows from Lemma \ref{errorbdlem} and the last two inequality follows from \eqref{bdbyfea} and the choice of $\lambda$. This implies that $\mathcal{P}_{{\rm FEA}(A,\bm{b},\sigma,1)}(\bm{x}_{\lambda}^*)$ is an $\epsilon$-minimizer of \eqref{lpl1model} and completes the proof.
\endproof

\end{APPENDICES}

% Acknowledgments here
\section*{Acknowledgments}

The authors are grateful to the editor and the anonymous referees for their valuable suggestions and comments, which have helped to improve the quality of this paper. The authors would also like to thank the CAS AMSS-PolyU Joint Laboratory of Applied Mathematics for its support while this research was being conducted. The research of Shuhuang Xiang was supported in part by the National Natural Science Foundation of China (Grant No. 11771454).

%%%%%%%%%%%%%%%%%%%%%%%%%%%%%%%%%%%
%%% References
\bibliographystyle{informs2014}
\bibliography{Ref_lplq}

\end{document}